\documentclass{amsart}
\usepackage{amssymb, latexsym}

\def\bi{\begin{itemize}}
\def\bs{\begin{split}}
\def\es{\end{split}}
\def\ba{\begin{align}}
\def\bas{\begin{align*}}
\def\ea{\end{align}}
\def\eas{\end{align*}}
\def\Im{{\hbox{Im}}}

\def\Re{{\hbox{Re}}}

\def\C{{\mathbb C}} 
\def\R{{{\mathbb R}}}
\def\Z{{{\mathbb Z}}}

\def\O{{{\mathcal O}}}
\def\lN{{\mathcal N}}
\def\p{{\phi}}
\def\tchi{{\tilde{\chi}}}
\def\mass{{m}}   
\def\momentum{{p}} 
\def\eps{\varepsilon}
\def \endprf{\hfill $\square$}
\def\emph#1{{\it #1}}
\def\textbf#1{{\bf #1}}
\def\divider#1{$\bullet\quad${\bf #1}}

\newcommand{\half}{\frac{1}{2}}

\newcommand{\nabb}{\mbox{$\nabla \mkern-13mu /$\,}}

\theoremstyle{plain}
\newtheorem{theorem}{Theorem}
\newtheorem{definition}[theorem]{Definition}
\newtheorem{remark}[theorem]{Remark}
\newtheorem{proposition}[theorem]{Proposition}
\newtheorem{lemma}[theorem]{Lemma}
\newtheorem{corollary}[theorem]{Corollary}

\numberwithin{equation}{section}
\numberwithin{theorem}{section}

\begin{document}

\title[Scattering for 3D critical NLS]{Global well-posedness and scattering for the energy-critical nonlinear Schr\"odinger equation in $\R^3$}
\author{J.~Colliander}
\thanks{J.C. is supported in part by N.S.F. Grant DMS 0100595,
  N.S.E.R.C. Grant R.G.P.I.N. 250233-03 and the Sloan Foundation.}
\address{University of Toronto and IAS}
\author{M.~Keel}
\thanks{M.K. was supported in part by N.S.F. Grant DMS-0303704; and by the
McKnight and Sloan Foundations.}
\address{University of Minnesota}
\author{G.~Staffilani}
\thanks{G.S. is supported in part by N.S.F. Grant DMS-0100375,
N.S.F. Grant DMS-0111298 through the IAS, and the Sloan Foundation.}
\address{MIT and IAS}
\author{H.~Takaoka}
\address{Kobe University}
\thanks{H.T. is supported in part by J.S.P.S. Grant No. 15740090 and by a J.S.P.S.      
 Postdoctoral Fellowship for Research Abroad.}
\author{T.~Tao}
\thanks{T.T. is a Clay Prize Fellow and is supported in part by grants
from the Packard Foundation.}
\address{University of California, Los Angeles and Australian National University}


\vspace{-0.3in}
\begin{abstract}
We obtain global well-posedness, scattering, and global $L^{10}_{t,x}$
spacetime bounds for energy-class solutions to the quintic defocusing
Schr\"odinger equation in $\R^{1+3}$, which is energy-critical. In particular,
this establishes global existence of classical solutions.  Our work extends
the results of Bourgain \cite{borg:scatter} and
Grillakis  \cite{grillakis:scatter}, which handled the radial case.
The method is similar in spirit to the induction-on-energy strategy
of Bourgain \cite{borg:scatter}, but we perform the induction analysis
in both frequency space and physical space simultaneously, and replace
 the Morawetz inequality by an interaction variant (first used in
\cite{ckstt:french},  \cite{ckstt:cubic-scatter}).  The principal
advantage of the  interaction Morawetz estimate is that it is not
localized to the
spatial origin and so is better able to handle nonradial solutions.
In particular, this interaction estimate, together with an almost-conservation
argument controlling the movement of $L^2$ mass in frequency space, rules out the
possibility of energy concentration.
\end{abstract}

\date{28 January 2004.}

\maketitle

\tableofcontents

\section{Introduction}\label{introduction-sec}

\subsection{Critical NLS and main result}

We consider the Cauchy problem for the quintic defocusing
Schr\"odinger equation in $\R^{1+3}$
\begin{equation}\label{nls}
\left\{
\begin{matrix}
iu_t + \Delta u = |u|^4 u \\
u(0,x) = u_0(x).
\end{matrix}
\right.
\end{equation}
where $u(t,x)$ is a complex-valued field in spacetime $\R_t \times
\R_x^3$. This equation has as Hamiltonian,
\begin{equation}\label{hamil-def}
E(u(t)) := \int \frac{1}{2} |\nabla u(t,x)|^2 + \frac{1}{6}
|u(t,x)|^6\ dx.
\end{equation}
Since the Hamiltonian \eqref{hamil-def} is preserved by the flow
\eqref{nls} we shall often refer to it as the {\em energy} and
write $E(u)$ for $E(u(t))$.


Semilinear  Schr\"odinger equations - with and without potentials, and
with various nonlinearities -   arise as
models for diverse physical phenomena, including Bose-Einstein
condensates \cite{bereference1, bereference2} and as a description
of the envelope dynamics of a general dispersive wave in
a weakly nonlinear medium (see e.g. the survey in  \cite{sulemsulem},
Chapter 1).  Our interest here in the defocusing quintic equation \eqref{nls}
is motivated mainly though by the fact that the problem is
critical with respect to the energy norm.  Specifically,   we map
a solution to another solution through the scaling  $u \mapsto u^\lambda$ defined by
\begin{equation}\label{scaling}
u^\lambda(t,x) := \frac{1}{\lambda^{1/2}} u(\frac{t}{\lambda^2},
\frac{x}{\lambda}),
\end{equation}
and this scaling leaves both terms in the
energy invariant.


The Cauchy problem for this equation has been intensively studied
(\cite{cwI}, \cite{grillakis:scatter}, \cite{borg:scatter},
\cite{borg:book},\cite{gv:localreference}, \cite{kato}). It is
known (see e.g. \cite{cw0,cwI}) that if the initial data $u_0(x)$
has finite energy, then the Cauchy problem is locally well-posed,
in the sense that there exists a local-in-time solution to
\eqref{nls} which lies in $C^0_t \dot H^1_x \cap L^{10}_{t,x}$,
and is unique in this class; furthermore the map from initial data
to solution is locally Lipschitz continuous in these norms.  If
the energy is small, then the solution is known to exist globally
in time, and scatters to a solution $u_\pm(t)$ to the free
Schr\"odinger equation $(i\partial_t + \Delta) u_\pm = 0$, in the
sense that $\| u(t) - u_\pm(t) \|_{\dot H^1(\R^3)} \to 0$ as $t
\to \pm \infty$. For \eqref{nls} with large initial data,
the arguments in \cite{cw0, cwI} do not
extend to yield global well-posedness, even
with the conservation of the energy \eqref{hamil-def}, because the
time of existence given by the local theory depends on the profile
of the data as well as on the energy\footnote{This is in constrast with
sub-critical equations such as the cubic equation $iu_t + \Delta u
= |u|^2 u$, for which one can use the local well-posedness theory
to yield global well-posedness and scattering even for large
energy data (see \cite{gv:scatter}, and the surveys
\cite{cazbook}, \cite{cazbooknew}).}.
For large finite energy data which is assumed to be in addition radially
symmetric, Bourgain \cite{borg:scatter} proved global existence
and scattering for \eqref{nls} in ${\dot{H}}^1(\R^3)$.
Subsequently Grillakis
\cite{grillakis:scatter} gave a different argument which recovered part of
\cite{borg:scatter} - namely, global existence from smooth, radial, finite energy data.
For general large data - in particular, general smooth data - global existence
and scattering were open.

Our main result is the following global well-posedness result for
\eqref{nls} in the energy class.

\begin{theorem}\label{main}  For any $u_0$ with finite energy,
$E(u_0) < \infty$, there exists a unique\footnote{In fact, uniqueness actually holds in the larger
space $C^0_t(\dot H^1_x)$ (thus eliminating the constraint that $u \in L^{10}_{t,x}$), as one can show by adapting the arguments of \cite{katounique}, \cite{twounique},
\cite{FPT_NLSunique}; see Section \ref{remarks-sec}.}
global solution $u \in
C^0_t(\dot H^1_x) \cap L^{10}_{t,x}$ to \eqref{nls} such that
\begin{equation}\label{l-10}
\int_{-\infty}^\infty \int_{\R^3}|u(t,x)|^{10}\ dx dt \leq C(E(u_0)).
\end{equation}
for some constant $C(E(u_0))$ that depends only on the energy.
\end{theorem}

As is well-known (see e.g. \cite{borg:book}, or \cite{ckstt:cubic-scatter} for
the sub-critical analogue), the $L^{10}_{t,x}$
bound above also gives scattering, asymptotic completeness,
and uniform regularity:

\begin{corollary}\label{scat-reg}  Let $u_0$ have finite energy.
Then there exists finite energy solutions $u_\pm(t,x)$ to the free
Schr\"odinger equation $(i\partial_t + \Delta) u_\pm = 0$ such that
$$ \| u_\pm(t) - u(t) \|_{\dot H^1} \to 0 \hbox{ as } t \to \pm
\infty.$$
Furthermore, the maps $u_0 \mapsto u_\pm(0)$ are
homeomorphisms from $\dot H^1(\R^3)$
to $\dot H^1(\R^3)$.
Finally, if $u_0 \in H^s$ for some $s > 1$, then $u(t) \in H^s$
for all time $t$, and one has the uniform bounds
$$ \sup_{t \in \R} \| u(t) \|_{H^s} \leq C( E(u_0), s ) \| u_0
\|_{H^s}.$$
\end{corollary}

It is also fairly standard to show that the $L^{10}_{t,x}$ bound \eqref{l-10} implies
further spacetime integrability on $u$, for instance $u$ obeys all the
Strichartz estimates that a free solution with the same regularity does
(see, for example, Lemma \ref{persistence} below).

The results here have analogs in previous work on second order wave equations
on $\R^{3+1}$ with energy-critical (quintic) defocusing nonlinearities.
Global in time existence for such equations from smooth data
was shown by Grillakis \cite{g_waveI},\cite{g_waveII} (for radial
data see Struwe \cite{struwe}, for small energy data see Rauch
\cite{rauch}); global in time solutions from finite energy data was
shown in Kapitanski \cite{kapitanski}, Shatah-Struwe\cite{shatah_struwe}.
For an analog of the scattering statement in Corollary \ref{scat-reg}
for the critical wave equation see Bahouri-Shatah \cite{bs}, Bahouri-G\'erard \cite{bg}; for the scattering statement for Klein-Gordon
equations see Nakanishi \cite{nak:imrn} (for radial data, see Ginibre-Soffer-Velo\cite{gsv}).  The existence
results mentioned here all involve an argument showing
that the solution's energy can't concentrate.  These energy nonconcentration
proofs combine {\em Morawetz inequalities} (\emph{a priori} estimates
for the nonlinear equation which  bound some quantity that scales
like energy) with careful analysis that strengthens the Morawetz
bound to control of energy. Besides the presence of
infinite propagation speeds, a main difference between
\eqref{nls} and the hyperbolic analogs is that here time
scales like $\lambda^2$, and as a consequence the quantity bounded
by the Morawetz estimate is supercritical with respect to energy.

Section \ref{overview-sec} below provides a fairly complete outline of the
proof of Theorem \ref{main}.  In this introduction we only briefly sketch
some of the ideas involved:   a suitable
modification of the Morawetz inequality for \eqref{nls},
along with the frequency localized $L^2$ almost-conservation
law that we'll ultimately use to prohibit energy concentration.

A typical example of a Morawetz inequality for \eqref{nls} is the following
bound due to Lin-Strauss \cite{linstrauss} who cite \cite{morawetz}
as motivation,
\begin{equation}\label{morawetz}
\int_I \int_{\R^3} \frac{|u(t,x)|^6}{|x|}\ dx dt \lesssim (\sup_{t
\in I} \| u(t) \|_{\dot H^{1/2}})^2
\end{equation}
for arbitrary time intervals $I$. (The estimate \eqref{morawetz} follows
from a computation showing the quantity,
\begin{equation}
\label{guts}
 \int_{\R^3} \Im(\bar{u} \nabla u \cdot \frac{x}{|x|}) dx
 \end{equation}
is monotone in time.)  Observe that the right-hand side of \eqref{morawetz}
will not grow in $I$ if the $H^1$ and $L^2$ norms are bounded, and
so this estimate gives a uniform bound on the left hand side where
$I$ is any interval on which we know the solution exists. However,
in the energy-critical problem \eqref{nls} there are two drawbacks
with this estimate.  The first is that the right-hand side
involves the $\dot H^{1/2}$ norm, instead of the energy $E$.
This is troublesome
since any Sobolev norm rougher than $\dot H^1$ is supercritical
with respect to the scaling \eqref{scaling}.  Specifically, the
right hand side of \eqref{morawetz} increases without bound when
we simply scale given finite energy initial data according to
\eqref{scaling} with $\lambda$ large.  The second difficulty is
that the left-hand side is localized near the spatial origin $x
= 0$ and does not convey as much information about the solution
$u$ away from this origin.  To get around the first difficulty
Bourgain \cite{borg:scatter} and Grillakis
\cite{grillakis:scatter} introduced a localized variant of the
above estimate:
\begin{equation}\label{morawetz-loc}
\int_I \int_{|x| \lesssim |I|^{1/2}} \frac{|u(t,x)|^6}{|x|}\ dx dt
\lesssim E(u) |I|^{1/2}.
\end{equation}
As an example of the usefulness of \eqref{morawetz-loc}, we observe that this estimate prohibits
the existence of finite energy (stationary) pseudosoliton solutions to \eqref{nls}. By a
(stationary) {\it{pseudosoliton}} we mean a solution such that $|u(t,x)| \sim 1$ for all $t \in
\R$ and $|x| \lesssim 1$; this notion includes soliton and breather type solutions. Indeed,
applying  \eqref{morawetz-loc} to such a solution,  we would see that the left-hand side grows by
at least $|I|$, while the right-hand side is $O(|I|^{\half})$, and so a pseudosoliton solution
will lead to a contradiction for $|I|$ sufficiently large.
A similar argument allows one
to use \eqref{morawetz-loc} to prevent ``sufficiently rapid''
concentration  of (potential) energy at the origin; for instance,
\eqref{morawetz-loc} can also be used to rule out self-similar
type blowup\footnote{This is not the only type of self-similar
blowup scenario; another type is when the energy concentrates in a
ball $|x| \leq A|t-t_0|^{1/2}$ as $t \to t_0^-$. This type of
blowup is consistent with the scaling \eqref{scaling} and is not
directly ruled out by \eqref{morawetz-loc}, however it can instead
be ruled out by spatially local mass conservation estimates. See
\cite{borg:scatter}, \cite{grillakis:scatter}.}, where the
potential energy density $|u|^6$ concentrates in the ball $|x| <
A|t-t_0|$ as $t \to t_0^-$ for some fixed $A > 0$. In
\cite{borg:scatter}, one main use of \eqref{morawetz-loc} was to
show that for each fixed time interval $I$, there exists at least
one time $t_0 \in I$ for which the potential energy was dispersed
at scale $|I|^{1/2}$ or greater (i.e. the potential energy could
not concentrate on a ball $|x| \ll |I|^{1/2}$ for all times in
$I$).

To summarize, the localized Morawetz estimate
\eqref{morawetz-loc} is very good at preventing $u$ from concentrating
near the origin; this is especially useful in the case of radial
solutions $u$, since the radial symmetry
(combined with conservation of energy) enforces decay of $u$ away
from the origin, and so resolves the second difficulty with the
Morawetz estimate mentioned earlier.  However, the estimate is
less useful when the solution is allowed to concentrate away from
the origin.  For instance, if we aim to preclude the existence of
a moving pseudosoliton  solution, in which $|u(t,x)| \sim 1$ when $|x -
vt| \lesssim 1$ for some fixed velocity $v$, then the left-hand
side of \eqref{morawetz-loc} only grows like $\log |I|$ and so one
does not necessarily obtain a contradiction\footnote{At first
glance it may appear that the global estimate \eqref{morawetz} is
still able to preclude the existence of such a pseudosoliton,
since the right-hand side does not seem to grow much as $I$ gets
larger. This can be done in the cubic problem (see e.g.
\cite{gv:scatter}) but in the critical problem one can lose
control of the $\dot H^{1/2}$ norm, by adding some very low
frequency components to the soliton solution $u$.  One might
object that one could use $L^2$ conservation to control the
$H^{1/2}$ norm, however one can rescale the solution to make the
$L^2$ norm (and hence the $\dot H^{1/2}$ norm) arbitrarily
large.}.

It is thus of interest to remove the $1/|x|$ denominator in
\eqref{morawetz}, \eqref{morawetz-loc}, so that these estimates
can more easily prevent concentration at arbitrary locations in
spacetime.  In \cite{ckstt:french}, \cite{ckstt:cubic-scatter}
this was achieved by translating the origin in the integrand of \eqref{guts} to
an arbitrary point $y$, and averaging against the $L^1$ mass density $|u(y)|^2\ dy$.
In particular, the following {\emph{interaction Morawetz}}  estimate\footnote{Strictly speaking, in
\cite{ckstt:french}, \cite{ckstt:cubic-scatter} this estimate was
obtained for the \emph{cubic} defocusing nonlinear Schr\"odinger
equation instead of the quintic, but the argument in fact works
for all nonlinear Schr\"odinger equations with a pure power
defocusing nonlinearity, and even for a slightly more general
class of repulsive nonlinearities satisfying a standard
monotonicity condition.  See \cite{ckstt:cubic-scatter} and
Section \ref{slimflim-sec} below for more discussion.}
\begin{equation}\label{morawetz-interaction}
\int_I \int_{\R^3} |u(t,x)|^4\ dx dt \lesssim \| u(0)\|_{L^2}^2 (\sup_{t
\in I} \| u(t) \|_{\dot H^{1/2}})^2
\end{equation}
was obtained.  (We have since learned that this averaging argument has an
analog in early work presenting and analyzing interaction functionals for one dimensional hyperbolic systems, e.g. \cite{glimmreferenceI, glimmreferenceII}.)  This $L^4_{t,x}$ estimate already gives a short proof
of scattering in the energy class (and below!) for the cubic
nonlinear Schr\"odinger equation (see \cite{ckstt:french},
\cite{ckstt:cubic-scatter}); however, like \eqref{morawetz}, this
estimate is not suitable for the critical problem because the
right-hand side is not controlled by the energy $E(u)$.  One could
attempt to localize \eqref{morawetz-interaction} as in
\eqref{morawetz-loc}, obtaining for instance a scale-invariant estimate such as
\begin{equation}\label{stupid-morawetz}
\int_I \int_{|x| \lesssim |I|^{1/2}} |u(t,x)|^4\ dx dt \lesssim
E(u)^2 |I|^{3/2}
\end{equation}
but this estimate, while true (in fact it follows immediately from
Sobolev and H\"older), is useless for such purposes as prohibiting
soliton-like behaviour, since the left-hand side grows like $|I|$
while the right-hand side grows like $|I|^{3/2}$.  Nor is this estimate
useful for preventing any sort of energy concentration.


Our solution to these difficulties proceeds in the context of an induction-on-energy
argument as in \cite{borg:scatter}:
assume for contradiction that Theorem \ref{main} is false, and consider a solution of minimal energy among all those solutions with $L^{10}_{x,t}$ norm above some threshhold.
We first show, without relying on any of the above Morawetz-type inequalities, that
such a {\emph {minimal energy blowup solution}} would have to be localized in both
frequency and in space at all times. Second, we prove that this localized
blowup solution satisfies Proposition \ref{Morawetz}, which
localizes \eqref{morawetz-interaction} in \emph{frequency} rather
than in space.  Roughly speaking, the {\em{frequency localized Morawetz
inequality}} of Proposition \ref{Morawetz} states that after throwing away
some small energy, low frequency portions of the blow-up solution, the
remainder obeys good $L^4_{t,x}$ estimates.  In principle, this
estimate should follow simply by repeating the proof of
\eqref{morawetz-interaction} with $u$ replaced by the high frequency
portion of the solution, and then controlling error terms; however
some of the error terms are rather difficult and the proof of the
frequency-localized Morawetz inequality is quite technical.
We emphasize that, unlike
the estimates \eqref{morawetz}, \eqref{morawetz-loc},
\eqref{morawetz-interaction}, the frequency-localized Morawetz
inequality \eqref{pm} is not an \emph{a priori} estimate valid for all
solutions of \eqref{nls}, but
instead applies only to
minimal energy blowup solutions; see Section \ref{overview-sec} for further discussion and precise definitions.

The strategy is then to
try to use Sobolev embedding to boost this $L^4_{t,x}$ control to $L^{10}_{t,x}$ control
which would contradict the existence of the blow-up solution.
There is, however, a remaining enemy, which is that the solution may
shift its energy from low frequencies to high, possibly causing the
$L^{10}_{t,x}$ norm to blow up while the $L^4_{t,x}$ norm stays
bounded.  To prevent this we look at what such a frequency evacuation
would imply for the location - in frequency space -  of the blow-up
solution's $L^2$ mass.  Specifically, we
prove a frequency localized $L^2$ mass estimate that gives
us information for longer time intervals than seems to be available from the {\em {spatially}} localized mass conservation laws used in the previous radial work (\cite{borg:scatter, grillakis:scatter}).  By combining this frequency localized mass estimate with
the $L^4_{t,x}$ bound and plenty of Strichartz estimate analysis, we
can control the  movement of energy and mass from one frequency range
to another, and prevent the low-to-high cascade from occurring.  The argument here is motivated by our previous low-regularity work involving
almost conservation laws (e.g. \cite{ckstt:cubic-scatter}).

The remainder of the paper is organized as follows:  Section \ref{conserv-sec}
reviews some simple, classical conservation laws for Schr\"odinger equations which
will be used througout, but especially in proving the frequency localized interaction
Morawetz estimate.  In section \ref{strichartz-sec} we recall some linear and multilinear Strichartz estimates, along with the useful nonlinear perturbation
statement of Lemma \ref{perturb}.  Section \ref{overview-sec}
outlines in some detail the argument behind our main Theorem, leaving the proofs of each
step to sections \ref{decoupled-sec} - \ref{evacuate-sec} of the paper.  Section
\ref{remarks-sec}  presents some miscellaneous remarks, including a proof
of the unconditional uniqueness statement alluded to above.

{\bf Acknowledgements}:


We thank the Institute for Mathematics and its Applications (IMA)
for hosting our collaborative meeting in July 2002. We thank Andrew Hassell,
Sergiu Klainerman,  and
Jalal Shatah for interesting discussions related to the
interaction Morawetz estimate, and Jean Bourgain for valuable comments on
an early draft of this paper, to Monica Visan and the anonymous referee for their thorough reading
of the manuscript and for many important corrections, and to Changxing Miao and Guixiang Xu for further corrections.   We thank  Manoussos Grillakis
for explanatory details related to  \cite{grillakis:scatter}.
Finally, it will be clear to the reader that our work here relies heavily in places
on arguments developed by J. Bourgain in \cite{borg:scatter}.

\subsection{Notation}\label{notation-sec}

If $X,Y$ are nonnegative quantities, we use $X \lesssim Y$ or $X = O(Y)$ to denote the estimate $X \leq C Y$
for some $C$ (which may depend on the critical energy $E_{crit}$ (see Section \ref{overview-sec})
but not on any other parameter such as
$\eta$), and $X \sim Y$ to denote the estimate $X \lesssim Y \lesssim X$.
We use $X \ll Y$ to mean $X \leq c Y$ for some small constant
$c$ (which is again allowed to depend on $E_{crit}$).

We use $C \gg 1$ to denote various large finite constants, and $0 < c
\ll 1$ to denote various small constants.

The Fourier transform on $\R^3$ is defined by
$$ \hat f(\xi) := \int_{\R^3} e^{-2\pi i x \cdot \xi} f(x)\ dx,$$
giving rise to the fractional differentiation operators $|\nabla|^s$,
$\langle \nabla \rangle^s$ defined by
$$ \widehat{|\nabla|^s f}(\xi) := |\xi|^s \hat f(\xi); \quad
\widehat{\langle \nabla \rangle^s f}(\xi) := \langle \xi \rangle^s \hat
f(\xi)$$
where $\langle \xi \rangle := (1 + |\xi|^2)^{1/2}$. In particular, we will use $\nabla$ to denote the spatial gradient $\nabla_x$. 
This in turn defines the Sobolev norms
$$ \| f \|_{\dot H^s(\R^3)} := \| |\nabla|^s f \|_{L^2(\R^3)}; \quad
\| f \|_{H^s(\R^3)} := \| \langle \nabla \rangle^s f \|_{L^2(\R^3)}.$$
More generally we define
$$ \| f \|_{\dot W^{s,p}(\R^3)} := \| |\nabla|^s f \|_{L^p(\R^3)}; \quad
\| f \|_{W^{s,p}(\R^3)} := \| \langle \nabla \rangle^s f \|_{L^p(\R^3)}$$
for $s \in \R$ and $1 < p < \infty$.

We let $e^{it\Delta}$ be the free Schr\"odinger propagator; in terms of the Fourier transform,
this is given by
\begin{equation}\label{fourier-rep}
 \widehat{e^{it\Delta} f}(\xi) = e^{-4\pi^2 i t |\xi|^2} \hat f(\xi)
\end{equation}
while in physical space we have
\begin{equation}\label{fundamental-soln}
e^{it\Delta} f(x) = \frac{1}{(4 \pi i t)^{3/2}} \int_{\R^3} e^{i|x-y|^2/4t} f(y)\ dy
\end{equation}
for $t \neq 0$, using an appropriate branch cut to define the complex square root. 
In particular the propagator preserves all the Sobolev norms $H^s(\R^3)$ and $\dot H^s(\R^3)$, and
also obeys the \emph{dispersive inequality}
\begin{equation}\label{dispersive}
\| e^{it\Delta} f \|_{L^\infty_x(\R^3)} \lesssim |t|^{-3/2} \| f\|_{L^1_x(\R^3)}.
\end{equation}
We also record \emph{Duhamel's formula}
\begin{equation}\label{duhamel}
 u(t) = e^{i(t-t_0)\Delta} u(t_0) -i  \int_{t_0}^t e^{i(t-s)\Delta} (iu_t + \Delta u)(s)\ ds
\end{equation}
for any Schwartz $u$ and any times $t_0, t \in \R$, with the convention that $\int_{t_0}^t = - \int_t^{t_0}$ if $t < t_0$.

We use the notation $\O(X)$ to denote an expression which is \emph{schematically} of the form $X$; this means that $\O(X)$ is a finite linear combination of expressions which look like $X$ but with some factors possibly replaced by their complex conjugates, thus for instance $3 u^2 \overline{v}^2 |v|^2 + 9 |u|^2 |v|^4 + 3 \overline{u}^2 v^2 |v|^2$ qualifies to be of the form $\O(u^2 v^4)$, and similarly we have
\begin{equation}\label{6-schematic}
|u+v|^6 = |u|^6 + |v|^6 + \sum_{j=1}^5 \O(u^j v^{6-j})
\end{equation}
and
\begin{equation}\label{5-schematic}
|u+v|^4 (u+v) = |u|^4 u + |v|^4 v + \sum_{j=1}^4 \O(u^j v^{5-j}).
\end{equation}

We will sometimes denote partial derivatives using subscripts:
$\partial_{x_j} u = \partial_j u = u_j$. We will also implicitly use
the summation convention when indices are repeated in expressions below.

We shall need the following Littlewood-Paley projection operators.
Let $\varphi(\xi)$ be a bump function adapted to the ball $\{ \xi \in
\R^3: |\xi| \leq 2 \}$ which equals 1 on the ball $\{ \xi \in \R^3:
|\xi| \leq 1 \}$.  Define a \emph{dyadic number} to be any number $N \in 2^\Z$ of the form
$N = 2^j$ where $j \in \Z$ is an integer.  For each dyadic number $N$,
we define the Fourier multipliers
\begin{align*}
\widehat{P_{\leq N} f}(\xi) &:= \varphi(\xi/N) \hat f(\xi)\\
\widehat{P_{> N} f}(\xi) &:= (1 - \varphi(\xi/N)) \hat f(\xi)\\
\widehat{P_N f}(\xi) &:= (\varphi(\xi/N) - \varphi(2\xi/N)) \hat f(\xi).
\end{align*}
We similarly define $P_{<N}$ and $P_{\geq N}$.
Note in particular the telescoping identities
$$
P_{\leq N} f = \sum_{M \leq N} P_M f; \quad
P_{> N} f = \sum_{M > N} P_M f; \quad
f = \sum_M P_M f
$$
for all Schwartz $f$, where $M$ ranges over dyadic numbers.  We also
define
$$ P_{M < \cdot \leq N} := P_{\leq N} - P_{\leq M} = \sum_{M < N' \leq N} P_{N'}$$
whenever $M \leq N$ are dyadic numbers.  Similarly define $P_{M \leq \cdot \leq N}$, etc.

The symbol $u$ shall always refer to a solution to the nonlinear Schr\"odinger equation
\eqref{nls}.  We shall use $u_N$ to denote the frequency piece $u_N := P_N u$ of $u$, and
similarly define $u_{\geq N} = P_{\geq N} u$, etc.  While this may cause some confusion with
the notation $u_j$ used to denote derivatives of $u$, the meaning of the subscript should be
clear from context.

The Littlewood-Paley operators commute with derivative operators
(including $|\nabla|^s$ and $i\partial_t + \Delta$), the propagator $e^{it\Delta}$, and
conjugation operations,
are self-adjoint, and are bounded on every Lebesgue space $L^p$
and Sobolev space $\dot H^s$ (if $1 \leq p \leq \infty$, of course).
Furthermore, they obey the following easily verified Sobolev (and
Bernstein) estimates for $\R^3$ with $s \geq 0$ and $1 \leq p \leq q
\leq \infty$:
\begin{align}
\| P_{\geq N} f \|_{L^p} &\lesssim N^{-s} \| |\nabla|^s P_{\geq N} f
\|_{L^p} \label{nabla-highfreq}\\
\| P_{\leq N} |\nabla|^s f \|_{L^p} &\lesssim N^{s} \| P_{\leq N} f
\|_{L^p} \label{nabla-lowfreq}\\
\| P_N |\nabla|^{\pm s} f \|_{L^p} &\sim N^{\pm s} \| P_N f \|_{L^p}
\label{nabla-block}\\
\| P_{\leq N} f \|_{L^q} &\lesssim N^{\frac{3}{p} - \frac{3}{q}} \|
P_{\leq N} f \|_{L^p} \label{sobolev}\\
\| P_N f \|_{L^q} &\lesssim N^{\frac{3}{p} - \frac{3}{q}} \| P_N f
\|_{L^p}.\label{bernstein}
\end{align}

\section{Local conservation laws}\label{conserv-sec}


In this section we record some standard facts about the (non)conservation of
mass, momentum and energy densities for general nonlinear Schr\"odinger equations of
the form\footnote{We will use $\p$ to denote general solutions to Schr\"odinger-type equations,
reserving the symbol $u$ for solutions to the quintic defocusing nonlinear Schr\"odinger equation \eqref{nls}.}
\begin{equation}
  \label{forced}
 i \partial_t \p + \Delta \phi = \mathcal{N}
\end{equation}
on the spacetime slab $I_0 \times \R^d$ with $I_0$ a compact
interval. Our primary interest is of course the quintic defocusing
case \eqref{nls} on $I_0 \times \R^3$ when $\lN = |\p|^4 \p$, but we will also
discuss here the $U(1)$-gauge invariant Hamiltonian
case, when $\lN = F' (|\p|^2) \p$ with $\R$-valued $F$.  Later on
we will consider various truncated versions of \eqref{nls} with
non-Hamiltonian forcing terms.  These local conservation laws
 will be used not only to imply the usual global conservation of mass
and energy, but also derive ``almost conservation'' laws for various
localized portions of mass, energy, and momentum, where the localization
is either in physical space or frequency space.  The localized momentum inequalities
are closely related to virial identities, and will be used later to deduce an
interaction Morawetz inequality which is crucial to our argument.

To avoid technicalities (and to justify all exchanges of derivatives and integrals), let us work purely with fields 
$\phi$, $\mathcal{N}$ which are smooth, with all derivatives rapidly decreasing in space; in practice, we can
then extend the formulae obtained here to more general situations by limiting arguments.  We begin by introducing
some notation which will be used to describe the mass and momentum
(non)conservation properties of \eqref{forced}.

\begin{definition}
  \label{emtensor}
Given a (Schwartz) solution $\p$ of \eqref{forced} we define the mass
density
\begin{equation*}
  T_{00}(t,x)  := |\p (t,x)|^2,
\end{equation*}
the momentum density
\begin{equation*}
  T_{0j}(t,x) := T_{j0} (t,x) := 2 \Im ({\overline{\phi}} \phi_j ),
\end{equation*}
and the (linear part of the) momentum current
\begin{equation*}
  L_{jk}(t,x) = L_{kj} (t,x) := - \partial_j \partial_k |\phi (t,x)|^2 + 4 \Re
  (\overline{\phi_j} \p_k).
\end{equation*}
\end{definition}

\begin{definition}
  \label{mmbrackets}
Given any two (Schwartz) functions $f,g: \R^d \rightarrow \C$, we
define the mass bracket
\begin{equation}
  \label{mass-flux-def}
  \{ f , g \}_m := \Im (f \overline{g} )
\end{equation}
and the momentum bracket
\begin{equation}
  \label{momentum-flux-def}
  \{f , g \}_p := \Re (f \nabla \overline{g} - g \nabla \overline{f}).
\end{equation}
Thus $\{ f , g \}_m$ is a scalar valued function, while $\{f , g \}_p
$ defines a vector field on $\R^d$. We will denote the $j$th component
of $\{f , g \}_p$ by $\{f , g \}_p^j$.
\end{definition}

With these notions we can now express the mass and momentum
(non)conservation laws for \eqref{forced}, which can be validated with
straightforward computations.

\begin{lemma}[Local conservation of mass and momentum]
\label{local-conserv}
If $\p$ is a (Schwartz) solution to \eqref{forced} then we have the
local mass conservation identity
\begin{equation}
  \label{local-mass-conserv}
  \partial_t T_{00} + \partial_j T_{0j} = 2 \{ \lN , \p \}_m
\end{equation}
and the local momentum conservation identity
\begin{equation}
  \label{local-momentum-conserv}
  \partial_t T_{0j} + \partial_k L_{kj} = 2 \{ \lN , \p \}_p^j.
\end{equation}
Here we adopt the usual\footnote{Repeated Euclidean coordinate indices are summed.
As the metric is Euclidean, we will not systematically match subscripts and superscripts.} summation conventions for the indices $j,k$.
\end{lemma}

Observe that the mass current coincides with
the momentum density in \eqref{local-momentum-conserv}, while the
 momentum current in \eqref{local-momentum-conserv}
has some ``positive  definite'' tendencies (think of $\Delta =
\partial_k \partial_k$ as a
negative definite operator, whereas the $\partial_j$ will eventually
be dealt with by integration by parts, reversing the sign).  These
two facts will underpin the interaction Morawetz estimate obtained
in Section \ref{slimflim-sec}.

We now specialize to the gauge invariant Hamiltonian case, when
$\lN = F' (|\p|^2) \p$; note that \eqref{nls} would correspond to
the case $F(|\p|^2) = \frac{1}{3} |\p|^6$.  Observe that
\begin{equation}
  \label{mass-cancel}
  \{ F' (|\p|^2) \p , \p \}_m = 0
\end{equation}
and
\begin{equation}
  \label{momentum-cancel-general}
  \{ F'(|\p|^2 ) \p , \p\}_{p} = -\nabla G( |\p|^2 )
\end{equation}
where $G(z) := z F' (z) - F(z).$ In particular, for the quintic case
\eqref{nls} we have
\begin{equation}
  \label{momentum-cancel}
  \{|\p|^4 \p , \p \}_p = -\frac{2}{3} \nabla {|\p|^6}.
\end{equation}
Thus, in the gauge invariant case we can reexpress \eqref{local-momentum-conserv} as
\begin{equation}
  \label{momentum-nl-conserv}
  \partial_t T_{0j} + \partial_k T_{jk} = 0
\end{equation}
where
\begin{equation}
  \label{Tjk}
  T_{jk} := L_{jk} + 2\delta_{jk} G(|\p|^2)
\end{equation}
is the (linear and nonlinear) {\it{momentum current}}.
Integrating \eqref{local-mass-conserv} and \eqref{momentum-nl-conserv}
in space we see that the total mass
$$ \int_{\R^d} T_{00} \ dx = \int_{\R^d} |\p(t,x)|^2\ dx$$
and the total momentum
$$ \int_{\R^d} T_{0j}\ dx = 2 \int_{\R^d} \Im( \overline{\p(t,x)} \partial_j \p(t,x) )\ dx$$
are both conserved quantities.  In this Hamiltonian setting one can also verify the local energy conservation law
\begin{equation}\label{energy-conserv}
 \partial_t \left[ \frac{1}{2} |\nabla \p|^2 + \frac{1}{2} F(|\p|^2) \right]
+ \partial_j \left[ \Im ( \overline{\p}_k \p_{kj}) -
   F'(|\p|^2) \Im ( \p \overline{\p}_j ) \right] = 0
\end{equation}
which implies conservation of total energy
$$ \int_{\R^d} \frac{1}{2} |\nabla \p|^2 + \frac{1}{2} F(|\p|^2)\ dx.$$
Note also that \eqref{Tjk} continues the tendency of the
right-hand side of \eqref{local-momentum-conserv}
to be ``positive definite''; this is a manifestation of the defocusing
nature of the equation. Later in our argument, however, we will be forced to deal with
frequency-localized versions of the nonlinear Schr\"odinger
equations, in which one does not have perfect conservation of mass and
momentum, leading to a number of unpleasant error terms in our analysis.

\section{Review of Strichartz theory in $\R^{1+3}$}\label{strichartz-sec}

In this section we review some standard (and some slightly less standard)
Strichartz estimates in three dimensions, and their application to the well-posedness
and regularity theory for \eqref{nls}.  We use $L^q_t L^r_x$ to denote the spacetime norm
$$ \| u \|_{L^q_t L^r_x(\R \times \R^3)} := (\int_\R (\int_{\R^3} |u(t,x)|^r\ dx)^{q/r}\ dt)^{1/q},$$
with the usual modifications when $q$ or $r$ is equal to infinity, or when the domain $\R \times \R^3$
is replaced by a smaller region of spacetime such as $I \times \R^3$.  When $q=r$ we abbreviate
$L^q_t L^q_x$ as $L^q_{t,x}$.

\subsection{Linear Strichartz estimates}


We say that a pair $(q,r)$ of exponents is \emph{admissible} if
$\frac{2}{q} + \frac{3}{r} = \frac{3}{2}$ and $2 \leq q, r \leq
\infty$; examples include $(q,r) = (\infty, 2)$, $(10,30/13)$,
$(5, 30/11)$, $(4,3)$, $(10/3,10/3)$, and
$(2, 6)$.

Let $I \times \R^3$ be a spacetime slab.  We define\footnote{The presence of
the Littlewood-Paley projections here may seem unusual, but they are necessary
in order to obtain a key $L^4_t L^\infty_x$ endpoint Strichartz estimate below.}
the \emph{$L^2$
Strichartz norm} $\dot S^0(I \times \R^3)$ by
\begin{equation}\label{strichartz-def}
\| u \|_{\dot S^0(I \times \R^3)} := \sup_{(q,r)\hbox{ admissible}}
(\sum_N \| P_N u \|_{L^q_t L^r_x(I \times \R^3)}^2)^{1/2}
\end{equation}
and for $k=1,2$ we then define the \emph{$\dot H^k$ Strichartz norm}
$\dot S^k(I \times \R^3)$ by
$$
\| u \|_{\dot S^k(I \times \R^3)} := \| \nabla^k u \|_{\dot S^0(I
\times \R^3)}.$$
We shall work primarily with the $\dot H^1$ Strichartz norm, but will need
the $L^2$ and $\dot H^2$ norms to control high frequency and low frequency portions
of the solution $u$ respectively.

We observe the elementary inequality
\begin{equation}\label{squarefunction}
\| (\sum_N |f_N|^2)^{1/2} \|_{L^q_t L^r_x(I \times \R^3)} \leq
(\sum_N \| f_N \|_{L^q_t L^r_x(I \times \R^3)}^2)^{1/2}
\end{equation}
for all $2 \leq q,r \leq \infty$ and arbitrary functions $f_N$; this is easy to verify in the extreme cases
$(q,r) = (2,2), (2,\infty), (\infty,2), (\infty,\infty)$, and the intermediate cases then follow
by complex interpolation.  In particular, \eqref{squarefunction} holds for all admissible
exponents $(q,r)$.  From this and the Littlewood-Paley inequality (see e.g. \cite{stein:small})
we have
\begin{align*}
 \| u \|_{L^q_t L^r_x(I \times \R^3)}
&\lesssim \| (\sum_N |P_N u|^2)^{1/2} \|_{L^q_t L^r_x(I \times \R^3)} \\
&\lesssim (\sum_N \| P_N u \|_{L^q_t L^r_x(I \times \R^3)}^2)^{1/2}\\
&\lesssim \| u \|_{\dot S^0(I \times \R^3)}
\end{align*}
and hence
\begin{equation}\label{sr}
 \| \nabla u \|_{L^q_t L^r_x(I \times \R^3)} \lesssim \| u \|_{\dot S^1(I \times \R^3)}.
\end{equation}
Indeed, the $\dot S^1$ norm controls the following spacetime norms:

\begin{lemma}\label{4-lemma}\cite{tao:focusing}  For any Schwartz function $u$ on $I \times \R^3$, we have
\begin{equation}\label{strichartz-components}
\begin{split}
&\| \nabla u \|_{L^\infty_t L^2_x}
+ \| \nabla u \|_{L^{10}_t L^{30/13}_x}
+ \| \nabla u \|_{L^5_t L^{30/11}_x}
+ \| \nabla u \|_{L^4_t L^3_x}
+ \| \nabla u \|_{L^{10/3}_{t,x}}
+ \| \nabla u \|_{L^2_t L^6_x} + \\
&
+ \| u \|_{L^4_t L^\infty_x} + \| u \|_{L^6_t L^{18}_x} + \| u
\|_{L^{10}_{t,x}} + \| u \|_{L^\infty_t L^6_x}
\lesssim \| u \|_{\dot S^1}.
\end{split}
\end{equation}
where all spacetime norms are on $I \times \R^3$.
\end{lemma}

\begin{proof} All of these estimates follow from \eqref{sr} and Sobolev embedding except for
the $L^4_t L^\infty_x$ norm, which is a little more delicate because endpoint Sobolev embedding does not work
at $L^\infty_x$.  
Write
$$ c_N := \| P_N \nabla u \|_{L^2_t L^6_x} + \| P_N \nabla u \|_{L^\infty_t L^2_x},$$
then by the definition of $\dot S^1$ we have
$$ (\sum_N c_N^2)^{1/2} \lesssim \| u \|_{\dot S^1}.$$
On the other hand, for any dyadic frequency $N$ we see from Bernstein's inequality \eqref{bernstein} and \eqref{nabla-block} that
$$ N^{\frac{1}{2}} \| P_N u \|_{L^2_t L^\infty_x} \lesssim c_N  \, \hbox{ and } \,
N^{-\frac{1}{2}} \| P_N u \|_{L^\infty_t L^\infty_x} \lesssim c_N.$$
Thus, if we write $a_N(t) := \| P_N u(t) \|_{L^\infty_x}$, we have
\begin{equation}\label{an-l2}
 (\int_I a_N(t)^2\ dt)^{1/2} \lesssim N^{-\frac{1}{2}} c_N.
\end{equation}
and
\begin{equation}\label{an-sup}
 \sup_{t \in I} a_N(t) \lesssim N^{\frac{1}{2}} c_N
\end{equation}
Let us now compute
$$ \| u \|_{L^4_t L^\infty_x}^4 \lesssim \int_I (\sum_N a_N(t))^4\ dt.$$
Expanding this out and using symmetry, we have
$$ \| u \|_{L^4_t L^\infty_x}^4 \lesssim \sum_{N_1 \geq N_2 \geq N_3 \geq N_4} \int_I a_{N_1}(t) a_{N_2}(t) a_{N_3}(t) a_{N_4}(t)\ dt.$$
Estimating the two highest frequencies using \eqref{an-l2} and the lowest two using \eqref{an-sup}, we can
bound this by
$$ \lesssim  \sum_{N_1 \geq N_2 \geq N_3 \geq N_4} \frac{N_3^{\frac{1}{2}} N_4^{\frac{1}{2}}}{N_1^{\frac{1}{2}} N_2^{\frac{1}{2}}} c_{N_1} c_{N_2} c_{N_3} c_{N_4}.$$
Let $\tilde c_N$ denote the quantity
$$ \tilde c_N := \sum_{N'} \min( N/N', N'/N )^{1/10} c_{N'}.$$
Clearly we can bound the previous expression by
$$ \lesssim  \sum_{N_1 \geq N_2 \geq N_3 \geq N_4} \frac{N_3^{\frac{1}{2}} N_4^{\frac{1}{2}}}{N_1^{\frac{1}{2}} N_2^{\frac{1}{2}}} \tilde c_{N_1} \tilde c_{N_2} \tilde c_{N_3} \tilde c_{N_4}.$$
But we have $\tilde c_{N_j} \lesssim (N_1/N_j)^{1/10} \tilde c_{N_1}$ for $j=2,3,4$, hence we can bound the above by
$$  \lesssim  \sum_{N_1 \geq N_2 \geq N_3 \geq N_4} \frac{N_3^{\frac{1}{2}} N_4^{\frac{1}{2}}}{N_1^{\frac{1}{2}} N_2^{\frac{1}{2}}} \tilde c_{N_1}^4
(N_1/N_2)^{1/10} (N_1/N_3)^{1/10} (N_1/N_4)^{1/10}.$$
Summing in $N_4$, then in $N_3$, then in $N_2$, this is bounded by
$$ \sum_{N_1} \tilde c_{N_1}^4 \lesssim (\sum_N \tilde c_N^2)^2.$$
But by Young's inequality this is bounded by $\lesssim (\sum\limits_N c_N^2)^2 \lesssim \| u \|_{\dot S^1}^4$,
and the claim follows.
\end{proof}

We have the following standard Strichartz estimates:

\begin{lemma}\label{disjointed-strichartz}  Let $I$ be a compact time interval, and let $u: I \times \R^3  \to \C$
be a Schwartz solution to the forced Schr\"odinger equation
$$ iu_t + \Delta u = \sum_{m=1}^M F_m$$
for some Schwartz functions $F_1, \ldots, F_M$.  Then we have
\begin{equation}\label{strichartz}
\| u \|_{\dot S^k(I \times \R^3)} \lesssim \| u(t_0) \|_{\dot H^k(\R^3)}
+ C \sum_{m=1}^M \| \nabla^k F_m\|_{L^{q_m'}_t L^{r_m'}_x(I \times \R^3)}
\end{equation}
for any integer $k \geq 0$, any time $t_0
\in I$, and any admissible exponents $(q_1,r_1), \ldots, (q_m,r_m)$, where we use
$p'$ to denote the dual exponent to $p$, thus $1/p' + 1/p = 1$.
\end{lemma}

\begin{proof}
We first observe that we may take $M=1$, since the claim for
general $M$ then follows from the principle of superposition (exploiting
the linearity of the operator $(i\partial_t + \Delta)$, or equivalently using
the Duhamel formula \eqref{duhamel}) and the triangle
inequality.  We may then take $k=0$, since the
estimate for higher $k$ follows simply by applying $\nabla^k$ to both
sides of the equation and noting that this operator commutes with
$(i\partial_t + \Delta)$.  The  Littlewood-Paley projections $P_N$
also commute with $(i \partial_t + \Delta)$, and so we have
$$ (i\partial_t + \Delta) P_N u = P_N F_1$$
for each $N$.  From the Strichartz estimates in \cite{tao:keel} we
obtain
$$ \| P_N u \|_{L^q_t L^r_x(I \times \R^3)}
\lesssim \| P_N u(t_0) \|_{L^2(\R^3)}
+ \| P_N F_1\|_{L^{q_1'}_t L^{r_1'}_x(I \times \R^3)}$$
for any admissible exponents $(q,r),~(q_1, r_1)$.
Finally, we square, sum this in $N$ and
use the dual of \eqref{squarefunction} to obtain the result.
\end{proof}

\begin{remark}
\label{strich-flex}  In practice we shall take $k=0,1,2$ and $M=1,2$, and
$(q_m,r_m)$ to be either $(\infty,2)$ or $(2,6)$, i.e. we shall measure part of the
inhomogeneity in $L^1_t \dot H^k_x$, and the other part in $L^2_t \dot W^{k, 6/5}_x$.
\end{remark}

\subsection{Bilinear Strichartz estimate}

It turns out that to control the interactions between very high frequency and very low frequency
portions of the Schr\"odinger solution $u$, Strichartz estimates are insufficient,
and we need the following bilinear refinement, which we state in arbitrary dimension (though
we need it only in dimension $d=3$).


\begin{lemma}\label{bil-strichartz}  Let $d \geq 2$.  For any spacetime slab $I_*
\times \R^d$, any $t_0 \in I_*$, and for any $\delta > 0$, we have
\begin{equation}
\label{bilinearendpoint}
\begin{split}
\| u v \|_{L^2_t L^2_x(I_* \times \R^d)} \leq C(\delta) & (\|
u(t_0) \|_{\dot{H}^{-1/2+\delta}} + \| (i \partial_t + \Delta) u
\|_{L^1_t \dot{H}^{-1/2+\delta}_x}) \\
&\times(\| v(t_0) \|_{\dot{H}^{{\frac{d-1}{2}}-\delta}} +
\| (i \partial_t + \Delta) v \|_{L^1_t \dot{H}^{{\frac{d-1}{2}}-\delta}_x}).
\end{split}
\end{equation}
\end{lemma}

This estimate is very useful when $u$ is high frequency and $v$ is low
frequency, as it moves plenty of derivatives onto the low frequency
term.  This estimate shows in particular that there is little
interaction between high and low frequencies; this heuristic will
underlie many of our arguments to come, especially when we begin to
control the movement of mass, momentum, and energy from high modes to
low or vice versa. This estimate is essentially the refined Strichartz
estimate of Bourgain in \cite{BRefine} (see also \cite{borg:book}).
We make the trivial remark that the $L^2_{t,x}$ norm of $uv$ is the
same as that of $u\overline{v}$, $\overline{u}v$, or $\overline{u}\overline{v}$,
thus the above estimate also applies to expressions of the form $\O(uv)$.

\begin{proof}  We fix $\delta$, and allow our implicit constants to depend on $\delta$.
We begin by addressing the homogeneous case, with $u(t) := e^{it\Delta} \zeta$
and $v(t) := e^{it\Delta} \psi$  and consider the more
general problem of proving
\begin{equation}
  \label{alphaendpoint}
  {{\| u v \| }_{L^2_{t,x}}} \lesssim {{\| \zeta
\|}_{\dot{H}^{\alpha_1}}}{{\| \psi \|}_{\dot{H}^{\alpha_2}}}.
\end{equation}
Scaling invariance of
this estimate demands that $\alpha_1 + \alpha_2 = \frac{d}{2} -1.$ Our
first goal is to prove this for $\alpha_1 = -\half + \delta$ and $\alpha_2 =
\frac{d-1}{2} -\delta.$
The estimate \eqref{alphaendpoint} may be recast using duality and
renormalization as
\begin{equation}
\label{alpharecast}
  \int g(\xi_1 + \xi_2 , |\xi_1|^2 + |\xi_2|^2 )
   {{| \xi_1 |}^{- \alpha_1 }} {\widehat{\zeta}} (\xi_1 )
  {{| \xi_2 |}^{- \alpha_2 }} {\widehat{\psi}} (\xi_2 ) d\xi_1 d\xi_2
\lesssim {{\| g \|}_{L^2}} {{\| \zeta \|}_{L^2}} {{\| \psi
\|}_{L^2}}   .
\end{equation}
Since $\alpha_2 \geq \alpha_1$, we may restrict attention to the
interactions with $|\xi_1 | \geq |\xi_2|$. Indeed, in the remaining
case
we can multiply by $(\frac{|\xi_2|}{|\xi_1|})^{\alpha_2 - \alpha_1}
\geq 1$ to return to the case under consideration. In fact, we may
further restrict attention to the case where $|\xi_1| > 4 |\xi_2|$
since, in the other case, we can move the frequencies between
the two factors and reduce to the case
where $\alpha_1 = \alpha_2$, which can be treated by $L^4_{t,x}$
Strichartz estimates\footnote{In one dimension $d=1$,
Lemma \ref{bil-strichartz} fails when $u, v$ have comparable
frequencies, but continues to hold when $u, v$ have separated frequencies; see
\cite{ckstt:dnls} for further discussion.} when $d \geq 2$.
Next, we decompose $|\xi_1|$ dyadically and $|\xi_2|$ in dyadic
multiples of the size of $|\xi_1|$ by rewriting the quantity to be
controlled as ($N, \Lambda$ dyadic):
\begin{equation*}
  \sum_N \sum_\Lambda \int \int g_N (\xi_1 + \xi_2 , |\xi_1|^2 +
  |\xi_2 |^2) {{| \xi_1 |}^{-\alpha_1 }} \widehat{\zeta_N} (\xi_1)
{{| \xi_2 |}^{-\alpha_2 }} \widehat{\psi_{\Lambda N}} (\xi_2) d\xi_1
d\xi_2 .
\end{equation*}
Note that subscripts on $g,\zeta,\psi$ have been inserted to evoke the
localizations to $|\xi_1 + \xi_2 | \thicksim N, |\xi_1 | \thicksim N,
|\xi_2 | \thicksim \Lambda N$, respectively. Note that in the
situation we are considering here, namely $|\xi_1 | \geq 4 |\xi_2 |$,
we have that $|\xi_1 + \xi_2 | \thicksim |\xi_1 |$ and this explains
why $g$ may be so localized.

By renaming components, we may assume that $|\xi_1^1 | \thicksim
|\xi_1 |$ and $|\xi_2^1 | \thicksim |\xi_2 |$. Write $\xi_2 = (\xi_2^1
, \underline{\xi_2}).$ We now change variables by writing $u= \xi_1 +
\xi_2,~v=|\xi_1|^2 +
|\xi_2|^2$ and $du dv = J d\xi_2^1 d\xi_1$. A calculation
then shows that $J = |2 (\xi_1^1 \pm \xi_2^1)| \thicksim |\xi_1 |.$
Therefore, upon changing variables in the inner two integrals, we
encounter
\begin{equation*}
  \sum_N N^{-\alpha_1 } \sum_{\Lambda \leq 1} (\Lambda N)^{- \alpha_2 }
  \int_{\R^{d-1}} \int_{\R} \int_{\R^d} g_N (u,v) H_{N,\Lambda } (u, v,
  \underline{\xi_2} ) du dv d\underline{\xi_2}
\end{equation*}
where
\begin{equation*}
  H_{N,\Lambda} ( u,v, \underline{\xi_2}) = \frac{ \widehat{\zeta_N}
  (\xi_1 ) \widehat{\psi_{\Lambda N}} (\xi_2)}{J}.
\end{equation*}
We apply Cauchy-Schwarz on the $u,v$ integration and change back to
the original variables to obtain
\begin{equation*}
  \sum_N N^{-\alpha_1 } {{\| g_N \|}_{L^2}}  \sum_{\Lambda \leq 1} (\Lambda N)^{- \alpha_2 }
  \int_{\R^{d-1}} {{\left[ \int_{\R} \int_{\R^{d-1}} \frac{ |\widehat{\zeta_N}
  (\xi_1 )|^2 \widehat{|\psi_{\Lambda N}} (\xi_2)|^2}{J} d\xi_1
  d\xi_2^1 \right] }^\half } d\underline{\xi_2} .
\end{equation*}
We recall that $J \thicksim N$ and use Cauchy-Schwarz in the
$\underline{\xi_2}$ integration, keeping in mind the localization
$|\xi_2 | \thicksim \Lambda N$, to get
\begin{equation*}
  \sum_N N^{- \alpha_1 - \half} {{\| g_N \|}_{L^2}} \sum_{\Lambda \leq
1} (\Lambda
  N)^{-\alpha_2 + \frac{d-1}{2}} {{\| \widehat{\zeta_N} \|}_{L^2}} {{\|
  \widehat{\psi_{\Lambda N} \|}_{L^2}}}.
\end{equation*}
Choose $\alpha_1 = -\half + \delta$ and $\alpha_2 = \frac{d-1}{2} -
\delta$ with $\delta >0$ to obtain
\begin{equation*}
  \sum_N {{\| g_N \|}_{L^2}} {{\| \widehat{\zeta_N} \|}_{L^2}}
  \sum_{\Lambda \leq 1} \Lambda^\delta  {{\| \widehat{\psi_{\Lambda
N} \|}_{L^2}}}
\end{equation*}
which may be summed up to give the claimed homogeneous estimate.

We turn our attention to the inhomogeneous estimate
\eqref{bilinearendpoint}. For simplicity we set
$F:=(i\partial_t +\Delta)u$ and $G:=(i\partial_t +\Delta)v$.
Then we use Duhamel's formula \eqref{duhamel} to write
$$u=e^{i(t-t_0)\Delta} u(t_0) - i\int_{t_0}^t e^{i(t-t')\Delta} F(t')\,dt', \, \, \
v=e^{i(t-t_0)\Delta} v(t_0) - i\int_{t_0}^t e^{i(t-t')\Delta} G(t').$$
We obtain\footnote{Alternatively, one can absorb the homogeneous components $e^{i(t-t_0)\Delta} u(t_0)$,
$e^{i(t-t_0)\Delta} v(t_0)$ into the inhomogeneous term by adding an artificial forcing term of $\delta(t-t_0) u(t_0)$
and $\delta(t-t_0) v(t_0)$ to $F$ and $G$ respectively, where $\delta$ is the Dirac delta.}
\begin{align*}
\|uv\|_{L^2} &\lesssim \left\|e^{i(t-t_0)\Delta} u(t_0) e^{i(t-t_0)\Delta} v(t_0)\right\|_{L^2}\\
&+\left\|e^{i(t-t_0)\Delta} u(t_0)\int_{t_0}^t e^{i(t-t')\Delta} G(t') \,dt' \right\|_{L^2}
+\left\|e^{i(t-t_0)\Delta} v(t_0)\int_{t_0}^t e^{i(t-t')\Delta} F(t')dt' \right\|_{L^2} \\
&+ \left\|\int_{t_0}^t e^{i(t-t')\Delta} F(t')dt' \int_{t_0}^t e^{i(t-t'')\Delta } G(x,t'')\,dt''
\right\|_{L^2}\\
&:= I_1 + I_2 + I_3 + I_4.
\end{align*}
The first term was treated in the first part of the proof.
The second and the third are similar so we consider only $I_2$. Using
the Minkowski inequality we have
$$I_2\lesssim \int_{\R} \|e^{i(t-t_0)\Delta} u(t_0) e^{i(t-t')\Delta} G(t')\|_{L^2}
\,dt',$$
and in this case the lemma follows from the homogeneous estimate
proved above.
Finally, again by Minkowski's inequality we have
$$I_4\lesssim \int_{\R} \int_{\R}\| e^{i(t-t')\Delta} F(t') e^{i(t-t'')\Delta} G(t'')\|_{L^2_x}
dt' dt'',$$
and the proof follows by inserting in the integrand the
homogeneous estimate above.
\end{proof}

\begin{remark}
\label{BRS}
In the situation where the initial data are dyadically localized in
frequency space, the estimate \eqref{alphaendpoint} is valid
\cite{BRefine} at the
endpoint $\alpha_1 = - \half, \alpha_2 =
\frac{d-1}{2}$. Bourgain's argument also establishes the result with $
\alpha_1 = -\half + \delta, \alpha_2 = \frac{d-1}{2} + \delta$, which
is not scale invariant.
However, the full estimate fails at the endpoint. This can be
seen by calculating the left and right sides of \eqref{alpharecast} in
the situation where $\widehat{\zeta_1} = \chi_{R_1}$ with $R_1 = \{
\xi : \xi_1 = N e^1 + O (N^\half) \}$ (where $e^1$ denotes the first
coordinate unit vector), $\widehat{\psi_2} (\xi_2) =
|\xi_2|^{-\frac{d-1}{2}} \chi_{R_2}$ where $R_2 = \{ \xi_2: 1 \ll
|\xi_2 | \ll N^\half, \xi_2 \cdot e^1 = O(1) \}$ and $g(u,v) =
\chi_{R_0} (u,v)$ with $R_0 = \{ (u,v): u = N e^1 + O(N^\half), v =
|u|^2 + O(N)\}.$ A calculation then shows that the left side of
\eqref{alpharecast} is of size $N^{\frac{d+1}{2}} \log N$ while the
right side is of size $N^{\frac{d+1}{2}} (\log N)^\half$. Note that
the same counterexample shows that the estimate
\begin{equation*}
  {{\| u \overline{v} \|}_{L^2_{t,x}}} \lesssim {{\| \zeta
  \|}_{\dot{H}^\alpha_1}}{{\| \psi \|}_{\dot{H}^\alpha_2}},
\end{equation*}
where $u(t) = e^{it\Delta} \zeta,~v(t) = e^{it\Delta} \psi$,
 also fails at the endpoint.
\end{remark}

\subsection{Quintilinear Strichartz estimates}\label{quint-sec}

We record  the following useful inequality:

\begin{lemma}\label{leibnitz-holder}  For any $k=0,1,2$ and any slab
$I \times \R^3$, and any smooth functions $v_1,\ldots,v_5$ on this
slab, we have
\begin{equation}\label{12}
\| \nabla^k\O(v_1 v_2 v_3 v_4 v_5) \|_{L^1_t L^2_x} \lesssim
\sum_{ \{a,b,c,d,e\} = \{1,2,3,4,5\} }
\| v_a \|_{\dot S^1} \| v_b \|_{\dot S^1} \| v_c \|_{L^{10}_{x,t}} \| v_d \|_{L^{10}_{x,t}}
|| v_e ||_{\dot S^k}
\end{equation}
where all the spacetime norms are on the slab $I \times \R^3$. In a similar spirit, we have
\begin{equation}\label{26}
\| \nabla\O(v_1 v_2 v_3 v_4 v_5) \|_{L^2_t L^{6/5}_x} \lesssim
\prod_{j=1}^5 \| \nabla v_j \|_{L^{10}_t L^{30/13}_x} \leq \prod_{j=1}^5 \| v_j \|_{\dot S^1}.
\end{equation}
\end{lemma}

\begin{proof}
Consider, for example, the $k=1$ case of \eqref{12}. Applying the
Leibnitz rule, we encounter various terms to control including
\begin{equation*}
  {{\| \O((\nabla v_1) v_2 v_3 v_4 v_5) \|}_{L^1_t L^2_x}} \lesssim
{\| \nabla v_1 \|}_{L^{\frac{10}{3}}_{t,x}}
{\| v_2 \|}_{L^4_t L^\infty_x}
{\| v_3 \|}_{L^4_t L^\infty_x}
{\| v_4 \|}_{L^{10}_{t,x}}
{\| v_5 \|}_{L^{10}_{t,x}}.
\end{equation*}
The claim follows then using \eqref{strichartz-components}.
The $k=2$ case of \eqref{12} follows similarly using estimates such as
\begin{equation*}
  {{\| \O((\nabla^2 v_1) v_2 v_3 v_4 v_5) \|}_{L^1_t L^2_x}} \lesssim {{\|
  \nabla^2 v_1 \|}}_{L^{\frac{10}{3}}_{t,x}}
\| v_2 \|_{L^4_t L^\infty_x}
\| v_3 \|_{L^4_t L^\infty_x}
{\| v_4 \|}_{L^{10}_{t,x}}
{\| v_5 \|}_{L^{10}_{t,x}}
\end{equation*}
and
\begin{equation*}
  {{\| \O((\nabla v_1) (\nabla v_2) v_3 v_4 v_5) \|}_{L^1_t L^2_x}} \lesssim {{\|
  \nabla v_1 \|}}_{L^{\frac{10}{3}}_{t,x}}
\| \nabla v_2 \|_{L^4_t L^\infty_x}
\| v_3 \|_{L^4_t L^\infty_x}
{\| v_4 \|}_{L^{10}_{t,x}}
{\| v_5 \|}_{L^{10}_{t,x}}
\end{equation*}
and the $k=0$ case is similar (omit all the $\nabla$s).

Finally, estimate \eqref{26} similarly follows from the Sobolev
embedding $\| u \|_{L^{10}_{t,x}} \lesssim \| \nabla u \|_{L^{10}_t L^{30/13}_x}$,
\eqref{strichartz-components}  and H\"older's inequality,
\begin{equation*}
\| \O(\nabla v_1 v_2 v_3 v_4 v_5) \|_{L^2_t L^{6/5}_x} \lesssim
\|\nabla v_1\|_{L^{10}_t L^{30/13}_x} 
\| v_2\|_{L^{10}_{t,x}}
\| v_3\|_{L^{10}_{t,x}}
\| v_4\|_{L^{10}_{t,x}}
\| v_5\|_{L^{10}_{t,x}}.
\end{equation*}

\end{proof}

We need a variant of the above lemma which also exploits the bilinear
Strichartz inequality in Lemma \ref{bil-strichartz} to obtain a gain when some of the factors are ``high frequency''
and others are ``low frequency''.

\begin{lemma}\label{quint-improved}  Suppose $v_{hi}$, $v_{lo}$ are functions on $I \times \R^3$ such that
\begin{align*}
\| v_{hi} \|_{\dot S^0} + \| (i\partial_t + \Delta) v_{hi} \|_{L^1_t L^2_x(I \times \R^3)} &\lesssim \eps K\\
\| v_{hi} \|_{\dot S^1} + \| \nabla (i\partial_t + \Delta) v_{hi} \|_{L^1_t L^2_x(I \times \R^3)} &\lesssim K\\
\| v_{lo} \|_{\dot S^1} + \| \nabla (i\partial_t + \Delta) v_{lo} \|_{L^1_t L^2_x(I \times \R^3)} &\lesssim K\\
\| v_{lo} \|_{\dot S^2} + \| \nabla^2 (i\partial_t + \Delta) v_{lo} \|_{L^1_t L^2_x(I \times \R^3)} &\lesssim \eps K
\end{align*}
for some constants $K > 0$ and $0 < \eps \ll 1$.  Then for any $j=1,2,3,4$, we have
$$
\| \nabla \O( v_{hi}^j v_{lo}^{5-j} ) \|_{L^2_t L^{6/5}_x(I \times \R^3)} \lesssim
\eps^{\frac{9}{10}} K^5.
$$
\end{lemma}


\begin{remark}  The point here is the gain of $\eps^{9/10}$, which cannot be obtained directly from the type of
arguments used to prove Lemma \ref{leibnitz-holder}.  As the proof will reveal, one can replace the
exponent $9/10$ with any exponent less than one, though for our purposes all that matters is that the power of $\eps$
is positive.  The $\dot S^0$ bound on $v_{hi}$ effectively restricts $v_{hi}$ to high frequencies (as the low and medium frequencies
will then be very small in $\dot S^1$ norm); similarly, the $\dot S^2$ control on $v_{lo}$ effectively restricts $v_{lo}$ to
low frequencies.  This Lemma is thus an assertion that the components of the nonlinearity in \eqref{nls} arising from
interactions between low and high frequencies are rather weak; this phenomenon underlies the important
frequency localization result in Proposition \ref{frequency-decoupling}, but the motif of controlling
the interaction between low and high frequencies underlies many other parts of our argument also, notably
in Proposition \ref{Morawetz} and Proposition \ref{energy-travel}.
\end{remark}

\begin{proof}
Throughout this proof all spacetime norms shall  be on $I \times
\R^3$.
We may normalize $K := 1$.  By the Leibnitz rule we have
$$ \| \nabla \O( v_{hi}^j v_{lo}^{5-j} ) \|_{L^2_t L^{6/5}_x} \lesssim
\| \O(v_{hi}^j v_{lo}^{4-j} \nabla v_{lo})\|_{L^2_t L^{6/5}_x}
+ \| \O(v_{hi}^{j-1} v_{lo}^{5-j} \nabla v_{hi}) \|_{L^2_t L^{6/5}_x}.$$
Consider the $\nabla v_{lo}$ terms first, which are rather easy.  By H\"older we have
$$
\| \O(v_{hi}^j v_{lo}^{4-j} \nabla v_{lo})\|_{L^2_t L^{6/5}_x}
\lesssim
\| \nabla v_{lo} \|_{L^\infty_t L^6_x} \| v_{hi}
\|_{L^\infty_t L^2_x}
 \| v_{lo} \|_{L^6_t L^{18}_x}^{4-j} \| v_{hi} \|_{L^6_t L^{18}_x}^{j-1}.$$
Applying \eqref{strichartz-components}, this is bounded by
$$ \lesssim
\| v_{lo} \|_{\dot S^2} \| v_{hi} \|_{\dot S^0} \| v_{lo} \|_{\dot S^1}^{4-j}
\| v_{hi} \|_{\dot S^1}^{j-1} \lesssim \eps^2$$
which is acceptable.

Now consider the $\nabla v_{hi}$ terms, which are more difficult.  First
consider the $j=2,3,4$ cases.  By H\"older we have
$$
\| \O(v_{hi}^{j-1} v_{lo}^{5-j} \nabla v_{hi}) \|_{L^2_t L^{6/5}_x}
\lesssim
\| \nabla v_{hi} \|_{L^2_t L^6_x} \| v_{lo} \|_{L^\infty_t
L^\infty_x}
\| v_{hi} \|_{L^\infty_t L^2_x}^{1/2} \| v_{lo} \|_{L^\infty_t L^6_x}^{4-j} \|
v_{hi} \|_{L^\infty_t L^6_x}^{j-3/2}.$$
Now observe (for instance from \eqref{bernstein}, \eqref{nabla-block} and dyadic
decomposition) that
$$ \| v_{lo} \|_{L^\infty_t L^\infty_x}
\lesssim  \| v_{lo} \|_{L^\infty_t L^6_x}^{1/2}
\| \nabla v_{lo} \|_{L^\infty_t L^6_x}^{1/2},$$
so by \eqref{strichartz-components} we have
$$ \| \O(v_{hi}^{j-1} v_{lo}^{5-j} \nabla v_{hi}) \|_{L^2_t L^{6/5}_x}
 \lesssim \| v_{hi} \|_{\dot S^1}^{j-1/2} \| v_{hi} \|_{\dot S^0}^{1/2}
\| v_{lo} \|_{\dot S^2}^{1/2} \| v_{lo} \|_{\dot S^1}^{9/2-j}$$
which is $O(\eps^{9/10})$, and is acceptable.

Finally consider the $j=1$ term.  For this term we must use dyadic decomposition, writing
$$ \| \O(v_{lo}^4 \nabla v_{hi}) \|_{L^2_t L^{6/5}_x}
\lesssim \sum_{N_1,N_2,N_3,N_4} \| \O((P_{N_1} v_{lo})
(P_{N_2} v_{lo}) (P_{N_3} v_{lo}) (P_{N_4} v_{lo})
 \nabla v_{hi}) \|_{L^2_t L^{6/5}_x}.$$
By symmetry we may take $N_1 \geq N_2 \geq N_3 \geq N_4$.  We then estimate this using H\"older by
$$ \sum_{N_1 \geq N_2 \geq N_3 \geq N_4}
\| \O( P_{N_1} v_{lo} \nabla v_{hi} ) \|_{L^2_t L^2_x}
\| P_{N_2} v_{lo} \|_{L^\infty_t L^6_x}
\| P_{N_3} v_{lo} \|_{L^\infty_t L^6_x}
\| P_{N_4} v_{lo} \|_{L^\infty_t L^\infty_x}.$$
The middle two factors can be estimated by $\| v_{lo} \|_{\dot S^1} = O(1)$.  The last factor can be estimated
using Bernstein \eqref{nabla-block} either as
$$ \| P_{N_4} v_{lo} \|_{L^\infty_t L^\infty_x} \lesssim N_4^{1/2}
\| P_{N_4} v_{lo} \|_{L^\infty_t L^6_x} \lesssim N_4^{1/2} \| v_{lo} \|_{\dot S^1} \lesssim N_4^{1/2}$$
or as
$$ \| P_{N_4} v_{lo} \|_{L^\infty_t L^\infty_x} \lesssim N_4^{-1/2}
\| \nabla P_{N_4} v_{lo} \|_{L^\infty_t L^6_x} \lesssim N_4^{-1/2} \| v_{lo} \|_{\dot S^2} \lesssim \eps N_4^{-1/2}.$$
Meanwhile, the first factor can be estimated using \eqref{bilinearendpoint} as
\begin{align*}
\| \O( P_{N_1} v_{lo} \nabla v_{hi} ) \|_{L^2_t L^2_x} &\lesssim
(\| \nabla v_{hi}(t_0) \|_{\dot{H}^{-1/2+\delta}} + \| (i \partial_t + \Delta) \nabla v_{hi}
\|_{L^1_t \dot{H}^{-1/2+\delta}_x}) \\
&\times(\| P_{N_1} v_{lo}(t_0) \|_{\dot{H}^{1-\delta}} +
\| (i \partial_t + \Delta) P_{N_1} v_{lo} \|_{L^1_t \dot{H}^{1-\delta}_x}),
\end{align*}
where $t_0 \in I$ is an arbitrary time and $0 < \delta < 1/2$ is an arbitrary exponent.
From the hypotheses on $v_{hi}$ and interpolation we see that
$$ \| \nabla v_{hi}(t_0) \|_{\dot{H}^{-1/2+\delta}} + \| (i \partial_t + \Delta) \nabla v_{hi}
\|_{L^1_t \dot{H}^{-1/2+\delta}_x} \lesssim \eps^{1/2-\delta}$$
while from the hypotheses on $v_{lo}$ and \eqref{nabla-block} we see that
$$ \| P_{N_1} v_{lo}(t_0) \|_{\dot{H}^{1-\delta}} +
\| (i \partial_t + \Delta) P_{N_1} v_{lo} \|_{L^1_t \dot{H}^{1-\delta}_x})
\lesssim N_1^{-\delta}.$$
Putting this all together, we obtain
$$ \| \O(v_{lo}^4 \nabla v_{hi}) \|_{L^2_t L^{6/5}_x}
\lesssim \sum_{N_1 \geq N_2 \geq N_3 \geq N_4} \eps^{1/2-\delta} N_1^{-\delta} \min(N_4^{1/2}, \eps N_4^{-1/2}).$$
Performing the $N_1$ sum, then the $N_2$, then the $N_3$, then the $N_4$, we obtain the desired bound of $O(\eps^{9/10})$,
if $\delta$ is sufficiently small.
\end{proof}

\subsection{Local well-posedness and perturbation theory}\label{perturb-sec}


It is well known (see e.g. \cite{borg:book}) that the equation \eqref{nls} is \emph{locally}
well-posed\footnote{In particular, we have uniqueness of this Cauchy
  problem, at least under the assumption that $u$ lies in
  $L^{10}_{t,x} \cap C^0_t  {\dot{H}}^1_x$, and so whenever we
  construct a  solution $u$ to \eqref{nls} with specified initial
data $u(t_0)$, we will refer to it as \emph{the} solution to
\eqref{nls} with this data.} in $\dot H^1(\R^3)$, and indeed that this well-posedness extends to any time interval on which
one has a uniform bound on the $L^{10}_{t,x}$ norm; this can already be seen from Lemma \ref{leibnitz-holder} and
\eqref{strichartz} (see also Lemma \ref{persistence} below).  In this section we detail some variants of the local
well-posedness
argument which
describe how we can perturb finite-energy solutions (or near-solutions)
to \eqref{nls} in the energy norm when we control the original solution in the
$L^{10}_{t,x}$ norm and the error of near-solutions in a dual Strichartz space. The arguments we give are similar to those in previous work such as
\cite{borg:book}.

 We begin with
a preliminary result where the near solution, the error of the near-solution, and
the free evolution of the perturbation are all assumed to be small in \emph{spacetime} norms,
but allowed to be large in energy norm.

\begin{lemma}[Short-time perturbations]\label{perturbtiny}  Let $I$ be a compact interval, and
let
$\tilde u$ be a function on $I \times \R^3$ which is a near-solution
to \eqref{nls}
in the sense that
\begin{equation}\label{approx-nls}
 (i \partial_t + \Delta) \tilde u = |\tilde u|^4 \tilde u + e
\end{equation}
for some function $e$.
Suppose that we also have the energy bound
$$ \| \tilde u \|_{L^\infty_t \dot H^1_x(I \times \R^3)} \leq E$$
for some $E > 0$.
Let $t_0 \in I$, and let $u(t_0)$ be close to $\tilde u(t_0)$ in the sense that
\begin{equation}\label{uu-energy}
 \| u(t_0) - \tilde u(t_0) \|_{\dot H^1_x} \leq E'
\end{equation}
for some $E' > 0$. 
Assume also that we have the  smallness conditions
\begin{align}
\label{assumetiny}
 \| \nabla \tilde u \|_{L^{10}_t L^{30/13}_x(I \times \R^3)} &\leq \eps_0 \\
\label{assumetiny-too}
 \| \nabla e^{i(t-t_0)\Delta} (u(t_0) - \tilde u(t_0)) \|_{L^{10}_t L^{30/13}_x(I \times \R^3)} &\leq \eps\\
\label{error-tiny-also}
\| \nabla e \|_{L^2_t L^{6/5}_x} &\leq \eps
\end{align}
for some $0 < \eps < \eps_0$, where $\eps_0$ is some constant $\eps_0 = \eps_0(E,E') > 0$.

We conclude that there exists a solution $u$ to \eqref{nls} on $I \times \R^3$ with
the specified initial data
$u(t_0)$ at $t_0$, and furthermore
\begin{align}
\label{showshow}
\| u - \tilde u \|_{\dot S^1(I \times \R^3)} &\lesssim E'\\
\label{showshow-1a}
\| u \|_{\dot S^1(I \times \R^3)} &\lesssim E' + E\\
\label{showshow-2}
\| u - \tilde u \|_{L^{10}_{t,x}(I \times \R^3)} \lesssim
\| \nabla (u - \tilde u) \|_{L^{10}_t L^{30/13}_x(I \times \R^3)} 
&\lesssim \eps\\
\label{showshow-3}
\| \nabla (i \partial_t + \Delta) (u - \tilde u) \|_{L^2_t L^{6/5}_x(I \times \R^3)} &\lesssim \eps.
\end{align}
\end{lemma}

Note that $u(t_0) - \tilde u(t_0)$ is allowed to have large energy, albeit
at the cost of forcing $\eps$ to be smaller, and worsening the bounds in \eqref{showshow}.
From the Strichartz estimate \eqref{strichartz}, \eqref{uu-energy} we see that the hypothesis \eqref{assumetiny-too} is redundant
if one is willing to take $E' = O(\eps)$.

\begin{proof}
By the well-posedness theory reviewed above, it suffice to prove
\eqref{showshow} - \eqref{showshow-3} as
 \emph{a priori} estimates\footnote{That is, we may assume the solution $u$
already exists and is smooth on the entire interval $I$.}.
We establish these bounds for $t \geq t_0$, since the corresponding bounds
for the $t \leq t_0$ portion
of $I$ are proved similarly.

First note that the Strichartz estimate (Lemma \ref{disjointed-strichartz}), Lemma
\ref{leibnitz-holder} and \eqref{error-tiny-also} give,
\begin{align*}
\|\tilde u \|_{\dot S^1 ( I \times \R^3)} & \lesssim E + \| \tilde u
\|_{L^{10}_{t,x}(I \times
\R^3)} \cdot \| \tilde u \|_{\dot S^1(I \times \R^3)}^4 + \eps.
\end{align*}
By \eqref{assumetiny} and Sobolev embedding we have
$ \| \tilde u \|_{L^{10}_{t,x}(I \times \R^3)} \lesssim \eps_0$. A standard continuity argument in $I$
then gives (if $\eps_0$ is sufficiently small depending on $E$)
\begin{align}
\| \tilde u \|_{\dot S^1(I \times \R^3)} & \lesssim E.
\label{tildeubound}
\end{align}

Define $v := u - \tilde u$.
For each $t \in I$ define the quantity
$$
S(t) := \| \nabla (i \partial_t + \Delta) v  \|_{L^2_t L^{6/5}_x([t_0,t] \times \R^3)}.
$$

From using Lemma \ref{4-lemma}, Lemma \ref{disjointed-strichartz}, \eqref{assumetiny-too}, we have
\begin{align}
\| \nabla v \|_{L^{10}_t L^{30/13}_x([t_0,t] \times \R^3)} &\lesssim \| \nabla(v - e^{i(t-t_0)\Delta}v(t_0)) \|_{L^{10}_t L^{30/13}_x([t_0,t] \times \R^3)} \\
&\quad+ \| \nabla e^{i(t-t_0)\Delta} v(t_0) \|_{L^{10}_t L^{30/13}_x([t_0,t] \times \R^3)} \nonumber \\
&\lesssim \| v - e^{i(t-t_0)\Delta} v(t_0) \|_{\dot S^1([t_0,t] \times \R^3)} + \eps\nonumber\\
& \lesssim S(t)+\eps.\label{QR-bound}
\end{align}
On the other hand, since $v$ obeys the equation
\begin{align*}
(i \partial_t + \Delta) v & = |\tilde u + v|^4 (\tilde u + v) - |\tilde u|^4 \tilde u - e = \sum_{j=1}^5
\O(v^j \tilde u^{5-j}) - e
\end{align*}
by \eqref{5-schematic}, we easily check using \eqref{26}, \eqref{assumetiny}, \eqref{error-tiny-also}, \eqref{QR-bound} that 
\begin{align*}
S(t)&\lesssim \eps + \sum_{j=1}^5 (S(t) + \eps)^j \eps_0^{5-j}.
\end{align*}
If $\eps_0$ is sufficiently small, a standard continuity argument then yields
the bound $S(t) \lesssim \eps$ for all $t \in I$.  This gives \eqref{showshow-3}, and \eqref{showshow-2} follows
from \eqref{QR-bound}.  Applying Lemma \ref{disjointed-strichartz}, \eqref{uu-energy} we then conclude
\eqref{showshow} (if $\eps$ is sufficiently small), and then from \eqref{tildeubound} and the triangle inequality
we conclude \eqref{showshow-1a}.
\end{proof}

We will actually be more interested in iterating the above Lemma\footnote{We are grateful to Monica Visan for pointing out an incorrect version of
Lemma \ref{perturb} in a previous version of this paper, and also in simplifying the proof of Lemma \ref{perturbtiny}.}
to deal with the more general situation of
near-solutions with finite but arbitrarily large
$L^{10}_{t,x}$ norms.

\begin{lemma}[Long-time perturbations]\label{perturb}  Let $I$ be a compact interval, and let
$\tilde u$ be a function on $I \times \R^3$ which obeys the bounds
\begin{equation}\label{ume}
\| \tilde u \|_{L^{10}_{t,x}(I \times \R^3)} \leq M
\end{equation}
and
\begin{equation}\label{ume-2}
 \| \tilde u \|_{L^\infty_t \dot H^1_x(I \times \R^3)} \leq E
 \end{equation}
for some $M, E > 0$.  Suppose also that $\tilde u$ is a near-solution
to \eqref{nls} in the sense that it solves \eqref{approx-nls} for some $e$.
Let $t_0 \in I$, and let $u(t_0)$ be close to $\tilde u(t_0)$ in the sense that
$$ \| u(t_0) - \tilde u(t_0) \|_{\dot H^1_x} \leq E'$$
for some $E' > 0$.  Assume also that we have the smallness conditions,
\begin{align}
\label{safetycheckgeneral}
 \| \nabla e^{i(t-t_0)\Delta} (u(t_0) - \tilde u(t_0)) \|_{L^{10}_t L^{30/13}_x(I \times \R^3)} &\leq \eps\\
\nonumber 
\| \nabla e \|_{L^2_t L^{6/5}_x(I \times \R^3)} &\leq \eps
\end{align}
for some $0 < \eps < \eps_1$, where $\eps_1$ is some constant $\eps_1 = \eps_1(E,E',M) > 0$.
We conclude
there exists a solution $u$ to \eqref{nls} on $I \times \R^3$ with
the specified initial data $u(t_0)$ at $t_0$, and furthermore
\begin{align*}
 \| u - \tilde u \|_{\dot S^1(I \times \R^3)} & \leq C(M, E, E') \\
 \| u \|_{\dot S^1(I \times \R^3)} & \leq C(M, E, E') \\
 \| u - \tilde u \|_{L^{10}_{t,x}(I \times \R^3)} \leq
 \| \nabla(u - \tilde u) \|_{L^{10}_t L^{30/13}_x(I \times \R^3)}
 & \leq C(M, E, E') \eps.
 \end{align*}
\end{lemma}

Once again, the hypothesis \eqref{safetycheckgeneral} is redundant
by the Strichartz estimate if
one is willing to take $E' = O(\eps)$;
 however it will be useful in our  applications to know
that this Lemma can tolerate a perturbation which is large in the
energy  norm but whose
free evolution is small in the $L^{10}_t \dot W^{1,30/13}_x$ norm.

This lemma is already useful in the $e=0$ case, as it says that one
has local well-posedness in the energy space whenever the
$L^{10}_{t,x}$ norm is bounded; in fact one has locally Lipschitz
dependence on the initial data.  For similar perturbative results
see \cite{borg:scatter}, \cite{borg:book}.


\begin{proof}  As in the  previous proof,
we may assume that $t_0$ is the lower bound of the interval $I$.
Let $\eps_0 = \eps_0(E,2E')$ be as in Lemma \ref{perturbtiny}. (We need to replace
$E'$ by the slightly larger $2E'$ as the $\dot H^1$ norm of $u - \tilde u$ is going
to grow slightly in time.)  

The first step is to establish a $\dot S^1$ bound on $\tilde u$.  Using \eqref{ume} we may subdivide
$I$ into $C(M, \eps_0)$ time intervals such that the $L^{10}_{t,x}$ norm of $\tilde u$ is at most $\eps_0$ on 
each such interval.  By using \eqref{ume-2} and Lemmas \ref{disjointed-strichartz}, \ref{leibnitz-holder}
as in the proof of \eqref{tildeubound} we see that the $\dot S^1$ norm of $\tilde u$ is $O(E)$ on each of these intervals.  Summing up
over all the intervals we conclude
$$ \| \tilde u \|_{\dot S^1(I \times \R^3)} \leq C(M,E,\eps_0)$$
and in particular by Lemma \ref{4-lemma}
$$ \| \nabla \tilde u \|_{L^{10}_t L^{30/13}_x(I \times \R^3)} \leq C(M,E,\eps_0).$$
We can then subdivide the interval 
$I$ into $N \leq C(M,E,\eps_0)$ subintervals
$I_j \equiv [T_j, T_{j+1}]$ so that on each $I_j$ we have,
$$
\| \nabla \tilde u\|_{L^{10}_t L^{30/13}_x(I_j \times \R^3)}  \leq \eps_0.
$$
We can then verify inductively using Lemma \ref{perturbtiny} for each $j$ that if $\eps_1$ is
sufficiently small depending on $\eps_0$, $N$, $E$, $E'$, then we have
\begin{align*}
\|u - \tilde u\|_{\dot S^1(I_j \times \R^3)} &\leq C(j) E'\\
\|u \|_{\dot S^1(I_j \times \R^3)} &\leq C(j) (E' + E)\\
\|\nabla(u - \tilde u)\|_{L^{10}_t L^{30/13}_x(I_j \times \R^3)} &\leq C(j) \eps\\
\|\nabla (i \partial_t + \Delta)(u - \tilde u)\|_{L^2_t L^{6/5}_x(I_j \times \R^3)} &\leq C(j) \eps
\end{align*}
and hence by Strichartz \eqref{strichartz} and \eqref{strichartz-components} we have
\begin{align*}
&\| \nabla e^{i(t - T_{j+1})\Delta} (u(T_{j+1}) - \tilde u(T_{j+1})) \|_{L^{10}_t L^{30/13}_x(I \times \R^3)} \\
&\quad \leq \| \nabla e^{i(t - T_{j})\Delta} (u(T_{j}) - \tilde u(T_{j})) \|_{L^{10}_t L^{30/13}_x(I \times \R^3)} + C(j) \eps
\end{align*}
and
$$ \| u(T_{j+1}) - \tilde u(T_{j+1}) \|_{\dot H^1} \leq
\| u(T_{j}) - \tilde u(T_{j}) \|_{\dot H^1} + C(j) \eps$$ 
allowing one to continue the induction
(if $\eps_1$ is sufficiently small depending on $E$, $N$, $E'$, $\eps_0$, then the quantity in
\eqref{uu-energy} will not exceed $2E'$).  The claim follows.
\end{proof}

\begin{remark}  The value of $\eps_1$ given by the above lemma deteriorates exponentially with $M$, or
more precisely it behaves like $\exp(-M^C)$ in its dependence\footnote{With respect to all its
parameters, $\eps_1(E, E', M) \approx \exp(- M^C \langle E \rangle^C \langle E' \rangle^C)$.} on
$M$. As this lemma is used quite often in our argument, this will cause the final bounds in
Theorem \ref{main} to grow extremely rapidly in $E$, although they will still of course be finite
for each $E$.
\end{remark}

A related result involves persistence of $L^2$ or $\dot H^2$
regularity:

\begin{lemma}[Persistence of regularity]\label{persistence}  Let $k = 0,1,2$, $I$ be a compact
time interval,
and let $u$ be a finite energy solution to \eqref{nls} on $I \times
\R^3$ obeying the bounds
$$ \| u \|_{L^{10}_{t,x}(I \times \R^3)} \leq M.$$
Then, if $t_0 \in I$ and $u(t_0) \in H^k$, we have
\begin{align}
\label{newark0}
 \| u \|_{\dot S^k(I \times \R^3)} &\leq C(M,E(u))
\| u(t_0) \|_{\dot H^k}.
\end{align}
\end{lemma}

In particular, once we control the $L^{10}_{t,x}$ norm of a finite
energy solution, we in fact control all the Strichartz norms in $\dot
S^1$, and can even control the $\dot S^2$ norm if the initial data is in $H^2(\R^3)$.
From this and standard iteration arguments, one can in fact show that a Schwartz solution
can be continued in time as long as the $L^{10}_{t,x}$ norm does not blow up to infinity.

\begin{proof}
By the local well-posedness theory it suffices to
prove \eqref{newark0} as an \emph{a priori} bound.

Applying Lemma
\ref{perturb} with $\tilde{u} := u$, $e := 0$, and $E' := 0$ gives us
the bound
\begin{align}
\label{newark1} \|u\|_{\dot{S}^1(I \times \R^3)} & \lesssim
C(M,E).
\end{align}
By \eqref{12} we also have
\begin{align}
\label{newark2} \| \nabla^k  \O(u^5)
\|_{L^1_t L^2_x} & \lesssim {{\|u \|}_{L^{10}_{x,t}}} \|u\|_{\dot
S^k}  \| u\|_{\dot S^1}^3;
\end{align}
the main thing to observe here is the presence of one factor of
$\|u\|_{L^{10}_{x,t}}$ on the right hand side.

As in the proof of Lemma \ref{perturb}, divide the time interval
$I$ into $N \approx \left(1 + \frac{M}{\delta} \right)^{10}$
subintervals $I_j := [T_j, T_j+1]$ on which we have
\begin{align}
\label{newark4} \| u \|_{L^{10}_{x,t} (I_j \times \R^3)} & \leq
\delta
\end{align}
where $\delta$ will be chosen momentarily.  We have on each $I_j$
by the Strichartz estimates (Lemma \ref{disjointed-strichartz})
and \eqref{newark2},
\begin{align*}
\| u \|_{\dot S^k (I_j \times \R^3)} & \leq C \left( \|
u(T_j)\|_{\dot H^k(\R^3)} + \| \nabla^k (|u|^4 u)
\|_{L^1_tL^2_x(I_j \times \R^3)} \right) \\
& \leq C \left(\| u(T_j)\|_{\dot H^k(\R^3)} + \|
u\|_{L^{10}_{x,t}(I_j \times \R^3)} \cdot \| u \|_{\dot S^k(I
\times \R^3)} \cdot \| u \|_{\dot S^1(I_j \times \R^3)}^3 \right).
\end{align*}
Choosing $\delta \leq (2 C  C(M,E))^{-3}$ gives,
\begin{align} \label{newarklast}
\| u\|_{\dot S^k(I_j \times \R^3)} & \leq 2 C \| u
(T_j)\|_{\dot{H}^k(\R^3)}.
\end{align}
The bound \eqref{newark0} now follows by adding up the bounds
\eqref{newarklast} we have on each subinterval.
\end{proof}

\section{Overview of proof of global spacetime
bounds}\label{overview-sec}

We now outline the proof of Theorem \ref{main},
breaking it down into a number of smaller Propositions.

\subsection{Zeroth stage: Induction on energy}

We say that a solution $u$ to \eqref{nls} is \emph{Schwartz} on a
slab $I \times \R^3$ if $u(t)$ is a Schwartz function for all $t \in
I$; note that such solutions are then also smooth in time as well as space,
thanks to \eqref{nls}.

The first observation is that in order to prove Theorem \ref{main},
it suffices to do so for Schwartz solutions.  Indeed, once one obtains a
uniform $L^{10}_{t,x}(I \times \R^3)$ bound for all Schwartz solutions
and all compact $I$, one can then
approximate arbitrary finite energy initial data by Schwartz initial data and use
Lemma \ref{perturb} to show that the corresponding sequence of solutions to \eqref{nls}
converge in $\dot S^1(I \times \R^3)$ to a finite energy solution to \eqref{nls}.  We omit
the standard details.

For every energy $E \geq 0$ we define the quantity $0 \leq M(E) \leq
+\infty$ by
$$ M(E) := \sup \{ \|u\|_{L^{10}_{t,x}(I_* \times \R^3)} \}$$
where $I_* \subset \R$ ranges over all compact time intervals, and
$u$ ranges over all Schwartz solutions to \eqref{nls} on $I_* \times
\R^3$ with $E(u) \leq E$.  We shall adopt the convention that $M(E) = 0$ for $E < 0$.
By the above discussion, it suffices to show that $M(E)$ is finite
for all $E$.

In the argument of Bourgain \cite{borg:scatter} (see also \cite{borg:book}),
the finiteness of $M(E)$ in the spherically symmetric case
is obtained by an induction on the energy $E$; indeed a bound of the form
$$ M(E) \leq C(E, \eta, M(E - \eta^4))$$
is obtained for some explicit $0 < \eta = \eta(E) \ll 1$ which does not collapse to $0$
for any finite $E$, and this easily implies via induction that $M(E)$
is  finite for all $E$.
Our argument will follow a similar induction on energy strategy,
however it will be convenient to run this induction in the contrapositive, assuming for
contradiction that $M(E)$ can be infinite, studying the minimal energy $E_{crit}$ for which
this is true, and then obtaining a contradiction using the ``induction hypothesis'' that $M(E)$ is
finite for all $E < E_{crit}$.  This will be more convenient for us, especially as we will require
more than one small parameter $\eta$.

We turn to the details.  We assume for contradiction that $M(E)$ is not always finite.
From Lemma \ref{perturb} we see that the set $\{ E: M(E) < \infty \}$
is open; clearly it is also connected and contains 0.  By our contradiction hypothesis, there must
therefore exist a \emph{critical energy} $0 < E_{crit} < \infty$ such that $M(E_{crit}) =
+\infty$,
but $M(E) < \infty$ for all $E < E_{crit}$.  One can think of $E_{crit}$ as the
minimal energy required to create a blowup solution.  For instance, we have

\begin{lemma}[Induction on energy hypothesis]\label{induction}  Let $t_0 \in \R$, and let $v(t_0)$ be a Schwartz function
such that $E(v(t_0)) \leq E_{crit} - \eta$ for some $\eta > 0$.  Then there exists a Schwartz global
solution $v: \R_t \times \R_x^3 \to \C$ to \eqref{nls} with initial data $v(t_0)$ at time $t=t_0$
such that $\| v \|_{L^{10}_{t,x}(\R \times \R^3)} \leq M(E_{crit}-\eta) = C(\eta)$.
Furthermore we have $\| v \|_{\dot S^1(\R \times \R^3)} \leq C(\eta)$.
\end{lemma}

Indeed, this Lemma follows immediately from the definition of $E_{crit}$, the local well-posedness theory in $L^{10}_{t,x}$,
and Lemma \ref{persistence}.

As in the argument in \cite{borg:scatter}, we will need a small parameter $0 < \eta = \eta(E_{crit}) \ll 1$
depending on $E_{crit}$.  In fact, our argument is somewhat lengthy and  we will
actually use  \emph{seven} such parameters
$$ 1 \gg \eta_0 \gg \eta_1 \gg \eta_2 \gg \eta_3 \gg \eta_4 \gg \eta_5 \gg \eta_6 > 0.$$
Specifically, we will need a small quantity $0 < \eta_0 = \eta_0(E_{crit})
\ll 1$ assumed to be sufficiently small depending on $E_{crit}$.  Then we need a smaller quantity
$0 < \eta_1 = \eta_1(\eta_0, E_{crit}) \ll 1$ assumed sufficiently small depending on $E_{crit}$, $\eta_0$
(in particular, it may be chosen smaller than positive quantities such as $M(E_{crit} - \eta_0^{100})^{-1}$).  We continue
in this fashion, choosing each $0 < \eta_j \ll 1$ to be sufficiently small depending on all previous
quantities $\eta_0, \ldots, \eta_{j-1}$ and the energy $E_{crit}$, all the way down to $\eta_6$ which is extremely small,
much smaller than any quantity depending on $E_{crit}, \eta_0, \ldots, \eta_5$ that will appear in our argument.
We will always assume implicitly that each $\eta_j$ has been chosen to be sufficiently small depending
on the previous parameters.  We  will often display the dependence of constants
on a parameter, e.g. $C(\eta)$ denotes a large constant
depending on $\eta$, and $c(\eta)$ will denote a small constant
depending upon $\eta$.  When $\eta_1 \gg \eta_2$, we will understand
$c(\eta_1) \gg c(\eta_2)$ and $C(\eta_1) \ll C(\eta_2)$.

Since $M(E_{crit})$ is infinite, it is in particular larger than $1/\eta_6$.  By definition of $M$, this means
that we may find a compact interval $I_* \subset \R$ and a smooth solution $u: I_* \times \R^3 \to \C$ to \eqref{nls}
with $E_{crit}/2 \leq E(u) \leq  E_{crit}$ so that
 $u$ is ridiculously large in the
sense that
\begin{equation}\label{m-bound}
\| u \|_{L^{10}_{t,x}(I_* \times \R^3)} > 1/\eta_6.
\end{equation}
We will show that this leads to a contradiction\footnote{Assuming, of course, that the parameters $\eta_0, \ldots, \eta_6$ are each chosen to be sufficiently small depending on previous parameters.  It is important to note however that the $\eta_j$
cannot be chosen to be small depending on the interval $I_*$ or the solution $u$; our estimates must be uniform with
respect to these parameters.}.  Although $u$ does not actually blow up (it is assumed smooth on all of the compact
interval $I_*$), it is still convenient to think of $u$ as almost\footnote{For instance, $u$ might genuinely blow up at some time $T_* > 0$, but $I_*$ is of the form $I_* = [0, T_* - \eps]$ for some very small $0 < \eps \ll 1$,
and thus $u$ remains Schwartz on $I_* \times \R^3$.} blowing up in $L^{10}_{t,x}$ in the sense
of \eqref{m-bound}.  We summarize the above discussion with the following,

\begin{definition} \label{minimal} A  {\it{minimal energy blowup solution}}
of \eqref{nls} is  a Schwartz solution on a time interval $I_*$ with energy\footnote{We could
modify our arguments below  to allow the assumption here $E(u) = E_{crit}$.  For example, the
arguments in the proof of Proposition \ref{frequency-decoupling} below also show that the function
$\tilde M(s) := \sup_{E(u)=s} \{\|u\|_{L^{10}_{x,t}}\}$ is a nondecreasing function of $s$. On first
reading, the reader may imagine $E(u) = E_{crit}$ in Definition \ref{minimal}.}, 
\begin{equation}
\label{energycritical}
\frac{1}{2} E_{crit} \leq
E(u)(t) = \int \frac{1}{2} |\nabla u(t,x)|^2 + \frac{1}{6}
|u(t,x)|^6\ dx \leq E_{crit}
\end{equation}
and $L^{10}_{x,t}$ norm enormous in the
sense of  \eqref{m-bound}.

\end{definition}

We remark that both conditions \eqref{m-bound}, \eqref{energycritical}
are invariant under the
scaling \eqref{scaling} (though of course the interval $I_*$ will be dilated by
$\lambda^2$ under this scaling).  Thus applying the scaling \eqref{scaling} to a minimal energy blowup solution produces another
minimal energy blowup solution.  Some of the proofs of the sub-propositions below will revolve around a specific frequency $N$; using this
scale invariance, we can then normalize that frequency to equal $1$ for the duration of that proof.  (Different parts of the argument involve different key frequencies, but we will not run into problems because we will only normalize one frequency at a time).

Henceforth we will not mention the $E_{crit}$ dependence of our constants
explicitly, as all our constants will depend on $E_{crit}$.  We shall need
however to keep careful track of the dependence of our argument on
$\eta_0, \ldots, \eta_6$.
Broadly speaking, we will start with the largest $\eta$, namely $\eta_0$, and
slowly ``retreat'' to increasingly smaller values of $\eta$ as the argument progresses
(such a retreat will for instance usually be required whenever the induction hypothesis Lemma \ref{induction}
is invoked).  However we will only retreat as far as $\eta_5$, not $\eta_6$, so that \eqref{m-bound}
will eventually lead to a contradiction when we show that
\begin{equation*}
  {{\| u \|}_{L^{10}_{t,x} (I^* \times \R^3 )}} \leq C( \eta_0, \dots
  , \eta_5).
\end{equation*}

Together with our assumption that we are considering a minimal
energy blowup solution $u$ as in Definition \ref{minimal}, Sobolev embedding
implies the bounds on kinetic energy
\begin{equation}\label{h1}
\| u \|_{L^\infty_t \dot H^1_x(I_* \times \R^3)} \sim 1
\end{equation}
and potential energy
\begin{equation}\label{l6}
\| u \|_{L^\infty_t L^6_x(I_* \times \R^3)} \lesssim 1
\end{equation}
(since our implicit constants are allowed to depend on $E_{crit}$).  Note that we do not
presently have any \emph{lower} bounds on the potential energy, but see below.


Having displayed our preliminary bounds on the kinetic and potential energy, we briefly discuss
the mass $\int_{\R^3} |u(t,x)|^2\ dx$, which is another conserved quantity.  Because of our \emph{a priori}
assumption that $u$ is Schwartz, we know that this mass is finite.  However, we cannot obtain uniform control
on this mass using our bounded energy assumption, because the very low frequencies of $u$ may simultaneously have
very small energy and very large mass.  Furthermore it is dangerous to rely too much on this conserved mass for this
energy-critical problem as the mass is not invariant under the natural scaling \eqref{scaling} of the equation
(indeed, it is super-critical with respect to that scaling).
On the other hand, from \eqref{h1} and \eqref{nabla-highfreq} we know that the \emph{high frequencies} of $u$ have
small mass:
\begin{equation}
  \label{nothigh}
  {{\| P_{>M} u \|}_{L^2 (\R^3)}} \lesssim \frac{1}{M} \hbox{ for all } M \in 2^\Z.
\end{equation}
Thus we will still be able to use the concept of mass in our estimates as long as we restrict our
attention to sufficiently high frequencies.

\subsection{First stage: Localization control on $u$}


We aim to show that a minimal energy blowup solution as in
Definition \ref{minimal} does not exist.  Intuitively, it seems reasonable to expect that
a minimal-energy blowup solution should be ``irreducible'' in the sense that it cannot
be decoupled into two or more components of strictly smaller energy that essentially do not interact with each other (i.e. each component also evolves via \eqref{nls} modulo small errors), since one of the components must then also blow up, contradicting the minimal-energy hypothesis.
In particular, we expect at every time that such a solution should be localized in both
frequency and space.

The first main step in the proof of Theorem \ref{main} is to make the above heuristics rigorous for our solution $u$.
Roughly speaking,
we would like to assert that at each time $t$, the solution $u(t)$ is localized
in both space and frequency to the maximum extent allowable under the uncertainty
principle (i.e. if the frequency is localized to $N(t)$, we would like to localize
$u(t)$ spatially to the scale $1/N(t)$).

These sorts of localizations already
 appear for instance in the argument of Bourgain
\cite{borg:scatter}, \cite{borg:book}, where the induction on energy argument
is introduced.  Informally\footnote{The heuristic that minimal energy blowup solutions should be strongly localized in both space and frequency has been employed in previous literature for a wide variety of nonlinear equations, including many of elliptic or parabolic type. Our formalizations of this heuristic, however, rely on the induction on energy methods of Bourgain and perturbation theory, as opposed to variational or compactness arguments.}
, the reason that we can expect such localization
is as follows.  Suppose for contradiction that at some time $t_0$ the solution $u(t_0)$
can be split into two parts $u(t_0) = v(t_0) + w(t_0)$ which are widely separated in either
space or frequency, and which each carry a nontrivial amount $O(\eta^C)$ of energy
for some $\eta_5 \leq \eta \leq \eta_0$.  Then by orthogonality we expect $v$ and $w$ to each have
strictly smaller energy than $u$, e.g. $E(v(t_0)), E(w(t_0)) \leq E_{crit} - O(\eta^C)$.
Thus by Lemma \ref{induction}
we can extend $v(t)$ and $w(t)$ to all of $I_* \times \R^3$ by evolving the nonlinear
Schr\"odinger equation \eqref{nls} for $v$ and $w$ separately, and furthermore we have
the bounds
$$ \| v \|_{L^{10}_{t,x}(I_* \times \R^3)}, \| w \|_{L^{10}_{t,x}(I_* \times \R^3)}
\leq M(E_{crit} - O(\eta^C)) \leq C(\eta).$$
Since $v$ and $w$ both solve \eqref{nls} separately, and $v$ and $w$ were assumed to be widely
separated, we thus expect $v+w$ to solve \eqref{nls} \emph{approximately}.  The idea
is then to use the perturbation theory from Section \ref{perturb-sec}
to  obtain a bound of the form
$\| u \|_{L^{10}_{t,x}(I_* \times \R^3)} \leq C(\eta)$, which contradicts \eqref{m-bound} if $\eta_6$
is sufficiently small.

A model example of this type of strategy occurs in Bourgain's
argument \cite{borg:book}, where substantial effort is invested in
locating a ``bubble'' - a small localized pocket of energy - which is
sufficiently isolated in physical space from the rest of the
solution.  One then removes this bubble, evolves the remainder of the
solution, and then uses perturbation theory, augmented with the
additional information about the isolation of the bubble, to place
the bubble back in.  We will use arguments similar to these in the sequel,
but first we need instead to show that a solution of \eqref{nls} which
is sufficiently delocalized in frequency space is globally spacetime bounded.
More precisely, we have:


\begin{proposition}[Frequency delocalization implies spacetime bound]\label{frequency-decoupling}
Let $\eta > 0$, and suppose there exists a dyadic
frequency $ N_{lo} > 0 $ and a
time $t_0 \in I_*$ such that we have the energy separation conditions
\begin{equation}\label{lo-energy}
\| P_{\leq N_{lo}} u(t_0) \|_{\dot H^1(\R^3)} \geq \eta
\end{equation}
and
\begin{equation}\label{hi-energy}
\| P_{\geq K(\eta) N_{lo}} u(t_0) \|_{\dot H^1(\R^3)} \geq \eta.
\end{equation}
If $K(\eta)$ is sufficiently large depending on $\eta$, i.e.
$$ K(\eta) \geq C(\eta)$$
then we have
\begin{equation}\label{10-bound}
\| u \|_{L^{10}_{t,x}(I_* \times \R^3)} \leq C(\eta).
\end{equation}
\end{proposition}

We prove this in Section \ref{decoupled-sec}.  The basic idea is
as outlined in previous discussion; the main technical tool we
need is the multilinear improvements to
Strichartz' inequality in Section \ref{quint-sec}
to control the interaction between the two
components and thus allow one to reconstruct the original solution
$u$.

Clearly the conclusion of Proposition \ref{frequency-decoupling} is in conflict with
the hypothesis \eqref{m-bound}, and so we should now expect the solution to be localized
in frequency for every time $t$.
This is indeed the case:

\begin{corollary}[Frequency localization of
energy at each time]\label{no-freq-dispersion}  A minimal energy blowup solution of \eqref{nls}  (see Definition \ref{minimal}) satisfies: For every time $t \in I_*$ there
exists a dyadic frequency $N(t) \in 2^{\Z}$ such that for every $\eta_5 \leq \eta \leq \eta_0$
we have small energy at frequencies $\ll N(t)$,
\begin{equation}\label{lowfreq-smallenergy}
\| P_{\leq c(\eta) N(t)} u(t) \|_{\dot H^1} \leq \eta,
\end{equation}
small energy at frequencies $\gg N(t)$,
\begin{equation}\label{highfreq-smallenergy}
\| P_{\geq C(\eta) N(t)} u(t) \|_{\dot H^1} \leq \eta,
\end{equation}
and large energy at frequencies $\sim N(t)$,
\begin{equation}\label{medfreq-largenergy}
\| P_{c(\eta) N(t) < \cdot < C(\eta) N(t)} u(t) \|_{\dot H^1} \sim 1.
\end{equation}
Here $0 < c(\eta) \ll 1 \ll C(\eta) < \infty$
are quantities depending on $\eta$.
\end{corollary}


Informally, this Corollary asserts that at every given time $t$ the solution $u$
is essentially concentrated at a single frequency $N(t)$.  Note
however that we do not presently have any information as to how
$N(t)$ evolves in time; obtaining long-term control on $N(t)$ will be a key objective of later
stages of the proof.

\begin{proof}
For each time $t \in I_*$, we define $N(t)$ as
$$ N(t) := \sup \{ N \in 2^\Z: \| P_{\leq N} u(t) \|_{\dot H^1} \leq \eta_0 \}.$$
Since $u(t)$ is Schwartz, we see that $N(t)$ is strictly larger than zero; from the lower bound in \eqref{h1}
we see that $N(t)$ is finite.   By definition of $N(t)$, we have
$$ \| P_{\leq 2N(t)} u(t) \|_{\dot H^1} > \eta_0.$$
Now let $\eta_5 \leq \eta \leq \eta_0$. Observe that we now have \eqref{highfreq-smallenergy} if
$C(\eta)$ is chosen sufficiently large, because if \eqref{highfreq-smallenergy} failed then
Proposition \ref{frequency-decoupling} would imply that $\| u \|_{L^{10}_{t,x}(I_* \times \R^3)}
\leq C(\eta)$, contradicting \eqref{m-bound} if $\eta_6$ is sufficiently small.  In particular we
have \eqref{highfreq-smallenergy} for $\eta = \eta_0$.  Since we also have
\eqref{lowfreq-smallenergy} for $\eta = \eta_0$ by construction of $N(t)$, we thus see from
\eqref{h1} that we have \eqref{medfreq-largenergy} for $\eta = \eta_0$, which of course then
implies (again by \eqref{h1}) the same bound for all $\eta_5 \leq \eta \leq \eta_0$.  Finally, we
obtain \eqref{lowfreq-smallenergy} for all $\eta_5 \leq \eta \leq \eta_0$ if $c(\eta)$ is chosen
sufficiently small, since if \eqref{lowfreq-smallenergy} failed then by combining it with
\eqref{medfreq-largenergy} and Proposition \ref{frequency-decoupling} we would once again imply
that $\| u \|_{L^{10}_{t,x}(I_* \times \R^3)} \leq C(\eta)$, contradicting \eqref{m-bound}.
\end{proof}

Having shown that any minimal energy blowup solution  $u$ must
be localized in frequency at each time, we now turn to showing that
such a $u$ is also localized in physical space.  This turns out to be somewhat more involved,
although it still follows the same general strategy.  We first borrow a useful trick
from \cite{borg:scatter}; since $u$ is Schwartz, we may divide the interval $I_*$ into
three consecutive pieces $I_* := I_- \cup I_0 \cup I_+$ where each of the three intervals contains a third of
the $L^{10}_{t,x}$ density:
$$ \int_I \int_{\R^3} |u(t,x)|^{10}\ dx dt = \frac{1}{3} \int_{I_*}\int_{\R^3} |u(t,x)|^{10}\ dx dt
\hbox{ for } I = I_-, I_0, I_+.$$
In particular from \eqref{m-bound} we have
\begin{equation}\label{m-bound-third}
 \| u \|_{L^{10}_{t,x}(I \times \R^3)} \gtrsim 1/\eta_6 \hbox{ for } I = I_-, I_0, I_+.
\end{equation}
Thus to contradict \eqref{m-bound} it suffices to obtain $L^{10}_{t,x}$ bounds on just one of the three
intervals $I_-$, $I_0$, $I_+$.

It is in the middle interval $I_0$ that we can obtain physical space localization; this shall be done in
several stages.  The first step is to ensure that the potential energy $\int_{\R^3} |u(t,x)|^6\ dx$ is
bounded from below.

\begin{proposition}[Potential energy bounded from below]\label{6-assume}
For any minimal energy blowup solution of \eqref{nls} (see Definition \ref{minimal})
we have for all $t \in I_0$,
\begin{equation}\label{6-large}
 \| u(t) \|_{L^6_x} \geq \eta_1.
\end{equation}
\end{proposition}

This is proven in Section \ref{six-small-sec}, and is inspired by a similar argument of Bourgain \cite{borg:scatter}.
Using \eqref{6-large} and some simple Fourier analysis, we can thus establish the following concentration result:

\begin{proposition}[Physical space concentration of energy at each time]\label{concentration} Any
minimal energy blowup solution of \eqref{nls} satisfies: For every $t \in I_0$, there
exists an $x(t) \in \R^3$ such that
\begin{equation}\label{kinetic-large}
 \int_{|x-x(t)| \leq C(\eta_1)/N(t)} |\nabla u(t,x)|^2\ dx \gtrsim c(\eta_1)
\end{equation}
and
\begin{equation}\label{lp-large}
\int_{|x-x(t)| \leq C(\eta_1)/N(t)} |u(t,x)|^p\ dx \gtrsim
c(\eta_1) N(t)^{\frac{p}{2} - 3}
\end{equation}
for all $1 < p < \infty$, where the implicit constant can depend on $p$.  In particular we have
\begin{equation}\label{potential-large}
 \int_{|x-x(t)| \leq C(\eta_1)/N(t)} |u(t,x)|^6\ dx \gtrsim c(\eta_1),
\end{equation}
\end{proposition}

This is proven in Section \ref{concentration-sec}.  Similar results were obtained in \cite{borg:scatter},
\cite{grillakis:scatter} in the radial case; see also \cite{BRefine}.  Informally, the above estimates
assert that $u(t,x)$ is roughly of size $N(t)^{1/2}$ on the average when $|x-x(t)| \lesssim 1/N(t)$; observe
that this is consistent with bounded energy \eqref{h1} as well as with Corollary \ref{no-freq-dispersion}
and the uncertainty principle.

It turns out that in our argument, it is not enough to know that the energy concentrates at one location $x(t)$ at
each time; we must also show that the energy is small at all other locations, where $|x-x(t)| \gg 1/N(t)$.
The main tool for achieving this is


\begin{proposition}[Physical space localization of energy at each time]\label{energy-decay}
For any minimal energy blowup solution of \eqref{nls} we have for every $t \in I_0$
\begin{equation}
\label{kinetic-localized}
 \int_{|x - x(t)| > 1 / (\eta_2 N(t))} |\nabla u(t,x)|^2\ dx \lesssim \eta_1.
\end{equation}
\end{proposition}


This is proven in Section \ref{energy-decay-sec}.  The proof follows a similar strategy to that used to
prove Proposition \ref{frequency-decoupling}; the main difference is that we now consider spatially
separated components of $u$ rather than frequency separated components, and
instead of using multilinear Strichartz estimates to establish the decoupling of these components,
we shall rely instead on approximate finite speed of propagation and on the pseudoconformal identity.

To summarize, at each time $t$ we have a location $x(t)$, around which the kinetic and potential
energy are large, and away from which the kinetic energy is small (and one can also show the
potential energy is small, although we will not need this).  From this and a little
Fourier analysis we obtain an important conclusion:

\begin{proposition}[Reverse Sobolev inequality]\label{reverse}
Assuming $u$ is a  minimal energy blowup solution (and hence \eqref{energycritical},
\eqref{lowfreq-smallenergy}-\eqref{kinetic-localized} hold), we have that
for every $t_0 \in I_0$, any $x_0 \in \R^3$, and any $R \geq 0$,
\begin{equation}\label{rs}
 \int_{B(x_0,R)} |\nabla u(t_0,x)|^2\ dx \lesssim \eta_1 + C(\eta_1,\eta_2)
\int_{B(x_0,{C}(\eta_1,\eta_2)R)} |u(t_0,x)|^6\ dx
\end{equation}
\end{proposition}

Thus, up to an error of $\eta_1$, we are able to control the kinetic
energy locally by the potential energy\footnote{Note that this is a special property of
the minimal energy blowup solution, reflecting the very strong
physical space localization properties of such
a solution; it is false in general, even for solutions to the free Schr\"odinger equation.  Of course,
Proposition \ref{6-assume} is similarly false in general, for instance for solutions of the free Schr\"odinger equation,
the $L^6_x$ norm goes to zero as $t \to \pm \infty$.}. This will be
proved in Section \ref{reverse-sec}.
This fact will be crucial in the interaction Morawetz portion of our argument when we have an error term
involving the kinetic energy, and control of a positive term which involves the potential energy; the reverse Sobolev
inequality is then used to control the former by the latter.


To summarize, the statements above tell us that any minimal energy blowup solution (Definition \ref{minimal}) to the equation
\eqref{nls} must be localized in both frequency and physical space
  at every time.
We are still far from done: we have not yet precluded blowup in finite time (which would happen if
$N(t) \to \infty$ as $t \to T_*$ for some finite time $T_*$), nor have we eliminated soliton or soliton-like
solutions (which would correspond, roughly speaking, to $N(t)$ staying close to constant for all time $t$).
To achieve this we need spacetime integrability bounds on $u$.  Our main tool for this is a frequency-localized
version of the interaction Morawetz estimate \eqref{morawetz-interaction}, to which we now turn.

\subsection{Second stage: Localized Morawetz estimate}

In order to localize the interaction Morawetz inequality, it turns out to be convenient to
work at the ``minimum'' frequency attained by $u$.

From \eqref{nabla-block} we observe that
$$ \| P_{c(\eta_0) N(t) < \cdot < C(\eta_0) N(t)} u(t) \|_{\dot H^1}
 \leq C(\eta_0) N(t) \| u \|_{L^\infty_t L^2_x}
$$
Comparing this with \eqref{medfreq-largenergy} we obtain the lower bound
$$ N(t) \geq c(\eta_0) \| u \|_{L^\infty_t L^2_x}^{-1}$$
for $t \in I_0$.  Since $u$ is Schwartz, the right-hand side is nonzero, and thus the quantity
$$ N_{min} := \inf_{t \in I_0} N(t)$$
is strictly positive.

From \eqref{lowfreq-smallenergy} we see that the low frequency
portion of the solution - where $|\xi| \leq c(\eta_0) N_{min}$ - has
small energy;
one might then hope to use Strichartz estimates to obtain some
spacetime control on these low frequencies.  However, we do not yet
have much control on the high
frequencies $|\xi| \geq c(\eta_0 ) N_{min}$, apart from the energy bounds
\eqref{h1} and
\eqref{l6} of
course.

Our initial spacetime bound in the high frequencies is provided by
the following
interaction Morawetz estimate.

\begin{proposition}[Frequency-localized interaction Morawetz
estimate]\label{Morawetz} Assuming $u$ is a  minimal energy blowup solution
of \eqref{nls} (and hence \eqref{energycritical},
\eqref{lowfreq-smallenergy}-\eqref{rs} all hold), we
have for all $N_* < c(\eta_3) N_{min}$
\begin{equation}\label{pm}
\int_{I_0} \int |P_{\geq N_*} u(t,x)|^4\ dx dt \lesssim \eta_1 N_*^{-3}.
\end{equation}
\end{proposition}

\begin{remark} The factor $N_*^{-3}$
on the right-hand side of \eqref{pm} is mandated by
scale-invariance considerations (cf. \eqref{scaling}). The $\eta_1$ factor
on the right side reflects our smallness assumption on $N_*$: if we think of
$N_*$ as being very small and then scale the solution so that $N_* = 1$, we
are pushing the energy to very high frequencies so heuristically it's
not unreasonable to expect the supercritical $L^4_{x,t}$ norm on the left hand
side to be small.

Regarding the size of $N_*$: write for
the  moment $\tilde c(\eta_3)$ as the constant
appearing in Corollary \ref{no-freq-dispersion} with
$\eta = \eta_3$. The constant $c(\eta_3)$ appearing in Proposition \ref{Morawetz}
is chosen so $c(\eta_3) \lesssim \tilde{c}(\eta_3) \cdot
\eta_3$,  hence at all times we know there is very little
energy at frequencies below $\frac{N_*}{\eta_3}$, and
(ignoring factors of $N_*$ which can be scaled
to 1) above frequency $N_*$ there is very little (at most $\eta_3/N_*$) $L^2$
mass.

This small $\eta_1$ factor will be used to close a
bootstrap argument in the proof of the important estimate on the movement of
energy to very low frequencies in Lemma
\ref{l2l1-bound} below. 
\end{remark}


 If one already had Theorem \ref{main}, then this Proposition would follow (but
with $\eta_1$ replaced by $C(E_{crit})$) from Lemma \ref{persistence},
since the $S^1$ norm will control $\| \nabla u \|_{L^4_t L^3_x}$
and hence $\||\nabla|^{3/4} u \|_{L^4_t L^4_x}$ by Sobolev
embedding. Of course, we will not prove Proposition \ref{Morawetz}
this way, as it would be circular. Instead, this Proposition is
based on the interaction Morawetz inequality developed in
\cite{ckstt:french}, \cite{ckstt:cubic-scatter} (see also a recent
extension in \cite{hassell}).  The key thing about this estimate
is that the right-hand side does not depend on $I_0$; thus for
instance it is already useful in eliminating soliton or
pseudosoliton solutions, at least for frequencies close to
$N_{min}$.  (Frequencies much larger than $N_{min}$ still cause
difficulty, and will be dealt with later in the argument).
Proposition \ref{Morawetz} roughly corresponds to the localized
Morawetz inequality used by Bourgain \cite{borg:scatter},
\cite{borg:book} and Grillakis \cite{grillakis:scatter} in the
radial case (see \eqref{morawetz-loc} above).  The main advantage
of \eqref{pm} is that it is not localized to near the spatial
origin, in contrast with the standard \eqref{morawetz} and
localized \eqref{morawetz-loc} Morawetz inequalities.

Although this Proposition is based on the interaction Morawetz inequality developed
in the references given above, there are significant technical difficulties in
truncating that inequality to the high frequencies.  As a consequence the proof of
this Proposition is somewhat involved and is given in Sections \ref{slimflim-sec}-\ref{flim-conclude}.
Also, we caution the reader that the above Proposition is not proved as an \emph{a priori}
estimate; indeed the proof relies crucially on the assumption that $u$ is a minimal energy blowup solution
in the sense of \eqref{m-bound}, and in particular verifies the reverse Sobolev
inequality \eqref{rs}.
See section \ref{slimflim-sec} for further remarks on the proof.

Combining Proposition \ref{Morawetz} with Proposition \ref{concentration} gives us the following
integral bound on $N(t)$.

\begin{corollary}\label{n-integral}  For any minimal energy
  blowup solution of \eqref{nls}, we have
\begin{equation}
\label{Nnotwild}
 \int_{I_0} N(t)^{-1}\ dt \lesssim C(\eta_1, \eta_3)
 N_{min}^{-3}.
\end{equation}
\end{corollary}


\begin{proof}
Let $N_* := c(\eta_3) N_{min}$ for
some sufficiently small $c(\eta_3)$.  Then from Proposition \ref{Morawetz}
we have
$$
 \int_{I_0} \int_{\R^3} |P_{\geq N_*} u(t,x)|^4\ dx dt \lesssim \eta_1 N_*^{-3}
\lesssim C(\eta_1,\eta_3) N_{min}^{-3}.
$$
On the other hand, from Bernstein \eqref{sobolev} and
\eqref{l6} we have for each $t \in I_0$ that
$$ \int_{|x-x(t)| \leq C(\eta_1)/N(t)} |P_{< N_*} u(t,x)|^4\ dx \lesssim
N(t)^{-3} \| P_{< N_*} u(t) \|_{L^\infty_x}^4 \lesssim C(\eta_1) N(t)^{-3} N_*^2,$$
so by \eqref{lp-large} and the triangle inequality we have (noting that
$N_* \leq c(\eta_3) N(t)$)
$$
 \int_{\R^3} |P_{\geq N_*} u(t,x)|^4\ dx \gtrsim c(\eta_1) N(t)^{-1}.
$$
Comparing this with the previous estimate, the claim follows.
\end{proof}

\begin{remark} The estimate \eqref{Nnotwild} is scale-invariant under the natural scaling \eqref{scaling}
($N$ has the units of $length^{-1}$, and $t$ has
the units of $length^2$).  In the radial case, a somewhat similar estimate was obtained
by Bourgain \cite{borg:scatter} and implicitly also by Grillakis \cite{grillakis:scatter}; in our notation,
this bound would be the assertion that
\begin{equation}\label{nt-bound}
\int_I N(t)\ dt \lesssim |I|^{1/2}
\end{equation}
for all $I \subseteq I_0$; indeed in the radial case (when
$x(t) = 0$) this bound easily follows from Proposition \ref{concentration} and \eqref{morawetz-loc}.
Both estimates are equally good at estimating the amount of time for which $N(t)$ is comparable to
$N_{min}$, but Corollary \ref{n-integral} is much weaker than \eqref{nt-bound} when it comes to controlling
the times for which $N(t) \gg N_{min}$.  Indeed if we could extend \eqref{nt-bound} to the nonradial case
one could obtain a significantly shorter proof of Theorem \ref{main}, however we were unable to prove this bound
directly\footnote{Comparing \eqref{pm} with \eqref{morawetz-loc} one also
sees that our control on how often the solution concentrates is weaker than
that in the radial arguments of Bourgain and Grillakis.  Heuristically:
\eqref{pm} allows
a long train in time of $N^3$ ``bubbles" at frequency $N \gg 1$, with size $\sim N^{\frac{1}{2}}$,
spatial extent $N^{-1}$, and individual duration $\sim N^{-2}$, so the total lifespan of the train is $\sim N$.  The
estimate \eqref{morawetz-loc}, on the other hand, restricts such bubbles to a set of dimension less
than $\frac{1}{2}$ in time.  Our proof of Theorem \ref{main} makes up for this
weakness in its next stage, specifically the relatively strong frequency localized
$L^2$ almost conservation estimate of Lemma \eqref{extralemma}.}, although it can be deduced from Corollary \ref{n-integral} and Proposition \ref{energy-travel} below).

\end{remark}

This Corollary allows us to obtain some useful $L^{10}_{t,x}$ bounds in the case when $N(t)$ is bounded from above.


\begin{corollary}[Nonconcentration implies spacetime bound]\label{bounded-snake-ok}
Let $I \subseteq I_0$, and suppose there exists an $N_{max} > 0$ such that $N(t) \leq N_{max}$ for all $t \in I$.
Then for any localized minimal energy blowup solution of \eqref{nls} we have
$$ \| u \|_{L^{10}_{t,x}(I \times \R^3)} \lesssim C( \eta_1, \eta_3, N_{max}/N_{min})$$
and furthermore
$$ \| u \|_{\dot S^1(I \times \R^3)} \lesssim C( \eta_1, \eta_3, N_{max}/N_{min}).$$
\end{corollary}


\begin{proof}  We may use scale invariance \eqref{scaling} to
rescale $N_{min} = 1$.  From Corollary \ref{n-integral} we obtain the useful bound
$$ |I| \lesssim C( \eta_1, \eta_3, N_{max}).$$
Let $\delta = \delta(\eta_0, N_{max}) > 0$ be a small number to be chosen later.  We
may partition $I$ into $O( |I|/\delta )$ intervals $I_1, \ldots, I_J$ of length at most $\delta$.  Let $I_j$ be any
of these intervals, and let $t_j$ be any time in $I_j$.  Observe from Corollary \ref{no-freq-dispersion}
and the hypothesis $N(t_j) \leq N_{max}$ that
$$ \| P_{\geq C(\eta_0) N_{max}} u(t_j) \|_{\dot H^1} \leq \eta_0$$
(for instance).  Now let $\tilde u(t) := e^{i(t-t_j)\Delta} P_{<C(\eta_0) N_{max}} u(t_j)$ be the free evolution of
the low and medium frequencies of $u$. The above estimate then becomes
$$ \| u(t_j) - \tilde u(t_j) \|_{\dot H^1} \leq \eta_0.$$
On the other hand, from Bernstein \eqref{sobolev} and \eqref{h1} we have
$$ \|\nabla \tilde u(t) \|_{L^{30/13}_x}
\lesssim C(\eta_0, N_{max}) \| \tilde u(t_j) \|_{\dot H^1} \lesssim C(\eta_0, N_{max})$$
for all $t \in I_j$, and hence
$$ \| \nabla \tilde u \|_{L^{10}_t L^{30/13}_x(I_j \times \R^3)} \lesssim
C(\eta_0, N_{max}) \delta^{1/10}.$$
Similarly we have
$$ \| \nabla (|\tilde u(t)|^4 \tilde u(t)) \|_{L^{6/5}_x}
\lesssim \| \nabla \tilde u(t) \|_{L^6_x} \| \tilde u(t) \|_{L^6_x}^4
\lesssim C(\eta_0,N_{max}) \| \tilde u(t_j) \|_{\dot H^1}^5 \lesssim C(\eta_0,N_{max})$$
and hence
$$ \| \nabla (|\tilde u(t)|^4 \tilde u(t)) \|_{L^2_t L^{6/5}_x(I_j \times \R^3)} \lesssim
C(\eta_0,N_{max})\delta^{1/2}.$$ From these two estimates, the energy bound \eqref{h1}, and Lemma
\ref{perturbtiny} with $e = - |\tilde{u}|^4 \tilde u$, we see (if $\delta$ is chosen sufficiently
small) that
$$ \| u \|_{L^{10}_{t,x}(I_j \times \R^3)} \lesssim 1$$
Summing this over each of the $O(|I|/\delta)$ intervals $I_j$ we obtain the desired $L^{10}_{t,x}$
bound.  The $\dot S^1$ bound then follows from Lemma \ref{persistence}.
\end{proof}

This above corollary gives the desired contradiction to \eqref{m-bound-third} when $N_{max}/N_{min}$ is bounded,
i.e. $N(t)$ stays in a bounded range.

\subsection{Third stage: Nonconcentration of energy}

Of course, any global well-posedness argument for \eqref{nls}
must eventually  exclude  a blowup scenario (self-similar or otherwise)
where $N(t)$ goes to infinity in finite time, and indeed by Corollary \ref{bounded-snake-ok}
this is the only remaining possibility for a minimal energy blowup solution.
Corollary \ref{no-freq-dispersion} implies that in such  a scenario the energy
must almost entirely evacuate the frequencies near $N_{min}$, and instead concentrate at frequencies  much larger than $N_{min}$.  While this scenario is consistent with
conservation of energy, it turns out to not be consistent with the time and frequency distribution of mass.

More specifically, we know there is a $t_{min} \in I_0$ so
that for all $t \in I_0, N(t) \geq N(t_{min}) := N_{min} > 0$.
By Corollary \ref{no-freq-dispersion}, at time $t_{min}$ the solution has the bulk of
its energy near the frequency $N_{min}$, and hence the medium frequencies at that time
have mass bounded below by,
\begin{align}
\label{late1}
\| P_{c(\eta_0) N_{min} \leq \cdot \leq C(\eta_0) N_{min}} u (t_{min}) \|_{L^2}
& \gtrsim c(\eta_0) N_{min}^{-1}.
\end{align}
The idea is to prove the following approximate mass conservation law
 for these high frequencies\footnote{It is necessary to truncate to
the high frequencies in order to exploit mass conservation because the low frequencies
contain an unbounded amount of mass.  This strategy of mollifying the solution in frequency space
in order to exploit a conservation law that would otherwise be unbounded or useless is inspired
by the ``$I$-method'' for sub-critical dispersive equations, see e.g. \cite{ckstt:cubic-scatter}.},
which states that while some mass might slip to very low frequencies as the solution moves to high frequencies, it can't all do so.

\begin{lemma}[Some mass freezes away from low frequencies] \label{extralemma}
Suppose $u$ is a minimal energy blowup solution of \eqref{nls}, and let $[t_{min}, t_{evac}] \subset I_0$ be such that
$N(t_{min}) = N_{min}$ and $N(t_{evac}) / N_{min} \geq C(\eta_5)$.
Then for all $t \in [t_{min}, t_{evac}]$,
\begin{align}
\label{late2}
\| P_{\geq \eta_4^{100} N_{min}} u(t) \|_{L^2} &  \gtrsim \eta_1 N_{min}^{-1}.
\end{align}
\end{lemma}


Lemma \ref{extralemma} will quickly show that the evacuation scenario - wherein
the solution cleanly concentrates energy to very high frequencies - cannot occur.
Instead the solution always leaves a nontrivial amount of mass and energy
behind at medium frequencies.  This ``littering'' of the solution
will serve (via Corollary \ref{no-freq-dispersion}) to keep $N(t)$ from escaping to
infinity\footnote{It is interesting to note that one must exploit conservation of energy, conservation of mass, \emph{and} conservation of momentum (via the Morawetz inequality) in order to prevent blowup for the equation \eqref{nls}; the same phenomenon occurs in the previous arguments \cite{borg:scatter}, \cite{grillakis:scatter} in the radial case, even though the details of those arguments are in many ways quite different to those here.} and gives us,
\begin{proposition}[Energy cannot evacuate from low frequencies]\label{energy-travel}
For any minimal energy
blowup solution of \eqref{nls} we have
\begin{equation}\label{delta}
N(t) \lesssim C(\eta_5) N_{min}
\end{equation}
for all $t \in I_0$.
\end{proposition}
We give the somewhat complicated proof of  Lemma \ref{extralemma} and
Proposition \ref{energy-travel} in Section \ref{evacuate-sec}.

By combining Proposition \ref{energy-travel} with Corollary
\ref{bounded-snake-ok}, we encounter a contradiction to
\eqref{m-bound-third} which completes the proof of Theorem \ref{main}.

The proofs of the above claims occupy the remainder of this paper.  Before
moving to these proofs, we summarize the role of the parameters $\eta_i,
i = 0, \ldots 5$ which have now all been introduced.  The number $\eta_1$
represents the amount of potential energy that must be present at every time
in a minimal energy blowup solution.  (Proposition
\ref{6-assume}); it also represents the extent of concentration of energy
(on the scale of $1/N(t)$) that must occur
in physical space at every time in a  minimal energy blowup solutions (Proposition \ref{concentration}).  The number $\eta_2$ is introduced in Proposition \ref{energy-decay},
where $1/\eta_2$ represents the extent
that there is localization (on the scale of $1/N(t)$) of energy
in a minimal energy blowup solution.
The number $\eta_3$  measures, on the
scale of the quantity $N_{min}$, what we mean by ``high frequency'' when we
say Proposition \ref{Morawetz} is an interaction Morawetz estimate localized
to high frequencies.  The number $\eta_4$ measures the frequency (on the scale
of $N_{min}$) below which the evolution can't move a certain portion (namely, $\eta_1$)
of the $L^2$ mass.  Finally, the number $\eta_0$ enters in Corollary  \ref{no-freq-dispersion}
and various other points in the paper where we simply use its value as a small, universal
constant.

\section{Frequency delocalized at one time $\implies$ spacetime
bounded}\label{decoupled-sec}

We now prove Proposition \ref{frequency-decoupling}.

Let $0 < \eps = \eps(\eta) \ll 1$ be a small number to be chosen
later.
 If $K(\eta)$ is sufficiently large depending
on $\eps$, then one can find at least $\eps^{-2}$ disjoint intervals
$[\eps^2 N_j, N_j/\eps^2]$, $j = 1, \ldots, \eps^{-2}$ contained inside $[N_{lo},
K(\eta) N_{lo}]$.  From \eqref{h1} and the pigeonhole principle, there must
therefore exist a $N_j$ so that the interval $[\eps^2 N_j, N_j/\eps^2]$
is mostly free of energy:
\begin{equation}\label{energy-free-zone}
\| P_{\eps^2 N_j \leq \cdot \leq N_j/\eps^2} u(t_0) \|_{\dot H^1}
\lesssim \eps.
\end{equation}
As the statement and conclusion of Proposition \ref{frequency-decoupling} is invariant
under the scaling \eqref{scaling} we may set $N_j := 1$.  Now define
$$ u_{lo}(t_0) := P_{\leq \eps} u(t_0); \quad u_{hi}(t_0) :=
P_{\geq 1/\eps} u(t_0).$$
The functions $u_{lo}(t_0)$, $u_{hi}(t_0)$ have strictly smaller energy than
$u$:

\begin{lemma}\label{energy-smaller}  If $\eps$ is sufficiently small
depending on $\eta$, then we have
$$ E(u_{lo}(t_0)), E(u_{hi}(t_0)) \leq E_{crit} - c\eta^C.$$
\end{lemma}

\begin{proof}
We prove this for $E(u_{lo}(t_0))$; the claim for $E(u_{hi}(t_0))$ is
similar.
Define $u_{hi'}(t_0) := P_{> \eps} u(t_0)$, so that $u(t_0) =
u_{lo}(t_0) + u_{hi'}(t_0)$, and consider the quantity
\begin{equation}\label{energy-gap}
 |E(u(t_0)) - E(u_{lo}(t_0)) - E(u_{hi'}(t_0))|.
\end{equation}
From \eqref{hamil-def} we can bound this by
\begin{equation}\label{energy-interaction}
\eqref{energy-gap}
\lesssim  |\langle \nabla u_{lo}(t_0), \nabla u_{hi'}(t_0) \rangle| + |\int
|u(t_0)|^6 - |u_{lo}(t_0)|^6 - |u_{hi'}(t_0)|^6\ dx|.
\end{equation}
The functions $u_{lo}$ and $u_{hi'}$ almost have disjoint supports, and their inner product is very close
to zero.  Indeed from Parseval and \eqref{energy-free-zone} we have
$$ |\langle \nabla u_{lo}(t_0), \nabla u_{hi'}(t_0) \rangle|  \lesssim
\eps^2.$$
Now we estimate the $L^6$-type terms.  From the pointwise estimate
$$ \left| |u(t_0)|^6 - |u_{lo}(t_0)|^6 - |u_{hi'}(t_0)|^6 \right|
\lesssim |u_{lo}|
|u_{hi'}| (|u_{lo}| + |u_{hi'}|)^4$$
(cf. \eqref{6-schematic}) and H\"older's inequality, we can bound the second term in
\eqref{energy-interaction} by
$$ \lesssim \| u_{lo} \|_\infty \| u_{hi'} \|_3 (\| u_{lo} \|_6 + \|
u_{hi'}\|_6)^4.$$
From \eqref{h1}, \eqref{energy-free-zone}, and Bernstein's inequality
\eqref{bernstein} we see that
\begin{align*}
\| u_{lo} \|_\infty
&\lesssim \sum_{N \leq \eps} \| P_N u \|_\infty \\
&\lesssim \sum_{N \leq \eps} N^{1/2} \| P_N u \|_{\dot H^1} \\
&\lesssim \sum_{N \leq \eps^2} N^{1/2} + \sum_{\eps^2 < N \leq
\eps } N^{1/2} \eps \\
&\lesssim \eps
\end{align*}
and
\bas
\| u_{hi'} \|_3
&\lesssim \sum_{N \geq \eps} \| P_N u \|_3 \\
&\lesssim \sum_{N \geq \eps} N^{-1/2} \| P_N u \|_{\dot H^1} \\
&\lesssim \sum_{N \geq 1/\eps} N^{-1/2} + \sum_{\eps \leq N \leq
\frac{1}{\eps}} N^{- 1/2} \eps \\
&\lesssim \eps^{1/2},
\end{align*}
and so from \eqref{l6} we can bound the second term in
\eqref{energy-interaction} by $O(\eps^{3/2})$.  Combining this with
the estimate obtained on the first piece of \eqref{energy-interaction}, we thus see that
$$ |E(u) - E(u_{lo}(t_0)) - E(u_{hi'}(t_0))| \lesssim \eps^{3/2}.$$
On the other hand, by hypothesis on $u$ we have $E(u) \leq E_{crit}$, while
from \eqref{hi-energy} and \eqref{hamil-def} we have $E(u_{hi'}(t_0))
\gtrsim \eta^C$.  The claim follows if $\eps$ is chosen sufficiently
small.
\end{proof}

From Lemma \ref{energy-smaller} and Lemma \ref{induction} we know that there exist Schwartz
solutions
$u_{lo}$, $u_{hi}$ on the slab $I_* \times \R^3$ with initial data
$u_{lo}(t_0)$, $u_{hi}(t_0)$ at time $t_0$, and furthermore
\begin{equation}\label{s1-bounds}
\| u_{lo} \|_{\dot S^1(I_* \times \R^3)},
\| u_{hi} \|_{\dot S^1(I_* \times \R^3)} \lesssim C(\eta).
\end{equation}
Let $\tilde u := u_{lo} + u_{hi}$.  We now claim that $\tilde u$ is
an approximate solution to \eqref{nls}:

\begin{lemma}\label{tildeu} We have
$$ i \tilde u_t + \Delta \tilde u = |\tilde u|^4 \tilde u - e$$
where the error $e$ obeys the bounds
\begin{equation}
\label{youwillbehave}
\| \nabla e \|_{L^2_t L^{6/5}_x(I_* \times \R^3)} \lesssim C(\eta)
\eps^{1/2}.
\end{equation}
\end{lemma}

\begin{proof} We begin by establishing further estimates on $u_{lo}$ and $u_{hi}$, beyond \eqref{s1-bounds}.  
For $u_{hi}$, we
observe from \eqref{nothigh} that $\| u_{hi}(t_0) \|_2 \lesssim \eps$,
so from Lemma \ref{persistence} we have
\begin{equation}\label{s0-bounds}
\| u_{hi} \|_{\dot S^0(I_* \times \R^3)} \lesssim C(\eta) \eps.
\end{equation}
Similarly, from \eqref{h1} and \eqref{nabla-lowfreq} we have
$\| u_{lo}(t_0) \|_{\dot H^2} \lesssim C \eps,$
and so from Lemma \ref{persistence} again we have
\begin{equation}\label{s2-bounds}
\| u_{lo}(t_0) \|_{\dot S^2(I_* \times \R^3)} \lesssim C(\eta) \eps.
\end{equation}

From Lemma \ref{leibnitz-holder} we also see that
\begin{align*}
\| |u_{hi}|^4 u_{hi} \|_{L^1_t L^2_x(I_* \times \R^3)} &\lesssim C(\eta) \eps\\
\| \nabla(|u_{hi}|^4 u_{hi}) \|_{L^1_t L^2_x(I_* \times \R^3)} &\lesssim C(\eta)\\
\| \nabla(|u_{lo}|^4 u_{lo}) \|_{L^1_t L^2_x(I_* \times \R^3)} &\lesssim C(\eta)\\
\| \nabla^2(|u_{lo}|^4 u_{lo}) \|_{L^1_t L^2_x(I_* \times \R^3)} &\lesssim C(\eta) \eps
\end{align*}

From Lemma \ref{quint-improved} we thus have
$$ \| \nabla \O(u_{hi}^j u_{lo}^{5-j}) \|_{L^2_t L^{6/5}_x(I_* \times \R^3)} \lesssim C(\eta)
\eps^{1/2}$$
for $j=1,2,3,4$.  Since $e = \sum_{j=1}^4 \O(u_{hi}^j u_{lo}^{5-j})$ by \eqref{5-schematic}, the claim follows.
\end{proof}

We now pass from estimates on $\tilde u$ to estimates on $u$ by
perturbation theory.  From \eqref{energy-free-zone} we have the
perturbation bound
$$ \| u(t_0 ) - \tilde u(t_0 ) \|_{\dot H^1} \lesssim \eps$$
while from \eqref{s1-bounds} we have
$$ \| \tilde u \|_{L^{10}_{t,x}(I_* \times \R^3)} \lesssim C(\eta).$$
Thus if $\eps$ is small enough depending on $\eta$, we may apply
Lemma \ref{perturb} and obtain the desired bound \eqref{10-bound}.  
The proof of Proposition \ref{frequency-decoupling} is now complete.
\endprf

\begin{remark}\label{huge-size} The dependence of constants in Proposition
\ref{frequency-decoupling} given by the above argument is quite poor.
Specifically, the separation $K(\eta)$ needs to be as large as
$$ K(\eta) \geq C \exp( C \eta^{-C} M(E_{crit} - \eta^C)^C )$$
(mainly in order for the pigeonhole argument to work) and the bound one obtains on the
$L^{10}_{t,x}$ norm at the end has a similar size.    This implies that the dependence of the
constants $C(\eta_j)$, $c(\eta_j)$ in Corollary \ref{no-freq-dispersion} is similarly similarly:
\begin{align*}
C(\eta_j)& \geq C \exp( C \eta_j^{-C} M(E_{crit} - \eta_j^C)^C ). \\
c(\eta_j) & \leq (C(\eta_j))^{-1}. \end{align*}
 This will force us to select each $\eta_{j+1}$
quite small depending on previous $\eta_j$; indeed in some cases the induction hypothesis is used
more than once and so $\eta_{j+1}$ is even smaller than the above expressions suggest.  If one
then runs the induction of energy argument in a direct way (rather than arguing by contradiction
as we do here), this leads to very rapidly growing (but still finite) bound for $M(E)$ for each
$E$, which can only be expressed in terms of multiply iterated towers of exponentials (the
Ackermann hierarchy). More precisely, if we use $X \uparrow Y$ to denote exponentiation $X^Y$,
$$X \uparrow \uparrow Y := X \uparrow (X \uparrow \ldots \uparrow X)$$
to denote the tower formed by exponentiating $Y$ copies of $X$,
$$X \uparrow \uparrow \uparrow Y := X \uparrow \uparrow (X \uparrow \uparrow \ldots \uparrow \uparrow X)$$
to denote the double tower formed by tower-exponentiating $Y$ copies of $X$, and so forth, then
we have computed our final bound for $M(E)$ for large $E$ to essentially be
$$M(E) \leq C \uparrow \uparrow \uparrow \uparrow \uparrow \uparrow \uparrow \uparrow (CE^C).$$
This rather Bunyanesque bound is mainly due to the large number of times
we invoke the induction hypothesis
Lemma \ref{induction}, and is presumably not best possible.  For instance, the best
bound known\footnote{Note added in proof: a bound of $M(E) \leq C \uparrow (C E^C)$ in the radial case
was recently obtained in \cite{tao}.} in the radial case is $M(E) \leq C \uparrow \uparrow (C E^C)$, where the induction hypothesis
is used only once; see \cite{borg:scatter}.  Finally, in the case of the subcritical cubic
nonlinear Schr\"odinger equation, the bound for the analogue of $M(E)$ is polynomial, $M(E) \leq C E^C$; see \cite{ckstt:cubic-scatter}. 
\end{remark}

\section{Small $L^6_x$ norm at one time $\implies$ spacetime bounded}\label{six-small-sec}

We now prove Proposition \ref{6-assume}.  The argument here is similar to the induction on energy arguments in
\cite{borg:scatter}.  The point is that the linear evolution of the solution must concentrate at some
point $(t_1,x_1)$ in spacetime (otherwise we could iterate using the small data theory).  If the solution does not concentrate in $L^6$ at time $t=t_0$, then $t_1$ must be far away from $t_0$.  The idea is then to remove the energy concentrating
at $(t_1,x_1)$ and induct on energy.

We turn to the details.  Assume for contradiction that \eqref{6-large} failed for some time $t_0 \in I_0$, so that
\begin{equation}\label{6-tiny}
 \| u(t_0) \|_{L^6_x} \leq \eta_1.
\end{equation}
Fix this $t_0$.  By rescaling using \eqref{scaling} we may normalize $N(t_0) = 1$.  Observe that if the linear solution
$e^{i(t-t_0)\Delta} u(t_0)$ had small $L^{10}_{t,x}$ norm, then the standard small data well-posedness
theory (based on Strichartz estimates and \eqref{12} or \eqref{26})
would already show that the nonlinear solution $u$ had bounded $L^{10}$ norm.  Thus we may assume that
$$ \| e^{i(t-t_0)\Delta} u(t_0) \|_{L^{10}_{t,x}(\R \times \R^3)} \gtrsim 1.$$
On the other hand, by Corollary \ref{no-freq-dispersion} we have
$$ \| P_{lo} u(t_0) \|_{\dot H^1_x}
+  \| P_{hi} u(t_0) \|_{\dot H^1_x} \lesssim \eta_0$$
where we define $P_{lo} := P_{<c(\eta_0)}$ and $P_{hi} := P_{>C(\eta_0)}$, so by Strichartz
(Lemma \ref{disjointed-strichartz})
$$ \| e^{i(t-t_0)\Delta} P_{lo} u(t_0) \|_{L^{10}_{t,x}(\R \times \R^3)}
+  \| e^{i(t-t_0)\Delta} P_{hi} u(t_0) \|_{L^{10}_{t,x}(\R \times \R^3)} \lesssim \eta_0$$
If we then define $P_{med} := 1 - P_{lo} - P_{hi}$, we must then have
$$ \| e^{i(t-t_0)\Delta} P_{med} u(t_0) \|_{L^{10}_{t,x}(\R \times \R^3)} \sim 1.$$
On the other hand, by \eqref{h1} we know that $P_{med} u(t_0)$ has bounded energy and has Fourier support in the region
$c(\eta_0) \lesssim |\xi| \lesssim C(\eta_0)$.  Thus by Strichartz \eqref{strichartz} we have that
$$ \| e^{i(t-t_0)\Delta} P_{med} u(t_0) \|_{L^{10/3}_{t,x}(\R \times \R^3)} \lesssim C(\eta_0)$$
(for instance).  From these two estimates and H\"older we see that the $L^\infty_{t,x}$ norm cannot be too small:
$$ \| e^{i(t-t_0)\Delta} P_{med} u(t_0) \|_{L^\infty_{t,x}(\R \times \R^3)} \gtrsim c(\eta_0).$$
In particular, there exists a time $t_1 \in \R$ and a point $x_1$ such that we have the concentration
$$ |e^{i(t_1-t_0)\Delta}(P_{med} u(t_0)) (x_1)| \gtrsim c(\eta_0).$$
By perturbing $t_1$ a little we may assume that $t_1 \neq t_0$; by time reversal symmetry we may take $t_1 < t_0$.

Let $\delta_{x_1}$ be the Dirac mass at $x_1$, and let $f(t_1) := P_{med} \delta_{x_1}$.  We extend $f$ to all of
$\R \times \R^3$ by the free evolution, thus $f(t) := e^{i(t-t_1)\Delta} f(t_1)$.  We record some explicit
estimates on $f$:

\begin{lemma}\label{f-est}  For any $t \in \R$ and any $1 \leq p \leq \infty$, we have
$$ \| f(t) \|_{L^p_x} \lesssim C(\eta_0) (1 + |t - t_1|)^{3/p - 3/2}.$$
\end{lemma}

\begin{proof}  We may translate $t_1 = x_1 = 0$.  From the unitarity of $e^{it\Delta}$ and Bernstein \eqref{bernstein}
we have
$$ \| f(t) \|_{L^\infty_x(\R^3)} \lesssim C(\eta_0)
\| f(t) \|_{L^2_x(\R^3)} = C(\eta_0) \| P_{med} \delta_{x_1} \|_{L^2_x(\R^3)} \lesssim C(\eta_0)$$
while from the dispersive inequality \eqref{dispersive} we have
$$ \| f(t) \|_{L^\infty_x(\R^3)} \lesssim |t|^{-3/2} \| P_{med} \delta_{x_1} \|_{L^1_x(\R^3)} \lesssim |t|^{-3/2}.$$
Combining these estimates we obtain the Lemma in the case $p = \infty$.  To obtain the other cases we need some
decay on $f(t,x)$ in the region $|x| \gg C(\eta_0) (1 + |t|)$.  For this we use the Fourier representation \eqref{fourier-rep}
to write
$$ f(t,x) = \int_{\R^3} e^{2\pi i (x \cdot \xi - 2\pi t |\xi|^2)} \varphi_{med}(\xi)\ d\xi$$
where $\varphi_{med}$ is the Fourier multiplier corresponding to $P_{med}$.  When $|x| \gg 1+|t|$, then
the phase oscillates in $\xi$ with a gradient comparable in magnitude to $|x|$.  By repeated integration by
parts (see e.g. \cite{stein:large}) we thus obtain a bound of the form $|f(t,x)| \lesssim |x|^{-100}$ in this
region.  Combining this with the previous $L^\infty$ bounds we obtain the result.
\end{proof}

In particular, from \eqref{6-tiny} and H\"older we have
$$ |\langle u(t_0), f(t_0) \rangle| \lesssim \eta_1 \| f(t_0) \|_{L^{6/5}_x(\R^3)}
\lesssim C(\eta_0) \eta_1 (1 + |t_0 - t_1|).$$
On the other hand, we have
$$ |\langle u(t_0), f(t_0) \rangle| = |\langle e^{i(t_1-t_0)\Delta} P_{med} u(t_0), \delta_{x_1} \rangle|
\gtrsim c(\eta_0).$$
Thus the concentration point $t_1$ must be far away from $t_0$ (recall that $\eta_1$ is much smaller
than $\eta_0$):
$$ |t_1 - t_0| \gtrsim c(\eta_0) \eta_1^{-1}.$$
In particular, the smallness of $\eta_1$ pushes the concentration time far away from the time
when $L^6_x$ is small.
Since $t_0 > t_1$ by hypothesis, we thus see from Lemma \ref{f-est} and the frequency localization of $f$
that $\nabla f$ has small $L^{10}_t L^{30/13}_x$ norm to the future of $t_0$:
\begin{equation}\label{f10-bound}
\begin{split}
\| \nabla f \|_{L^{10}_t L^{30/13}_x([t_0,+\infty) \times \R^3)} 
&\lesssim C(\eta_0) \| f \|_{L^{10}_t L^{30/13}([t_0,+\infty) \times \R^3)}\\
&\lesssim C(\eta_0) \|  (1 + |t - t_1|)^{-2/10} \|_{L^{10}_t([t_0,+\infty))}\\
&\lesssim C(\eta_0) | t_0 - t_1 |^{-1/10}\\
&\lesssim C(\eta_0) \eta_1^{1/10}.
\end{split}
\end{equation}
We now use the induction hypothesis (inspired by a similar argument in \cite{borg:scatter}).  We split
$u(t_0) := v(t_0) + w(t_0)$ where $w(t_0) := \delta e^{i\theta} \Delta^{-1} f(t_0)$, and $\delta = \delta(\eta_0) > 0$ is a small number to be chosen shortly, and $\theta$ is a phase to be chosen shortly.
We now claim that if $\delta$ and $\theta$ are chosen correctly, then
$v(t_0)$ has slightly smaller energy than $u$.  Indeed, we have by integration by parts and definition of $f$ that
\begin{align*}
\frac{1}{2} \int_{\R^3} |\nabla v(t_0)|^2 &=
\frac{1}{2} \int_{\R^3} |\nabla u(t_0) - \nabla w(t_0)|^2\\
&= \frac{1}{2} \int_{\R^3} |\nabla u(t_0)|^2 - \delta \Re
\int_{\R^3} e^{-i\theta} \overline{\nabla \Delta^{-1} f(t_0)} \cdot \nabla u(t_0)
+ O( \delta^2 \| \Delta^{-1} f(t_0) \|_{\dot H^1}^2 )\\
&\leq E_{crit} + \delta \Re e^{-i\theta} \langle u(t_0), f(t_0) \rangle
+ O(\delta^2 C(\eta_0)).
\end{align*}
Since $|\langle u(t_0), f(t_0) \rangle|$ was already shown to have magnitude at least $c(\eta_0)$,
we see if $\delta = \delta(\eta_0)$ and $\theta$ are chosen correctly that we can ensure that
$$ \frac{1}{2} \int_{\R^3} |\nabla v(t_0)|^2 \leq E_{crit} - c(\eta_0).$$
Meanwhile, another application of Lemma \ref{f-est} shows that
$$ \| w(t_0) \|_{L^6_x} \lesssim C(\eta_0)
\| f(t_0) \|_{L^{6}_x} \lesssim C(\eta_0) \langle t_0 - t_1 \rangle^{-1} \lesssim C(\eta_0) \eta_1$$
so by \eqref{6-tiny} and the triangle inequality we have
$$ \int_{\R^3} |v(t_0)|^6 \lesssim C(\eta_0) \eta_1^6.$$
Thus if $\eta_1$ is sufficiently small depending on $\eta_0$ we have
$$ E(v(t_0)) \leq E_{crit} - c(\eta_0).$$
By Lemma \ref{induction} we may extend $v(t_0) $ into a solution of
the nonlinear Schr\"odinger equation \eqref{nls}
on $[t_0,+\infty)$ such that
\begin{equation}\label{v10}
 \| v \|_{L^{10}_{t,x}([t_0,+\infty) \times \R^3)} \lesssim M(E_{crit} - c(\eta_0)) = C(\eta_0).
\end{equation}
On the other hand, from \eqref{f10-bound} and the frequency localization we have
$$ \| \nabla e^{i(t-t_0)\Delta} w(t_0) \|_{L^{10}_t L^{30/13}_x([t_0,+\infty) \times \R^3)} \lesssim C(\eta_0) \eta_1^{1/10}.$$
Thus if $\eta_1$ is sufficiently small depending on $\eta_0$ then we
may use Lemma \ref{perturb} (with $\tilde u := v$ and $e=0$) to conclude that $u$ extends to all of
$[t_0,+\infty)$ with
$$ \| u \|_{L^{10}_{t,x}([t_0,+\infty) \times \R^3)} \lesssim C(\eta_0, \eta_1).$$
Since this time interval $[t_0,+\infty)$ contains $I_+$, this contradicts \eqref{m-bound-third}.  (If
we had $t_0 < t_1$ instead, we would obtain a similar contradiction involving $I_-$).
This concludes the proof of Proposition \ref{6-assume}.
\endprf

\section{Spatial concentration of energy at every time}\label{concentration-sec}

We now prove Proposition \ref{concentration}.  Fix $t$.  By scaling using \eqref{scaling} we may take $N(t) = 1$.
By Corollary \ref{no-freq-dispersion} this implies that
\begin{equation}\label{conc-local}
\| P_{> C(\eta_1)} u(t) \|_{\dot H^1} + \| P_{<c(\eta_1)} u(t) \|_{\dot H^1}
\lesssim \eta_1^{100}
\end{equation}
(for instance).  In particular by Sobolev we have
$$\| P_{> C(\eta_1)} u(t) \|_{L^6_x} + \| P_{<c(\eta_1)} u(t) \|_{L^6_x}
\lesssim \eta_1^{100}$$
and hence by \eqref{6-large}
$$ \| P_{med} u(t) \|_{L^6_x} \gtrsim \eta_1$$
where $P_{med} := P_{c(\eta_1) < \cdot < C(\eta_1)}$.  On the other hand, from \eqref{h1}
we have
$$ \| P_{med} u(t) \|_{L^2_x} \lesssim C(\eta_1)$$
and thus by H\"older's inequality we have
$$ \| P_{med} u(t) \|_{L^\infty_x} \gtrsim c(\eta_1).$$
Thus there exists $x(t) \in \R^3$ such that
\begin{equation}
\label{xoftexists}
c(\eta_1) \lesssim |P_{med} u(t,x(t))|.
\end{equation}
Let $K_{med}$ denote the kernel associated to the operator
  $P_{med} \nabla \Delta^{-1}$, and let $R > 0$ be a radius to be chosen later. Then \eqref{xoftexists} can be continued with
  \begin{align*}
    c(\eta_1 ) & \lesssim | K_{med}
    * \nabla u (t, x(t)) | \\
& \lesssim \int |K_{med}(x(t) -x)| |\nabla u(t,x)| dx \\
& \thicksim \int_{|x - x(t)| < R )} |K_{med}(x(t) -x )|
|\nabla
    u(t,x)| dx +\int_{|x - x(t)| \geq R} |K_{med}(x(t) -x )| |\nabla
    u(t,x)| dx \\
& \lesssim C(\eta_1) {\left( \int_{|x - x(t)| < R}
|\nabla u(t,x)|^2 dx
    \right)^\half} + C(\eta_1) {\left( \int_{|x - x(t)| \geq R}
    \frac{|\nabla u (t,x)|}{|x - x(t)|^{100}} dx \right)}.
  \end{align*}
Here we used Cauchy-Schwarz and the fact that $K_{med}$ is a Schwartz function. Using the fact
that $(\int |\nabla u |^2 dx )^{1/2}$ is bounded uniformly, we obtain
\begin{equation*}
  c(\eta_1) \lesssim {\left( \int_{|x-x(t)| < R} |\nabla
  u(t,x)|^2 dx \right)}^{\half} + C(\eta_1) R^{-10}
\end{equation*}
(say), which proves \eqref{kinetic-large} after setting $R := C(\eta_1)$ sufficiently large. Similarly, writing $\tilde{K}_{med}$
for the kernel associated to $P_{med}$, we have for all $1 < p < \infty$,
\begin{align*}
  c(\eta_1 ) & \lesssim \int_{|x - x(t)| < R} |\tilde{K}_{med} (x(t) -x
  )| |u(t,x)| dx + \int_{|x - x(t)| \geq R} |\tilde{K}_{med} (x(t) -x
  )| |u(t,x)| dx \\
& \lesssim C(\eta_1) {\left( \int_{|x-x(t)| \leq R} |u(t,x)|^p
  dx \right)}^{\frac{1}{p}} + \int_{|x - x(t)| \geq R}
  \frac{|u(t,x)|}{{|x- x(t)|}^{100}} dx \\
& \lesssim C(\eta_1) {\left( \int_{|x - x(t)| \leq R}
  |u(t,x)|^p dx \right)}^{\frac{1}{p}} + {\left(
  \int_{|x-x(t)| > R}\frac{1}{{|x-x(t)|}^{100 \times \frac{6}{5}}} dx \right)}^{\frac{5}{6}}
  {{\| u(t) \|}_{L^6_x}} \\
& \lesssim C(\eta_1 ) {\left( \int_{|x-x(t)| \leq R} |u(t,x)|^p
  dx \right)}^{\frac{1}{p}} + C(\eta_1 ) R^{-10},
\end{align*}
where we used \eqref{l6}.  This proves
\eqref{lp-large} upon rescaling, if the radius $R = C(\eta_1)$ was chosen sufficiently large.
\endprf

\section{Spatial delocalized at one time $\implies$ spacetime bounded}\label{energy-decay-sec}


We now prove Proposition \ref{energy-decay}.  This is the spatial
analogue of the frequency delocalization result in Proposition
\ref{frequency-decoupling}.  The role of the bilinear Strichartz
estimate in that Proposition will be played here by finite speed
of propagation and pseudoconformal identity estimates. We follow
the same basic strategy as in Proposition
\ref{frequency-decoupling} (but now played out in the arena of
physical space rather than frequency space).  More specifically,
we assume that there is a large amount of energy away from the
concentration point, and then use approximate finite speed of
propagation to decouple the solution into two nearly
noninteracting components of strictly smaller energy, which can
then be handled by the induction hypothesis and perturbation
theory.

We turn to the details.  We first need a number of large
quantities.  Specifically, we need a large integer $J = J(\eta_1) \gg 1$
to be chosen later, and then a large
frequency\footnote{More precisely, this is a ratio of two
frequencies, but as we will shortly normalize $N(0) = 1$ the
distinction between a frequency and a frequency ratio becomes
irrelevant.  Similarly the radii given below should really be
ratios of radii.} $N_0 = N_0(\eta_1, J) \gg 1$ to be
chosen later, and then a large radius $R_0 = R_0(\eta_1,
N_0, J) \gg 1$ to be chosen later.

Suppose for contradiction that the Proposition were false. Then there must exist a time $t_0 \in
I_0$ such that
$$ \int_{|x - x(t_0)| > 1 / (\eta_2 N(t_0))} |\nabla u(t_0,x)|^2\ dx \gtrsim \eta_1.$$
Here $x(t)$ is the quantity constructed in Proposition \ref{concentration} (see \eqref{xoftexists}).
We may normalize $t_0 = x(t_0) = 0$, and rescale so that $N(0) =
1$, thus
\begin{equation}
\label{honebigoutRJ}
 \int_{|x| > 1/\eta_2} |\nabla u(0,x)|^2\ dx \gtrsim \eta_1.
\end{equation}
On the other hand, if $R_0 = R_0(\eta_1)$ is chosen large enough then we see
from Proposition \ref{concentration} that
\begin{equation}
\label{honebiginRzero}
 \int_{|x| < R_0} |\nabla u(0,x)|^2\ dx \gtrsim c(\eta_1)
\end{equation}
and
\begin{equation}
\label{lsixbiginRzero}
 \int_{|x| < R_0} |u(0,x)|^6\ dx \gtrsim c(\eta_1).
\end{equation}
We then define the radii $R_0 \ll R_1 \ll \ldots \ll R_J$
recursively by $R_{j+1} := 100 R_j^{100}$.

The region $R_0 < |x| < R_J$ can be partitioned into $J$ dyadic
shells of the form $R_j < |x| < R_{j+1}$.  By
\eqref{h1},\eqref{l6} and the pigeonhole principle we may find $0
\leq j < J$ such that
\begin{equation}\label{cut-small}
 \int_{R_j < |x| < R_{j+1}} |\nabla u(0,x)|^2 + |u(0,x)|^6\ dx \lesssim \frac{1}{J}.
\end{equation}
We fix this $j$.  We now introduce cutoff functions $\chi_{inner}$, $\chi_{outer}$, where
$\chi_{inner}$ is adapted to the ball $B(0, 2R_j)$ and equals one on $B(0,R_j)$, whereas
$\chi_{outer}$ is a bump function adapted to $B(0,R_{j+1})$ which equals one on $B(0,R_{j+1}/2)$.
We then define $v(0)$, $w(0)$ as
\begin{align} \label{herehere}
 v(0,x) & := P_{1/N_0 \leq \cdot \leq N_0} (\chi_{inner}
u(0)); \quad w(0) \; := \; P_{1/N_0 \leq \cdot \leq N_0} ((1-\chi_{outer}) u(0)).
\end{align}
By
\eqref{cut-small} we easily see that
$$ \| P_{1/N_0 \leq \cdot \leq N_0} u(0) - v(0) - w(0) \|_{\dot H^1}
\lesssim \| (\chi_{outer} - \chi_{inner}) u(0) \|_{\dot H^1}
\lesssim \frac{1}{J^{1/2}}.$$ Also, if $N_0$ is chosen sufficiently
large depending on $J$ we see from the normalization $N(0) = 1$
and Corollary \ref{no-freq-dispersion} that
$$ \| u(0) - P_{1/N_0 \leq \cdot \leq N_0} u(0) \|_{\dot H^1} \lesssim \frac{1}{J^{1/2}}$$
and thus
\begin{equation}\label{uvw}
 \| u(0) - v(0) - w(0) \|_{\dot H^1} \lesssim \frac{1}{J^{1/2}}.
\end{equation}
We also know that $v$, $w$ have slightly smaller energy than $u$:

\begin{lemma}  For the the functions  $v,w$ defined in \eqref{herehere} we have
$$ E(v(0)), E(w(0)) \leq E_{crit} - c(\eta_1).$$
\end{lemma}


\begin{proof} The argument here is completely analogous to that
in Lemma \ref{energy-smaller},  except that now
we work in physical space, exploiting the fact that $v$ is mostly supported in 
the region $|x| < 3R_j$ and $w$ is mostly supported in the region $|x| > R_{j+1}/2$.

We begin by estimating the quantity
$$ \left| E(v(0)+w(0)) - E(v(0)) - E(w(0)) \right|.$$
Expanding out the definition of energy, we can bound this by
$$ \lesssim \int_{\R^3} |\nabla v(0,x)| |\nabla w(0,x)|
+ |v(0,x)| |w(0,x)|^5 + |v(0,x)|^5 |w(0,x)|\ dx.$$
We subdivide $\R^3$ into the regions $|x| \leq R_{j+1}/2$ and $|x| > R_{j+1}/2$.  Since $v$ and $w$ are bounded
in $\dot H^1$ and hence in $L^6$, we can use H\"older's inequality to estimate the previous expression by
$$ \lesssim \| \nabla v \|_{L^2(|x| > R_{j+1}/2)} + \| v \|_{L^6(|x| > R_{j+1}/2)} + \| \nabla w \|_{L^2(|x| \leq R_{j+1}/2)}
 + \|w\|_{L^6(|x| \leq R_{j+1}/2)}.$$
Consider for instance the quantity $\|\nabla v \|_{L^2(|x| > R_{j+1}/2)}$.  Let $K$ be the convolution kernel associated with
$P_{1/N_0\leq \cdot \leq N_0}$, then we have $\nabla v = (\nabla (\chi_{inner} u(0))) * K$.  Since $\chi_{inner}$ is supported
on the region $|x| \leq 2R_j \leq R_{j+1}/4$, we may restrict $K$ to the region $|x| > R_{j+1}/4$.  On this region, $K$ decays rapidly
and in fact has an $L^1$ norm of at most $C(N_0) / R_{j+1}^{100} \leq C(N_0) / R_0^{100}$ (for instance).  
Since $\nabla (\chi_{inner} u(0))$ is bounded in $L^2$, we thus have
$$ \|\nabla v \|_{L^2(|x| > R_{j+1}/2)} \leq C(N_0) / R_0^{100}.$$
The other three terms above can be estimated similarly.  Thus we have
$$ \left| E(v(0)+w(0)) - E(v(0)) - E(w(0)) \right| \leq C(N_0) / R_0^{100}.$$
On the other hand, from \eqref{uvw}, the boundedness of $u(0), v(0), w(0)$ in $\dot H^1$ and $L^6$, and H\"older's inequality, we have
$$ \left| E(u(0)) - E(v(0)+w(0)) \right| \lesssim J^{-1/2}.$$
Thus we have
\begin{equation}\label{squeezinsqueezin}
 |E(u(0)) - E(v(0)) - E(w(0))| \lesssim J^{-1/2} + C(N_0) / R_0^{100}.
\end{equation}
On the other hand, by \eqref{honebiginRzero} and
\eqref{honebigoutRJ} (choosing $\eta_2$ to be smaller than $1/R_J$), we
know that $E(v)(0), E(w)(0) \gtrsim c(\eta_1)$. Together with \eqref{squeezinsqueezin},
this yields the Lemma when $J = J(\eta_1)$ and $R_0 = R_0(\eta_1,N_0,J)$ are chosen sufficiently large.
\end{proof}

From the above Lemma and Lemma \ref{induction} we may then extend $v$ and $w$ by the nonlinear Schr\"odinger equation
\eqref{nls} to all of $\R \times \R^3$, so that
\begin{equation}\label{vw-l10}
 \| v \|_{L^{10}_{t,x}} + \| w \|_{L^{10}_{t,x}} \lesssim M(E_{crit} - c(\eta_1))
= C(\eta_1).
\end{equation}
From this, Lemma \ref{persistence} and the frequency localization of $v,w$ we thus obtain
the Strichartz bounds
\begin{equation}
\label{inductstrich}
 \| v \|_{\dot S^k} + \| w \|_{\dot S^k} \lesssim C(\eta_1,
 N_0)
\end{equation}
for $k=0,1,2$.

The idea is now to use our perturbation lemma to approximate $u$ by $v+w$.  To do this we need to
ensure that $v$ and $w$ do not interact.  This is the objective of the next two lemmas.

\begin{lemma} \label{speed} Let $v(x,t), w(x,t)$ be the evolutions according to \eqref{nls} of the functions defined
in \eqref{herehere}.
  For times $|t| \leq R_j^{10}$, we have the ``finite speed of propagation'' estimate
\begin{equation}
\label{voutsidesmalll2}
 \int_{|x| \gtrsim R_j^{50}} |v(t,x)|^2\ dx \lesssim C(\eta_1,
 N_0) R_j^{-20}
\end{equation}
and for times $|t| \geq R_j^{10}$ we have the decay estimate
\begin{equation}
\label{vdecays} \int |v(t,x)|^6\ dx \lesssim R_j^{-10}.
\end{equation}
(The powers of $R_j$ are far from sharp).  Meanwhile, the mass
density of  $w$ obeys the finite speed of propagation estimate
\begin{equation}
\label{winsidesmalll2} \int_{|x| \lesssim R_j^{50}} |w(t,x)|^2\ dx
\lesssim C(\eta_1, N_0) R_j^{-20}\end{equation} for all
$|t| \leq R_j^{10}$ 
and, similarly,
the energy density of $w$ obeys the finite speed of propagation
estimate
\begin{equation}
  \label{winsidesmalle}
  \int_{|x|\lesssim R_j^{50}} [ \frac{1}{2}
|\nabla w (t,x)|^2 + \frac{1}{6} |w(t,x)|^6] dx
  \lesssim R_j^{-20} C(\eta_1, N_0) ,
\end{equation}
for all $|t| < R_j^{10}$. We also have
\begin{equation}
  \label{voutsidesmalle}
  \int_{|x|\gtrsim R_j^{50}} [\frac{1}{2} |\nabla v (t,x)|^2 + \frac{1}{6} |v(t,x)|^6] dx
  \lesssim R_j^{-20} C(\eta_1, N_0)
\end{equation}
\end{lemma}

Thus at short times $t = O(R_j^{10})$, $v$ and $w$ are separated
in space, whereas at long times $v$ has decayed (while $w$ is
still bounded in Strichartz norms).

\begin{proof}
The estimates \eqref{voutsidesmalll2} and \eqref{vdecays} follow
from the pseudoconformal law following arguments from
\cite{borg:book}. Recall the pseudoconformal conservation law for
sufficiently regular and decaying solutions of \eqref{nls}:
\begin{equation}
  \label{pclaw}
  {{\| (x + 2it\nabla) u(t)\|}_{L^2_x}^2} + \frac{4}{3} t^2
{{\| u(t) \|}_{L^6_x}^6} = {{\| |x| u_0 \|}_{L^2_x}^2} -
\frac{16}{3} \int_0^t s {{\| u(s)\|}^6_{L^6_x}} ds.
\end{equation}
Thus, since $v$ solves \eqref{nls}, we have that
\begin{equation*}
  \int_{|x| \gtrsim R_j^{50}} |x|^2 |v(t,x)|^2 dx \lesssim t^2 {{\|
  \nabla v(t) \|}_{L^2_x}^2} + t^2 {{\| v(t) \|}_{L^6_x}^6} + {{\| |x|
  v_0 \|}_{L^2_x}^2} + \int_0^t s {{\| v(s)\|}^6_{L^6_x}} ds.
\end{equation*}
We restrict to times $|t| \leq R_j^{10}$ and have
\begin{align*}
  R_j^{100}& \int_{|x|  \gtrsim R_j^{50}} |v(t,x)|^2 dx  \lesssim
  R_j^{20} {{\| \nabla v (t) \|}^2_{L^\infty_{|t| \leq R_j^{10}} L^2_x}}
  + R_j^{20}  {{\| v (t) \|}^6_{L^\infty_{|t| \leq R_j^{10}} L^6_x}} +
  R_j^2 {{\| v_0 \|}_{L^2_x}}^2 \\
& \lesssim R_j^{20} \left( {{\| \nabla v (t) \|}^2_{L^\infty_{|t|
\leq
  R_j^{10}} L^2_x}} + {{\| v (t) \|}^6_{L^\infty_{|t| \leq R_j^{10}}
  L^6_x}}  \right) + R_j^2 N_0^{2}  {{\| \nabla v_0 \|}_{L^2_x}^2} \\
& \lesssim C(N_0) R_j^{20} E(u_0),
\end{align*}
which proves \eqref{voutsidesmalll2}.

From \eqref{pclaw}, we observe that
\begin{equation*}
  {{\| v(t) \|}_{L^6_x}^6} \lesssim  \frac{R_j^2 N_0^{2} E (u_0)}{t^2},
\end{equation*}
so, for times $|t| > R_j^{10}$, we obtain \eqref{vdecays}.

We control the $L^2_x$-mass of $w$ in the ball $|x| < 1000
R_j^{50}$ using a virial identity. Let $\zeta$ denote a
nonnegative smooth bump function equaling 1 on $B(0,1000
R_j^{50})$ and supported on $B(0, 2000 R_j^{50})$. Note that
$\zeta$ has been chosen so that the support of $\nabla \zeta$ does
not intersect the support of $\chi_{inner}(0)$ or the support of $(1-\chi_{outer})(0)$. From
\eqref{local-mass-conserv}, \eqref{mass-cancel} and integration by
parts we have the identity
$$
\partial_t \int \zeta (x) |w(t,x)|^2 dx
 = - 2 \int \zeta_j (x) \Im  (w \overline{w}_j )(t,x) dx.
$$
Thus,
\begin{equation*}
  | \partial_t \int \zeta (x) |w(t,x)|^2 dx | \lesssim R_j^{-50} {{\|
    \nabla w (t) \|}_{L^2_x}} {{ \| w(t) \|}_{L^2_x}},
\end{equation*}
and we have, using the support properties and \eqref{h1},
\begin{align*}
  \sup_{|t| < R_j^{10} } \int \zeta (x) |w(t,x)|^2 dx & \lesssim \int
  \zeta(x) |w(0,x)|^2 dx + R_j^{-40} \sup_{|t| < R_j^{10} } {{\|
    \nabla w (t) \|}_{L^2_x}} {{ \| w(t) \|}_{L^2_x}}\\
&\lesssim R_j^{-40} + R_j^{-40} (E(u_0))^2 N_0
\end{align*}
(say), which proves \eqref{winsidesmalll2}.

We now control the energy density of $w$,
\begin{equation}
  \label{energydensity}
  e(w)(t,x) := \frac{1}{2} |\nabla w (t,x)|^2 + \frac{1}{6} |w(t,x)|^6,
\end{equation}
on the ball $|x| < R_j^{50}$ by a similar argument.  From
\eqref{energy-conserv} and integration by parts we have
\begin{equation*}
  \frac{d}{dt} \int \zeta e(w) dx = \int \zeta_j [ \Im
  (\overline{w}_k w_{kj}) - \delta_{jk} |w|^4 \Im (w \overline{w}_k
  ) ] dx,
\end{equation*}
which implies that
\begin{align*}
  \int \zeta  e(w)(T) dx \lesssim \int \zeta & e(w)(0) dx  + \int_0^T
  \int |\nabla \zeta | |\nabla w | |\nabla \nabla w| dx dt \\
& + \int_0^T \int |\nabla \zeta | |w|^5 |\nabla w | dx dt.
\end{align*}
We will control the three terms on the right to obtain
\eqref{winsidesmalle}. The first term vanishes due to support
properties of $\zeta$ and $1-\chi_{outer}$. The second term is crudely
bounded using \eqref{h1} by
\begin{equation*}
  R_j^{-50} R_j^{10} {{ \| \nabla^2 w \|}_{L^\infty_{|t| <
  R_j^{10}} L^2_x }}.
\end{equation*}
By the induction hypothesis we have \eqref{inductstrich} and, in
particular,
\begin{equation*}
  {{ \| \nabla^2 w \|}_{L^\infty_{t} L^2_x }} \lesssim C(\eta_1 , N_0 ).
\end{equation*}
The third term is bounded using H\"older by, say
\begin{equation*}
\lesssim  R_j^{-50} {{\| w \|}_{L^{10}_{t,x}}^2} {{\| w
\|}_{L^{6}_{t,x}}^3} {{\| \nabla w \|}_{L^{10/3}_{t,x}}}.
\end{equation*}
Again, the global $L^{10}_{t,x}$ bound and Lemma \ref{persistence}
gave us \eqref{inductstrich} which, by interpolation, controls all
the norms appearing here. This proves \eqref{winsidesmalle}.

Replacing $\zeta$ by $1-\zeta$ and $w$ by $v$ in the discussion
just completed establishes \eqref{voutsidesmalle}.
\end{proof}

\begin{corollary} For $v,w$ as in Lemma \ref{speed}, we have
$$ \| \nabla (|v+w|^4 (v+w) - |v|^4 v - |w|^4 w) \|_{L^2_t L^{6/5}_x(\R \times \R^3)} \lesssim C(\eta_1,N_0) R_j^{-5/6}
\lesssim C(\eta_1,N_0) R_0^{-5/6}.$$
\end{corollary}

\begin{proof}
By \eqref{5-schematic}, the task is to control terms of the form $\O (v^j w^{4-j} \nabla
w)$ and $\O (w^j v^{4-j} \nabla v)$, for $j = 1,2,3,4$, in $L^2_t L^{6/5}_x$. Separate the
analysis into three spacetime regions based on the estimates in Lemma \ref{speed}: (short time,
near origin) $|t| < R_j^{10}, ~|x| < 2000 R_j^{50}$; (short time, far from origin) $|t| <
R_j^{10}, ~|x| \geq 2000 R_j^{50}$; (long time) $|t| \geq R_j^{10}$. In all but one of these
cases, an application of a variant of \eqref{26},
\begin{equation*}
  {{\| \nabla u_1 u_2 u_3 u_4 u_5 \|}_{L^2_t L^{6/5}_x}} \lesssim
{{\| \nabla u_1 \|}_{L^\infty_t L^2_x }}{{\| u_2 \|}_{L^4_t
L^\infty_x
  }}{{\| u_3 \|}_{L^4_t L^\infty_x}}{{\| u_4 \|}_{L^\infty_t L^6_x
  }}{{\| u_5 \|}_{L^\infty_t L^6_x}},
\end{equation*}
together with \eqref{inductstrich} and
the decay properties in Lemma \ref{speed}
establishes the claimed estimate controlling the interaction of
$v$ and $w$. The term  $\O (w^4 \nabla v)$ in the long time regime
$|t| \geq R_j^{10}$ presents an exceptional case since we do not
directly encounter the available long time decay estimate
\eqref{vdecays}. This situation is treated separately with the
following argument. By H\"older and interpolation, we have
\begin{align*}
  {{\| \O ( w^4 \nabla v) \|}_{L^2_t L^{6/5}_x}} & \lesssim {{\| \nabla v \|}_{L^\infty_t L^3_x }} {{\| \O (w^4) \|}_{L^2_t L^2_x}} \\
& \lesssim {{\| v \|}_{L^\infty_t L^6_x }^{1/2}} {{\| \nabla^2 v \|}_{L^\infty_t L^2_x }^{1/2}}  {{\| \O (w^4) \|}_{L^2_t L^2_x}} \\
& \lesssim {{\| v \|}_{L^\infty_t L^6_x }^{1/2}} {{\|  v \|}_{\dot{S}^2 }^{1/2}}  {{\| \O (w^4) \|}_{L^2_t L^2_x}} \\
\end{align*}
and note the appearance of \eqref{vdecays} which contributes the
decay $R_j^{-5/6}$. We complete the proof by estimating,
\begin{align*}
  {{\| \O (w^4) \|}_{L^2_t L^2_x}} & \lesssim {{\| \O (w^3) \|}_{L^2_t L^6_x}} {{\| w \|}_{L^\infty_t L^3_x}} \\
& \lesssim {{\| w \|}_{L^6_t L^{18}_x}^3} {{\| w \|}_{L^\infty_t L^2_x }^{1/2}} {{\| w \|}_{L^\infty_t L^6_x}^{1/2}} \\
& \lesssim {{\| w \|}_{\dot{S}^1}^3}  {{\| w(0) \|}_{L^2_x}^{1/2}} E^{1/4}\\
& \lesssim C(\eta_1 , \eta_2 ) N_0^{1/2}\\
&\lesssim C(\eta_1 , \eta_2 , N_0 ).
\end{align*}
\end{proof}

In light of this corollary, \eqref{uvw}, \eqref{vw-l10}, and the observation
that $u,v,w$ all have bounded energy, we see from Lemma
\ref{perturb} (with $\tilde u := v+w$ and $e := |v+w|^4 (v+w) - |v|^4 v - |w|^4 w$)
that if $J$ is sufficiently
large depending on $\eta_1$, and $R_0$ sufficiently large
depending on $\eta_1, J, N_0$, then we have
$$ \| u \|_{L^{10}_{t,x}(I_* \times \R^3)} \lesssim C(\eta_1),$$
which contradicts \eqref{m-bound}. This proves Proposition
\ref{energy-decay}.
\endprf

\section{Reverse Sobolev inequality}\label{reverse-sec}

We now prove Proposition \ref{reverse}.  Fix $t_0, x_0, R$.  We may normalize $x(t_0) = 0$
and $N(t_0) = 1$.  Then by Proposition \ref{energy-decay} we have
\begin{equation} \label{decemberstar} \int_{|x| > 1/\eta_2} |\nabla u(t_0,x)|^2\ dx \lesssim \eta_1.
\end{equation}
Now suppose for contradiction that we had
\begin{equation}\label{rs-contradiction}
 \int_{B(x_0,R)} |\nabla u(t_0,x)|^2\ dx \gg \eta_1 + K(\eta_1,\eta_2)
\int_{B(x_0,K(\eta_1,\eta_2) R)} |u(t_0,x)|^6\ dx
\end{equation}
for some large $K(\eta_1,\eta_2)$ to be chosen later.  From \eqref{rs-contradiction}
and \eqref{decemberstar} we see that $B(x_0,R)$ cannot be completely contained inside
the region $|x| > 1/\eta_2$; so we have
\begin{equation}\label{x-far}
 |x_0| \lesssim R + 1/\eta_2.
\end{equation}
Next, we obtain a lower bound on $R$.  Recall from the normalization $N(t_0) = 1$ and
Corollary \ref{no-freq-dispersion} that
$$
 \| P_{> C(\eta_1)} u(t_0) \|_{\dot H^1} \lesssim \eta_1.
$$
On the other hand, from \eqref{rs-contradiction} we have
$$ \int_{B(x_0,R)} |\nabla u(t_0,x)|^2\ dx \gg \eta_1 $$
we thus see from the triangle inequality that
$$
\int_{B(x_0,R)} |\nabla P_{\leq C(\eta_1)} u(t_0,x)|^2\ dx \gg \eta_1.$$
But by H\"older, Bernstein \eqref{bernstein}, and \eqref{h1} we can bound the left-hand side by
$$ \lesssim R^3 \|\nabla P_{\leq C(\eta_1)} u(t_0)\|^2_{L^\infty_x}
\lesssim C(\eta_1) R^3 $$
and thus we have
$$ R \gtrsim c(\eta_1).$$
Combining this with \eqref{x-far} we see that the ball $B(x_0,
K(\eta_1,\eta_2) R)$ will contain $B(0, 1/\eta_2)$ (and hence any
ball of the form $B(0, C(\eta_1))$)
if the constant
$K(\eta_1,\eta_2)$ is large enough. In particular, from
Proposition \ref{concentration} we have
$$ \int_{B(x_0,K(\eta_1,\eta_2)R)} |u(t_0,x)|^6\ dx
\gtrsim c(\eta_1),$$
which inserted into
\eqref{rs-contradiction} contradicts the energy bound \eqref{h1}
if $K(\eta_1,\eta_2)$ is
chosen sufficiently large. This proves Proposition \ref{reverse}.
\endprf

\section{Interaction Morawetz: generalities}\label{slimflim-sec}

We shall shortly begin the proof of Proposition \ref{Morawetz}, which is a variant of the
interaction Morawetz inequality \eqref{morawetz-interaction}.  As noted above, this inequality
cannot be applied directly to our situation because the right-hand side of
\eqref{morawetz-interaction} can be very large due to low frequency  contributions to $u$.  It is
then natural (in light of \eqref{nothigh}) to try to adapt the interaction Morawetz inequality to
only deal with the high frequencies $u_{\geq 1}$, but this turns out to not quite be enough
either.  The trouble is that the inhomogeneous Schr\"odinger equation satisfied by $u_{\geq 1}$ is
not Lagrangian - and in particular it no longer enjoys the usual $L^2$ conservation.  Hence  when
we apply the argument from \cite{ckstt:cubic-scatter}, \cite{ckstt:french}  (which gave
\eqref{morawetz-interaction}) to this $u_{\geq 1}$ equation, we get new terms arising from the
fact that the right side of \eqref{local-mass-conserv} is no longer zero.  We can find no
appropriate bounds for these new terms. Our solution to this problem is to localize the previous
interaction Morawetz arguments  in space\footnote{See \eqref{localizedazero}, where the novelty
over the arguments from \cite{ckstt:cubic-scatter}, \cite{ckstt:french} is now the presence of the
spatial cut-off function.}, yielding a much more complicated version (Theorem \ref{slim-flim}) of
the inequality \eqref{morawetz-interaction}.

To summarize: the increased complexity in the right hand side of \eqref{full-morawetz} below is due to the fact that  we have localized the argument in frequency (because \eqref{nls} is critical) and space (because
of the error terms introduced by the frequency localization.)  Of course, all of these
extra terms will have to somehow be shown to be bounded, and to this end the second term on the {\emph{left}} side of \eqref{full-morawetz} is very important.  An analogue
of this term - where the $x$ integration is taken over
all of $\R^3$ - can also be included on
the left side of \eqref{morawetz-interaction} (see \cite{ckstt:cubic-scatter}, \cite{ckstt:french}),
but we had previously found no use for this term.  In what follows, the second term on the left side of \eqref{slim-flim} will be used to absorb - via the reverse Sobolev inequality of Proposition \ref{reverse} and an averaging argument -  some of the most troublesome
terms involving kinetic energy that appear on the right side of \eqref{full-morawetz}.\footnote{The
discussion here gives another way to frame our regularity argument which was sketched
in Section \ref{overview-sec}.  We only bother to show that a minimal energy blowup solution must
be localized in space in order that we can apply the reverse Sobolev inequality to such solutions.  The reverse Sobolev inequality is needed in the proof of the
frequency localized $L^4_{x,t}$ bound.}

The above argument  will be carried out in the next section; in this
section we prepare for our work involving the non-Lagrangian equation satisfied
by $u_{\geq 1}$ by discussing interaction Morawetz inequalities in more general situations than
the quintic NLS \eqref{nls}.  In particular, we shall consider general solutions $\phi$
to the equation \eqref{forced}, where $\mathcal{N}$ is an arbitrary nonlinearity.

\subsection{Virial-type Identity}

We introduce two related quantities which average the mass and
momentum densities (see Definition \ref{emtensor}) against a weight
function $a(x)$.

\begin{definition}
  Let $a(x)$ be a function\footnote{In other contexts it's useful to consider also
  time dependent weight functions $a(t,x)$.}
  defined on the spacetime slab $I_0 \times
  \R^3$. We define the associated virial potential
  \begin{equation}
    \label{Vsuba}
    V_a (t) = \int_{\R^3} a(x) |\p (t,x) |^2 dx
  \end{equation}
and the associated Morawetz action
\begin{equation}
  \label{Msuba}
  M_a (t) = \int_{\R^3} a_j 2 \Im (\overline{\phi}
  \p_j) dx.
\end{equation}
\end{definition}
A calculation using Lemma \ref{local-conserv} shows that
\begin{equation}
  \label{Vsubadot}
  \partial_t V_a = M_a + 2 \int_{\R^3} a \{ \lN, \phi \}_m dx,
\end{equation}
so $M_a = \partial_t V_a$ when $\lN = F' ( |\p|^2 ) \p$.
Using Lemma \ref{local-conserv}, a longer but similar calculation
establishes,

\begin{lemma}[Virial-type  Identity]
\label{genvirial}
Let $\p$ be a (Schwartz) solution of \eqref{forced}. Then
\begin{align}
  \label{Madot}
  \partial_t M_a = & \int_{\R^3} (- \Delta \Delta a) |\p|^2
+ 4 a_{jk} \Re (\overline{\p_j }  \p_k ) + 2 a_j \{ \lN , \p \}_p^j dx .
\end{align}
\end{lemma}

We now infer  a useful identity by choosing the weight function
$a(x)$ above to be,
\begin{equation}
  \label{localizedazero}
  a(x) = |x| \chi (|x|).
\end{equation}
where $\chi (r)$ denotes a
smooth nonnegative bump function defined on $r \geq 0$, supported on
$0 \leq r \leq 2$ and satisfying $\chi(r) = 1$ for $0 \leq r \leq 1.$
We calculate,
\begin{align*}
  a_j (x) &= \frac{x^j}{|x|} \tchi(|x|) ~{\mbox{where}}~ \tchi (r) =
  \chi(r) + r \chi' (r), \\
a_{jk} (x) &= \frac{1}{|x|} \left( \delta_{jk} - \frac{x^j}{|x|}
  \frac{x^k}{|x|} \right) \tchi (|x|) + \frac{x^j}{|x|} \frac{x^k}{|x|}
  \tchi' (|x|), \\
\Delta a(x) &= \frac{2}{|x|} \tchi(|x|) + \tchi' (|x|), \\
\Delta \Delta a(x) &= 2 \Delta \left( \frac{1}{|x|} \right)
  \tchi(|x|) + \psi (|x|),
\end{align*}
where $\psi (|x| )$ is smooth and supported in $1 \leq |x| \leq 2$.
Define now the notation $M^0 := M_a$ when $a(x)$ is chosen as in \eqref{localizedazero}.
(Later we will localize around a different fixed point $y \in \R^3$, in which
case we'll write the Morawetz action as $M^y$. Note that the letter $a$ is dropped completely
now from the notation for the Morawetz action.)
By the definition \eqref{Msuba},
\begin{equation}
  \label{lmp}
  M^0 (t) = 2 \Im \int_{\R^3} \frac{x^j}{|x|} \tchi ( |x| ) \phi_j (t,x)
  \overline{\p} (t,x) dx
\end{equation}
Note that\footnote{ One can also show that $|M^0 (t) | \lesssim {{\| \phi(t)
    \|}^2_{\dot{H}^{1/2}}}$; see Lemma 2.1 in \cite{ckstt:cubic-scatter}.}\label{half-footnote}
\begin{equation}
  \label{Mabound}
  |M^0 (t)| \leq 2 \| \phi (t) \|_{L^2_x }  {{\| \phi (t) \|}_{\dot{H}^1}}.
\end{equation}
Inserting the calculation below \eqref{localizedazero} above  into \eqref{Madot}, and using the
notation $\partial_{r(0)}$ to denote the radial part of the gradient with center at $0 \in \R^3$
(so $\partial_{r(0)} \equiv \frac{x}{|x|} \cdot \nabla$) we get,
\begin{align*}
  \partial_t M^0 & = - 2 \int_{\R^3} \Delta ( \frac{1}{|x|} )
  | \p (x) |^2 \tchi(|x|) dx \\
& + 4 \int_{\R^3} [ |\nabla \p (x)|^2 - |\partial_{r(0)} \p (x) |^2 ]
 \frac{1}{|x|} \tchi (|x|) dx \\
& + 2 \int_{\R^3} \frac{(x)}{|x|} \cdot \{ \lN, \phi \}_p \tchi
(|x|)
  dx \\
& - \int_{\R^3} |\p (x) |^2 \psi (|x|) dx + 4 \int_{\R^3} |\partial_{r(0)} \p |^2 \tchi'
 (|x|) dx
\end{align*}

As remarked above, we
translate the origin and choose, for fixed $y \in \R^3$
\begin{equation}
  \label{localizeday}
  a(x) = |x-y| \chi (|x-y|),
\end{equation}
instead of \eqref{localizedazero}, in which case the Morawetz action
\eqref{Msuba} is written $M^y$.   Then the preceding formula adjusts
to give the following  spatially localized virial-type identity,
where we write $\tchi$ for another bump function with the same properties
as $\chi$,
\begin{align}
\label{localizedvirialy}
  \partial_t M^y = & \\
\label{oney}
& - 2 \int_{\R^3} \Delta ( \frac{1}{|x-y|} ) {|\p (x )|^2} \tchi(|x-y|)
dx \\
\label{twoy}
& + 4 \int_{\R^3} |\nabb_y \p (x)|^2 \frac{1}{|x-y|} \tchi (|x-y|) dx \\
\label{threey}
& + 2 \int_{\R^3} \frac{(x-y)}{|x-y|} \cdot \{ \lN , \phi \}_p \tchi(|x-y|)
dx \\
\label{Ey}
& + O \left( \int_{\R^3} ( |\p (x) |^2 + |\partial_{r(y)} \p (x) |^2) |\psi
  (x-y) | dx \right).
\end{align}
We have used here the notation $\partial_{r(y)}$
to denote the radial portion of the gradient centered at $y$  and
$\nabb_y$ for the rest of the gradient.  That is,
\begin{equation*}
\partial_{r(y)} := \frac{x-y}{|x-y|} \cdot \nabla \quad \quad \text{and } \;
\nabb_y :=  \nabla - \frac{x-y}{|x-y|}( \frac{x-y}{|x-y|} \cdot  \nabla).
 \end{equation*}
We have also taken the liberty
to dismiss some of the structure in \eqref{Ey} using the fact that
$\tchi' $ has the same support properties as $\psi$.

\subsection{Interaction virial identity and general interaction Morawetz
 estimate for general equations}

When we choose $a(x) = |x| \chi(x)$ above , the virial potential reads
$V_a(t) = \int_{\R^3} |\phi(x,t)|^2 |x| \chi(x)dx$ and hence
$M^0(t) := \frac{d}{dt} V_a(t)$ might be thought of as measuring the extent
to which the mass of $\phi$ (near the origin at least) is moving away from the origin
at time $t$.  Similarly, for fixed $y \in \R^3$, $M^y(t)$ gives some measure of
the mass movement away from the point $y$.

Since we are ultimately interested in global decay and scattering properties of $\phi$, it's
reasonable to look for some measure of how the mass is moving away from (or {\em{interacting}}
with) {\em{itself}}. We might therefore sum over all $y \in \R^3$  the extent to which mass is
moving away from $y$ (that is, $M^y(t)$) multiplied by the amount of mass present at that point
$y$ (that is, $|\phi(y,t)|^2 dy$).  The result is the following quantity which we'll call the
{\it{spatially localized Morawetz interaction potential}}.
\begin{align}
  \label{slimp}
M^{interact} (t) & = \int_{\R^3_y} |\p (t,y)|^2 M^y (t) dy \\
& = 2\Im \int_{\R^3_y } \int_{\R^3_x} |\p (t,y)|^2
\tchi{(|x-y|)} \frac{(x-y)}{|x-y|} \cdot [\nabla \phi (t,x) ]
\overline{\p }(t,x) dx dy .
\end{align}
Note that, using \eqref{Mabound},
\begin{equation}
  \label{Mchibound}
  |M^{interact} (t) | \lesssim {{\| \phi (t) \|}_{L^2_x}^3} {{\| \phi (t)
   \|}_{\dot{H}^1_x}}.
\end{equation}

We calculate, using \eqref{localizedvirialy} and
\eqref{local-mass-conserv},
\begin{align}
\partial_t M^{interact} = & \int_{\R_y^3} |\p (y)|^2 \partial_t M^y dy \\
\label{intermed}
& + \int_{\R^3_y} [ 2 \partial_{y^k} \Im (\phi \overline{\p}_k ) (y) +
2 \{ \lN , \phi \}_m ] M^y (t) dy.
\end{align}
The $\partial_{y^k}$ appearing in \eqref{intermed} will now be
integrated by parts. Thus, using Lemma \ref{local-conserv} and
the fact that on $\R^3, \Delta \frac{1}{|x|} = - 4 \pi \delta$,
we have our spatially localized interaction
virial-type  identity
\begin{align}
  \label{slivi}
\partial_t M^{interact} = & \\
\label{slivione}
& 8 \pi \int_{\R^3}|\p (t,y) |^4  dy \\
\label{slivitwo}
& + 4\int_{\R^3_y} \int_{\R^3_x } |\p (t,y)|^2 \left[ \frac{1}{|x-y|}
\tchi (|x-y|) \right] |\nabb_y \p (t,x)|^2 dx dy \\
\label{slivithree}
& + 2 \int_{\R^3_y } \int_{\R^3_x } |\p (t,y)|^2 \left[\tchi (|x-y|)
  \frac{(x-y)}{|x-y|}\right]
\cdot \{
\lN , \phi \}_p (t,x) dx dy \\
\label{slivifour}
& - 4 \int_{\R^3_y} \int_{\R^3_x } \Im (\phi \overline{\p}_k ) (t,y) \partial_{y^k} \left[
  \frac{(x-y)^j}{|x-y|} \tchi (|x-y|) \right] \Im (\phi_j
  \overline{\phi} (t,x)) dx dy \\
\label{sliviE}
& + O \left( \int_{\R^3_y} \int_{\R^3_x } |\p (t,y)|^2 |\psi (|x-y|) | [ |\p(t,x) |^2 +
  |\partial_{r(y)} \p (t,x) |^2]  dx dy \right) \\
\label{slivimass}
& + 4 \int_{\R^3_y} \int_{\R^3_x } \{ \lN , \p \}_m (t,y) \left[ \tchi (|x-y|)
  \frac{(x-y)}{|x-y|} \right] \cdot \Im (\overline{\p} \nabla \p )
(t,x) dx dy.
\end{align}

The $\partial_{y^k}$ differentiation in \eqref{slivifour} produces two
terms. When $\partial_{y^k}$ falls on $\tchi(|x-y|)$, we encounter a
term controlled by \eqref{sliviE}. Indeed, we get a term bounded by
\begin{equation*}
  \int_{\R^3_y} \int_{\R^3_x } |\p (x) | |\partial_{r(y)} \p (x) | |\p(y) | |\partial_{r(x)}
  \p(y)| \tchi' (|x - y|) dx dy.
\end{equation*}
Upon grouping the terms in the integrand as $[|\p (x) ||\partial_{r(x)}
  \p(y)|][|\p (y) ||\partial_{r(y)}
  \p(x)|]$ and using $|ab| \lesssim |a|^2 + |b|^2$, we find the second
  part of the expression \eqref{sliviE}.
When the derivative falls on the unit
vector, we encounter a term of indeterminate sign but which is bounded
from below by \eqref{slivitwo}. We present the details (which also
  appear in Proposition 2.2 of \cite{ckstt:cubic-scatter}). The term we
are considering is (with $t$ dependence supressed)
\begin{align*}
  -4 & \int_{\R^3_y} \int_{\R^3_x} \Im (\phi \overline{\phi}_k )(y)
\left[ -\delta_{jk} + \frac{(x-y)^j (x-y)^k}{|x-y|^2} \right]
\Im (\phi_j \overline{\phi})(x) \left[ \frac{1}{|x-y|} \tchi(|x-y|)
\right] dx dy \\
& \geq -4 \int_{\R^3_x} \int_{\R^3_y} \left| \Im (\phi \nabb_x
  \overline{\phi} ) (y) \cdot \Im (\overline{\phi} \nabb_y \phi) (x)
\right|\left[ \frac{1}{|x-y|} \tchi(|x-y|)
\right] dx dy \\
& \geq -4  \int_{\R^3_x} \int_{\R^3_y} |\phi(y) | |\nabb_x \phi(y)|
|\phi(x)| |\nabb_y \phi (x) | \left[ \frac{1}{|x-y|} \tchi(|x-y|)
\right] dx dy \\
& \geq -2 \int_{\R^3_y} \int_{\R^3_x} (|\phi (y)|^2 |\nabb_y \phi(x)|^2
+ |\phi(x)|^2 |\nabb_x \phi(y)|^2) \left[ \frac{1}{|x-y|} \tchi(|x-y|)
\right] dx dy \\
& \geq \eqref{slivitwo}.
\end{align*}
Thus, apart from another term of the form \eqref{sliviE}, we have
shown that $-$\eqref{slivitwo} and \eqref{slivifour} together contribute a
nonnegative term.

Restricting the above calculations to the time interval $I_0 $,
we have the following useful estimate.

\begin{proposition}[Spatially Localized Interaction Morawetz
  Inequality]
\label{localized-morawetz}
Let $\phi$ be a (Schwartz) solution to the equation \eqref{forced} on a
spacetime slab $I_0 \times \R^3$ for some compact interval $I_0$. Then
we have
\begin{align*}
  8 \pi & \int_{I_0} \int_{\R^3_y} |\p (t,y)|^4 dy dt \\
& + 2 \int_{I_0} \int_{R^3_y} \int_{R^3_x}  |\p (t,y)|^2 \left[ \tchi(|x-y|)
  \frac{(x-y)}{|x-y|} \right] \cdot \{ \lN, \phi \}_{p} (t,x) dx dy dt
\\
\leq & 2 {{\| \p \|}^3_{L^\infty_t L^2_x ( I_0 \times \R^3 )}}
{{\| \p
      \|}_{L^\infty_t {\dot{H}}^1_x (I_0 \times \R^3 )}} + \\
& + 4 \int_{I_0} \int_{\R^3_y} \int_{\R^3_x}
 | \{ N, \p \}_m (t,y)|  |\tchi(|x-y|)| |\p
(t,x)| |\nabla \p (t,x) | dx dy dt \\
& + O \left(  \int_{I_0} \int_{\R^3_y} \int_{\R^3_x}
   |\p (t,y)|^2 |\psi (|x - y |)|
[ |\p (t,x) |^2 + |\partial_{r(y)} \p (t,x)|^2] dx dy dt \right) .
\end{align*}
\end{proposition}
The proof follows directly by integrating \eqref{slivi} over the
time interval $I_0$ using \eqref{Mchibound}.

\begin{remark}
\label{Rscaling}
If we replace \eqref{localizeday} by
\begin{equation}
  \label{localizedayR}
  a(x) = |x-y| \chi \left( \frac{|x-y|}{R} \right)
\end{equation}
then there are adjustments to the inequality obtained in
Proposition \ref{localized-morawetz}. Of course, $\tchi (\cdot) $
is replaced by $\tchi (\cdot/R)$. The annular cutoff $\psi$ in the
final term arises in the analysis above when derivatives fell on
$\chi$ or $\tchi$. Reviewing the derivation shows that this term
adjusts under \eqref{localizedayR} into
\begin{equation}
  \label{bigohtermR}
 O \left(  \int_{I_0} \int_{R^3_y} \int_{R^3_x}
   |\p (t,y)|^2 \left|\psi \left(\frac{|x - y |}{R}\right) \right|
\left[ \frac{1}{R^3} |\p (t,x) |^2 +
\frac{1}{R} |\partial_{r(y)} \p (t,x)|^2 \right] dx dy dt \right).
\end{equation}
If we send $R \to \infty$ and specialize to solutions of
\eqref{nls}, then we can apply \eqref{mass-cancel},
\eqref{momentum-cancel} to obtain the bound
$$  \int_{I_0} \int_{\R^3_y} |u (t,y)|^4 dy dt \\
 + \int_{I_0} \int_{R^3_y} \int_{R^3_x}  \frac{|u (t,y)|^2 |u(t,x)|^6}{|x-y|}
  dx dy dt \lesssim \| u \|_{L^\infty_t L^2_x(I_0 \times \R^3)}^3
\| u \|_{L^\infty_t \dot H^1_x(I_0 \times \R^3)},$$ which is
basically \eqref{morawetz-interaction} (see Footnote
\ref{half-footnote}). However for the purposes of proving Proposition
\ref{Morawetz}, it turns out to not be feasible to send $R \to
\infty$, as one of the error terms (specifically, the portion of
the mass bracket $\{ \lN, \phi \}_m$ which looks schematically
like $u_{hi}^5 u_{lo}$) becomes too difficult to estimate.
\end{remark}

\section{Interaction Morawetz: The setup and an averaging argument}
\label{interactionsetup}

 Having discussed interaction Morawetz
inequalities in general, we are now ready to begin the proof of
Proposition \ref{Morawetz}.

From the invariance of this Proposition under the scaling
\eqref{scaling} we may normalize $N_* = 1$.  Since we are assuming
$1 = N_* < c(\eta_3) N_{min}$, we have in particular that $1 <
c(\eta_3) N(t)$ for all $t \in I_0$. From Corollary
\ref{no-freq-dispersion} and Sobolev we have the low frequency
estimate
\begin{equation}\label{lowmedium-small}
 \| u_{< 1/\eta_3} \|_{L^\infty_t \dot H^1_x(I_0 \times \R^3)} + \| u_{< 1/\eta_3} \|_{L^\infty_t L^6_x(I_0 \times \R^3)}\lesssim \eta_3
\end{equation}
(for instance), if $c(\eta_3)$ was chosen sufficiently small. By \eqref{medfreq-largenergy} we thus see
that our choice of $N_*=1$ has
forced $N_{min} \geq c(\eta_0) \eta_3^{-1}$.
Define $P_{hi} := P_{\geq 1}$ and $P_{lo} := P_{<1}$, and then define $u_{hi} := P_{hi} u$ and
$u_{lo} := P_{lo} u$.  Observe from \eqref{lowmedium-small} that $u_{lo}$ has small energy,
\begin{equation}\label{low-small}
 \| u_{lo} \|_{L^\infty_t \dot H^1_x(I_0 \times \R^3)}
+ \| u_{lo} \|_{L^\infty_t L^6_x(I_0 \times \R^3)}\lesssim \eta_3.
\end{equation}
while from \eqref{lowmedium-small}, \eqref{nabla-highfreq} and \eqref{nothigh} we see that $u_{hi}$ has small
mass:
\begin{equation}\label{hi-small-l2}
 \| u_{hi} \|_{L^\infty_t L^2_x(I_0 \times \R^3)}
\lesssim \eta_3.
\end{equation}

We wish to prove \eqref{pm}, or in other words
\begin{equation}\label{Q-target}
\| u_{hi} \|_{L^4_t L^4_x(I_0 \times \R^3)}
\lesssim \eta_1^{1/4}.
\end{equation}
By a standard continuity argument\footnote{Strictly speaking, one needs to prove that
\eqref{Q-bootstrap} implies \eqref{Q-target} whenever $I_0$ is
replaced by any  sub-interval $I_1$ of $I_0$, in order to run the
continuity argument correctly, but it will be clear that the argument
below works not only for $I_0$ but also for all sub-intervals of
$I_0$.}, it will  suffice to show this under the bootstrap hypothesis
\begin{equation}\label{Q-bootstrap}
\| u_{hi} \|_{L^4_t L^4_x(I_0 \times \R^3)} \leq (C_0 \eta_1)^{1/4}
\end{equation}
where $C_0 \gg 1$ is a large constant (depending only on the
energy, and not on any of  the $\eta$'s). In practice we will
overcome this loss of $C_0$ with  a positive power of $\eta_3$ or
$\eta_1$.

We now use Proposition \ref{localized-morawetz} to obtain a Morawetz
estimate for  $\phi := u_{hi}$.

\begin{theorem}[Spatially and frequency localized interaction Morawetz
  inequality]\label{slim-flim}  Let the notation and assumptions be as
  above.  Then   for
any $R \geq 1$ we have
\begin{equation}\label{full-morawetz}
\int_{I_0} \int_{\R^3} |u_{hi}|^4\ dx dt
+ \int_{I_0} \int\int_{|x-y| \leq 2R} \frac{|u_{hi}(t,y)|^2
  |u_{hi}(t,x)|^6}{|x-y|}\ dx dy dt
\lesssim X_R,
\end{equation}
where $X_R$ denotes the quantity
\begin{align}
X_R &:= \eta_3^3 \label{morawetz-0} \\
&+ \int_{I_0} \int\int_{|x-y| \leq 2R} \frac{|u_{hi}(t,y)|^2
  |u_{lo}(t,x)|^5  |u_{hi}(t,x)|}{|x-y|}\ dx dy dt \label{morawetz-1} \\
&+ \sum_{j=0}^4 \int_{I_0} \int\int_{|x-y| \leq 2R}
|u_{hi}(t,y)| |P_{hi} \O(u_{hi}^j u_{lo}^{5-j})(t,y)|
|u_{hi}(t,x)| |\nabla u_{hi}(t,x)| \ dx dy dt \label{morawetz-2}\\
&+ \int_{I_0} \int\int_{|x-y| \leq 2R}  |u_{hi}(t,y)| |P_{lo} \O(u_{hi}^5) (t,y)|
|u_{hi}(t,x)| |\nabla u_{hi}(t,x)|
\ dx dy dt\label{morawetz-3}\\
&+ \eta_3^{1/10} \frac{1}{R} \int_{I_0} (\sup_{x \in \R^3}
\int_{B(x,2R)}  |u_{hi}(t,y)|^2\ dy) dt
\label{morawetz-4}\\
&+ \frac{1}{R} \int_{I_0} \int\int_{|x-y| \leq 2R} |u_{hi}(t,y)|^2
(|\nabla  u_{hi}(t,x)|^2 + |u_{hi}(t,x)|^6)\ dx dy dt.
\label{morawetz-5}
\end{align}
\end{theorem}

\begin{remark} This should be compared with \eqref{morawetz-interaction}.
The terms \eqref{morawetz-1}--\eqref{morawetz-5} may look fearsome, but
most of these  terms are
manageable, because of the spatial localization $|x-y| \leq 2R$, and
because there are  not too many derivatives on
the high-frequency term $u_{hi}$; the only truly difficult terms will
be the last two  \eqref{morawetz-4}, \eqref{morawetz-5}.  Observe from
\eqref{h1},  \eqref{l6} that we could control \eqref{morawetz-5} by
\eqref{morawetz-4}  if we dropped the $\eta_3^{1/10}$
factor from \eqref{morawetz-4}, however this type of factor is
indispensable in closing  our bootstrap argument, and
so we must treat
\eqref{morawetz-5} separately.  The idea is to use the reverse Sobolev
inequality,  Proposition \ref{reverse},
to control \eqref{morawetz-5} by the second term in
\eqref{full-morawetz}, plus an  error of the form
\eqref{morawetz-4}.  This can be done, but requires us to apply
Theorem \ref{slim-flim}  not just for $R=1$
(which would be the most natural choice of $R$) but rather for a range
of $R$ and then average over such $R$; see the discussion after the
proof of this Theorem.
\end{remark}

\begin{proof}  We apply $P_{hi}$ to \eqref{nls} to obtain
$$ (i\partial_t + \Delta) u_{hi} = P_{hi} (|u|^4 u),$$
and then apply Proposition \ref{localized-morawetz} with $\phi := u_{hi}$ and
$F := P_{hi}(|u|^4 u)$, to obtain
\begin{align*}
&c_1 \int_{I_0} \int_{\R^3} |u_{hi}(t,x)|^4\ dx dt \\
&+ c_2 \int_{I_0}
\int_{\R^3}\int_{\R^3} |u_{hi}(t,y)|^2 \tilde \chi(\frac{x-y}{R}) \frac{x-y}{|x-y|} \cdot \{P_{hi}(|u|^4 u),u_{hi}\}_\momentum(t,x)\ dx dy dt \\
&\lesssim \| u_{hi} \|_{L^\infty_t L^2_x(I_0 \times \R^3)}^3 \| u_{hi} \|_{L^\infty_t \dot H^1_x(I_0 \times \R^3)} \\
&+ \int_{I_0} \int\int_{|x-y| \leq 2R} |\{P_{hi}(|u|^4 u),u_{hi}\}_\mass(t,y)| |u_{hi}(t,x)| |\nabla u_{hi}(t,x)| \ dx dy dt\\
&+ \frac{1}{R} \int_{I_0} \int\int_{|x-y| \leq 2R} |u_{hi}(t,y)|^2 (\frac{1}{R^2} |u_{hi}(t,x)|^2 + |\nabla u_{hi}(t,x)|^2)\ dx dy dt.
\end{align*}

We now estimate the right-hand side by $X_R$.  Observe from \eqref{h1}, \eqref{hi-small-l2} that
$$ \| u_{hi} \|_{L^\infty_t L^2_x(I_0 \times \R^3)}^3 \| u_{hi} \|_{L^\infty_t \dot H^1_x(I_0 \times \R^3)}
\lesssim \eta_3^3 = \eqref{morawetz-0} \leq X_R,$$
while from \eqref{hi-small-l2} again we have
$$ \int_{\R^3} \frac{1}{R^2} |u_{hi}(t,x)|^2\ dx \lesssim \eta^2_3 /R^2 \leq \eta_3$$
and hence
$$
\frac{1}{R} \int_{I_0} \int\int_{|x-y| \leq 2R } |u_{hi}(t,y)|^2 \frac{1}{R^2} |u_{hi}(t,x)|^2 \ dx dy dt
\lesssim \eqref{morawetz-4} \leq X_R.$$
Similarly we have
$$
\frac{1}{R} \int_{I_0} \int\int_{|x-y| \leq 2 R } |u_{hi}(t,y)|^2 |\nabla u_{hi}(t,x)|^2 \ dx dy dt
\lesssim \eqref{morawetz-5} \leq X_R.$$

Now we deal with the mass bracket term.  We take advantage of the cancellation
\eqref{mass-cancel} to write
$$ \{P_{hi}(|u|^4 u),u_{hi}\}_\mass =
\{P_{hi}(|u|^4 u) - |u_{hi}|^4 u_{hi},u_{hi}\}_\mass.$$
We can write, using the notation \eqref{5-schematic},
$$ P_{hi}(|u|^4 u) - |u_{hi}|^4 u_{hi} =
P_{hi}(|u|^4 u - |u_{hi}|^4 u_{hi})
- P_{lo}(|u_{hi}|^4 u_{hi}) =
\sum_{j=0}^4 P_{hi}\O(u_{hi}^j u_{lo}^{5-j})
+ P_{lo} \O(u_{hi}^5 ).$$
Thus these terms can be bounded by $O(\eqref{morawetz-2} + \eqref{morawetz-3}) = O(X_R)$
(where we take absolute values everywhere).  To summarize
so far, we have shown that
\begin{align*}
c_1 \int_{I_0} \int_{\R^3} |u_{hi}(t,x)|^4\ dx dt &\\
+ c_2 \int_{I_0}
\int_{\R^3}\int_{\R^3} |u_{hi}(t,y)|^2 \tilde \chi(\frac{x-y}{R}) \frac{x-y}{|x-y|} \cdot \{P_{hi}(|u|^4 u),u_{hi}\}_\momentum(t,x)\ dx dy dt &\lesssim X_R.
\end{align*}
We now deal with the momentum bracket term, which is a little more involved as it requires a little
more integration by parts, and we will need to exploit the positivity of one of the components of this term.  In order to exploit the cancellation \eqref{momentum-cancel}, we break up the momentum bracket
into three pieces
\begin{align*}
 \{ P_{hi}(|u|^4 u), u_{hi} \}_\momentum &=
 \{ |u|^4 u, u_{hi} \}_\momentum -  \{ P_{lo}(|u|^4 u), u_{hi} \}_\momentum \\
&=  \{ |u|^4 u, u \}_\momentum -  \{ |u|^4 u, u_{lo} \}_\momentum  -  \{ P_{lo}(|u|^4 u), u_{hi} \}_\momentum\\
&= \{ |u|^4 u, u \}_\momentum - \{ |u_{lo}|^4 u_{lo}, u_{lo} \}_\momentum
- \{ |u|^4 u - |u_{lo}|^4 u_{lo}, u_{lo} \}_\momentum  -  \{ P_{lo}(|u|^4 u), u_{hi} \}_\momentum\\
&= -\frac{2}{3} \nabla (|u|^6 - |u_{lo}|^6)
- \{ |u|^4 u - |u_{lo}|^4 u_{lo}, u_{lo} \}_\momentum  -  \{ P_{lo}(|u|^4 u), u_{hi} \}_\momentum.
\end{align*}
We first deal with $\{ P_{lo}(|u|^4 u), u_{hi} \}_\momentum$; we estimate the contribution of this term
crudely in absolute values as
$$ O\left( \int_{I_0}
\int\int_{|x-y| \leq 2R} |u_{hi}(t,y)|^2 |\{u_{hi}, P_{lo}(|u|^4
u)\}_\momentum(t,x)|\ dx dy  dt \right).$$
We wish to bound this term by
$O(\eqref{morawetz-4}) = O(X_R)$.  
Recalling that any positive power of
$\eta_3$ overwhelms any loss in $R$ powers, this follows from

\begin{lemma}  We have $\int_{\R^3}|\{ u_{hi}, P_{lo}(|u|^4 u) \}_p|\ dx \lesssim \eta_3^{1/2}$.
\end{lemma}

\begin{proof}
By \eqref{h1}, H\"older and Bernstein \eqref{sobolev} we have
\begin{align*}
\int_{\R^3}|\{ u_{hi}, P_{lo}(|u|^4 u) \}_p|\ dx
&\lesssim \int_{\R^3} |\nabla u_{hi}| |P_{lo}(|u|^4 u)|\ dx +
 \int_{\R^3} |u_{hi}| |\nabla P_{lo}(|u|^4 u)|\ dx \\
&\lesssim \| \nabla u_{hi}\|_{L^2_x} \| P_{lo}(|u|^4 u) \|_{L^2_x} +
 \| u_{hi}\|_{L^2_x} \| \nabla P_{lo}(|u|^4 u) \|_{L^2_x} \\
&\lesssim \| P_{lo}(|u|^4 u) \|_{L^2_x}.
\end{align*}
Now decompose $u = u_{hi} + u_{lo}$, and use \eqref{5-schematic} to then decompose
$$ P_{lo}( |u|^4 u ) = \sum_{j=0}^5 P_{lo} \O(u_{hi}^j u_{lo}^{5-j}).$$
The terms $j=0,1,2,3,4$ can be estimated by using Bernstein \eqref{sobolev} and then H\"older to estimate
$$ \sum_{j=0}^4\| P_{lo} \O(u_{hi}^j u_{lo}^{5-j}) \|_{L^2_x}
\lesssim \sum_{j=0}^4 \| \O(u_{hi}^j u_{lo}^{5-j}) \|_{L^{6/5}_x}
\lesssim \sum_{j=0}^4 \| u_{hi} \|_{L^6_x}^j \| u_{lo} \|_{L^6_x}^{5-j} \lesssim \eta_3$$
thanks to \eqref{l6}, \eqref{lowmedium-small}.  For the $j=5$ term, we argue similarly, indeed we have
$$ \| P_{lo} \O(u_{hi}^5) \|_{L^2_x} \lesssim
\| \O(u_{hi}^5) \|_{L^1_x} \lesssim \| u_{hi} \|_{L^6_x}^{9/2} \| u_{hi} \|_{L^2_x}^{1/2} \lesssim \eta_3^{1/2}$$
by \eqref{h1}, \eqref{lowmedium-small}.
\end{proof}

Now we deal with the second term in the momentum bracket, namely
$\{ |u|^4 u - |u_{lo}|^4 u_{lo}, u_{lo} \}_\momentum$.  We
first move the derivative in \eqref{momentum-flux-def} into a more favorable
position, using the identity
$$ \{ f, g \}_\momentum = \nabla \O(fg) + \O(f \nabla g)$$
and the identity $|u|^4 u - |u_{lo}|^4 u_{lo} = \sum_{j=1}^5 \O(u_{hi}^j u_{lo}^{5-j})$
from \eqref{5-schematic}
to write
$$\{ |u|^4 u - |u_{lo}|^4 u_{lo}, u_{lo} \}_\momentum
= \sum_{j=1}^5
(\nabla \O( u_{hi}^j u_{lo}^{6-j} ) +
\O( u_{hi}^j u_{lo}^{5-j} \nabla u_{lo} )).$$
For the second term, we argue crudely again, estimating this contribution
in absolute value as
\begin{equation}
\label{topline}
 O(\int_{I_0}
\sum_{j=1}^5 \int\int_{|x-y| \leq 2R} |u_{hi}(t,y)|^2
|u_{hi}(t,x)|^j |u_{lo}(t,x)|^{5-j} |\nabla u_{lo}(t,x)|\ dx dy
dt).
\end{equation}
But by \eqref{l6}, \eqref{lowmedium-small} and the low
frequency localization of $u_{lo}$ we have
$$
 \int_{\R^3} |u_{hi}|^j |u_{lo}|^{5-j} |\nabla u_{lo}|\ dx
\lesssim \| u_{hi} \|_{L^6_x}^j \| u_{lo} \|_{L^6_x}^{5-j} \| \nabla u_{lo} \|_{L^6_x} \lesssim  \| u_{hi} \|_{L^6_x}^j \| u_{lo} \|_{L^6_x}^{6-j}  \lesssim \eta_3,
$$
and so this term can also be bounded by $O(\eqref{morawetz-4}) = O(X_R)$.
For the first term, we can integrate
by parts and then take absolute values to estimate the contribution of this term by
\begin{equation}\label{mixed}
 \sum_{j=1}^5
O(\int_{I_0}
\int\int_{|x-y| \leq 2R} \frac{|u_{hi}(t,y)|^2 |u_{hi}(t,x)|^j |u_{lo}(t,x)|^{6-j}}{|x-y|} \ dx dy
dt).
\end{equation}
Note that the terms where the integration by parts hits the cutoff is of the same type, but with the $\frac{1}{|x-y|}$
factor replaced by $O(\frac{1}{R})$ on the region where $|x-y| \sim R$; it is clear that this term
is dominated by \eqref{mixed} (or by a variant of \eqref{topline}, where we remove the $\nabla$ from $u_{lo}$).
We hold off on estimating this term for now, and turn to the
first term in the momentum bracket: $-\frac{2}{3} \nabla (|u|^6 - |u_{lo}|^6)$.  After an
integration by parts, this term can be written as
\begin{align*}
& c_3
\int_{I_0}
\int\int_{|x-y| \leq 2R} \frac{|u_{hi}(t,y)|^2 (|u(t,x)|^6 - |u_{lo}(t,x)|^6)}{|x-y|} \ dx dy  dt\\
&+ O(\frac{1}{R} \int_{I_0} \int\int_{|x-y| \leq 2R}
|u_{hi}(t,y)|^2 ||u(t,x)|^6 - |u_{lo}(t,x)|^6|| \ dx dy  dt),
\end{align*}
for some explicit constant $c_3 > 0$; note that the error incurred by removing the cutoff $\chi(\frac{x-y}{R})$ by
the cruder cutoff $|x-y| \leq 2R$ can be controlled by \eqref{morawetz-5} which is acceptable.
To control the error term, we use
\eqref{6-schematic} to split
\begin{equation}\label{hi-lo-split}
 |u(t,x)|^6 - |u_{lo}(t,x)|^6 = |u_{hi}(t,x)|^6 + \sum_{j=1}^5 \O(u_{hi}^j(t,x) u_{lo}^{6-j}(t,x)),
\end{equation}
The $|u_{hi}|^6$ term is of course bounded by $O(\eqref{morawetz-5}) = O(X_R)$.  The remaining terms
have an $L^1_x$ norm of $O(\eta_3)$ by \eqref{l6}, \eqref{lowmedium-small}, and H\"older, so
this error term is bounded by $O(\eqref{morawetz-4}) = O(X_R)$.  For the main term, we again
use \eqref{hi-lo-split}
and observe that the contribution of the error terms are bounded by \eqref{mixed}.  Collecting all these
estimates together, we have now shown that
\begin{align*}
c_1 \int_{I_0} \int_{\R^3} |u_{hi}(t,x)|^4\ dx dt &\\
+ c_3
\int_{I_0}
\int\int_{|x-y| \leq 2R} \frac{|u_{hi}(t,y)|^2 |u_{hi}(t,x)|^6}{|x-y|} \ dx dy  dt
&\lesssim X_R \\
+  \sum_{j=1}^5
O(\int_{I_0}
\int\int_{|x-y| \leq 2R} \frac{|u_{hi}(t,y)|^2 |u_{hi}(t,x)|^j |u_{lo}(t,x)|^{6-j}}{|x-y|}& \ dx dy  dt).
\end{align*}
We can eliminate the $j=2,3,4,5$ terms using the elementary inequality
$$ |u_{hi}(t,x)|^j |u_{lo}(t,x)|^{6-j} \leq \eps |u_{hi}(t,x)|^6 + C(\eps) |u_{hi}(t,x)| |u_{lo}(t,x)|^5$$
for some small absolute constant $\eps$; this allows us to control
the $j=2,3,4,5$ terms by the $j=1$ term, plus a small multiple of
the $j=6$ term which can then be absorbed by the second term on
the left-hand side (shrinking $c_3$ slightly).  The Theorem
follows.
\end{proof}

We would like to use Theorem \ref{slim-flim} to prove \eqref{Q-target}.  If it were not for the error terms
\eqref{morawetz-1}-\eqref{morawetz-5} then this estimate would follow immediately from  Theorem \ref{slim-flim},
since we can discard the second term in \eqref{full-morawetz} as being positive.  One would then hope
to estimate all of the error terms \eqref{morawetz-1}-\eqref{morawetz-5} by $O(\eta_1)$ using the bootstrap
hypothesis \eqref{Q-bootstrap} and the other estimates available (e.g. \eqref{h1}, \eqref{l6}, \eqref{lowmedium-small}, \eqref{hi-small-l2}).  It turns out that this strategy works (setting $R=1$) for the first four error terms
\eqref{morawetz-1}-\eqref{morawetz-4} but not for the fifth term \eqref{morawetz-5}.  To estimate this term
we need the reverse Sobolev inequality, which effectively replaces the $|\nabla u_{hi}|^2$ density here by
$|u_{hi}|^6$, at which point one can hope to control this term by the second positive term in \eqref{full-morawetz}.
But to do so it turns out that one cannot apply Theorem \ref{slim-flim} for a single value of $R$, but must instead
average over a range of $R$.

We turn to the details.  We let $J = J(\eta_1,\eta_2) \gg 1$ be a large\footnote{But it is not \emph{too}
large; any factor of $\eta_3$ shall easily be able to overcome any losses depending on $J$.} integer, and apply
Theorem \ref{slim-flim} to all dyadic radii $R = 1, 2, \ldots, 2^J$ separately.  We then average over $R$
to obtain
\begin{equation}\label{full-morawetz-average}
\int_{I_0} \int_{\R^3} |u_{hi}|^4\ dx dt + Y \lesssim X,
\end{equation}
where $Y$ is the positive quantity
$$ Y := \frac{1}{J} \sum_{1 \leq R \leq 2^J}
\int_{I_0} \int\int_{|x-y| \leq 2R} \frac{|u_{hi}(t,y)|^2 |u_{hi}(t,x)|^6}{|x-y|}\ dx dy dt $$
($R$ is always understood to sum over dyadic numbers), and $X$ is the quantity
\begin{align}
X &:= \eta_3^3 \label{smorawetz-0} \\
&+ \sup_{1 \leq R \leq 2^J} \int_{I_0} \int\int_{|x-y| \leq 2R} \frac{|u_{hi}(t,y)|^2 |u_{lo}(t,x)|^5 |u_{hi}(t,x)|}{|x-y|}\ dx dy dt \label{smorawetz-1} \\
&+ \sup_{1 \leq R \leq 2^J} \sum_{j=0}^4 \int_{I_0} \int\int_{|x-y| \leq 2R}
|u_{hi}(t,y)| |P_{hi}\O(u_{hi}^j u_{lo}^{5-j})(t,y)|
|u_{hi}(t,x)| |\nabla u_{hi}(t,x)| \ dx dy dt \label{smorawetz-2}\\
&+ \sup_{1 \leq R \leq 2^J} \int_{I_0} \int\int_{|x-y| \leq 2R}  |u_{hi}(t,y)| |P_{lo}\O(u_{hi}^5)(t,y)|
|u_{hi}(t,x)| |\nabla u_{hi}(t,x)|
\ dx dy dt\label{smorawetz-3}\\
&+ \eta_3^{1/10} \sup_{1 \leq R \leq 2^J} \frac{1}{R} \int_{I_0} (\sup_{x \in \R^3} \int_{B(x,2R)} |u_{hi}(t,y)|^2\ dy) dt
\label{smorawetz-4}\\
&+ \frac{1}{J} \sum_{1 \leq R \leq 2^J} \frac{1}{R} \int_{I_0} \int\int_{|x-y| \leq 2R} |u_{hi}(t,y)|^2 (|\nabla u_{hi}(t,x)|^2 + |u_{hi}(t,x)|^6)\ dx dy dt,
\label{smorawetz-5}
\end{align}
where we have estimated the average $\frac{1}{J} \sum_{1 \leq R \leq 2^J}$ by the supremum in those terms for
which the averaging is not important\footnote{Indeed, in those terms we will extract a gain of $\eta_3$ which
will easily absorb any losses relating to $R = O(2^J)$, since $J$ depends only on $\eta_1$ and $\eta_2$.}.

The terms \eqref{smorawetz-0}-\eqref{smorawetz-5} are roughly arranged in increasing order of difficulty to estimate.
But let us use the reverse
Sobolev inequality already obtained in Proposition \ref{reverse} to deal with the most difficult term \eqref{smorawetz-5},
replacing it by easier terms.

\begin{lemma} \label{11.4} We have
$$ \eqref{smorawetz-4} + \eqref{smorawetz-5} \lesssim
\eta_1^{1/100} (Y+W)$$
where
\begin{equation}\label{W-def}
 W := \sup_{1 \leq R \leq C(\eta_1,\eta_2) 2^J}
\frac{1}{R} \int_{I_0} (\sum_{x \in \frac{R}{100} \Z^3} (\int_{B(x,3R)} |u_{hi}(t,y)|^2\ dy)^{100})^{1/100} dt,
\end{equation}
and $\frac{R}{100} \Z^3$ is the integer lattice $\Z^3$ dilated by $R/100$.
\end{lemma}

\begin{proof}  We first handle \eqref{smorawetz-4}.  For every $x \in \R^3$ there exists $x' \in \frac{R}{100}\Z^3$
such that $B(x,2R)$ is contained in $B(x',3R)$.  Thus
$$ \frac{1}{R} \sup_{x \in \R^3} \int_{B(x,2R)} |u_{hi}(t,y)|^2\ dy \lesssim
\frac{1}{R}( \sum_{x' \in \frac{R}{100} \Z^3} (\int_{B(x' ,3R)}
|u_{hi}(t,y)|^2\ dy)^{100})^{1/100},$$ 
and
the claim follows by noting that $\eta_3^{1/10} \ll \eta_1^{1/100}
$.

Now we consider \eqref{smorawetz-5}.  We write this term as
$$ \frac{1}{J} \sum_{1 \leq R \leq 2^J} \frac{1}{R} \int_{I_0}
\int_{\R^3} |u_{hi}(t,y)|^2 e_{hi} (t,y,2R)\ dy dt,$$ where $e_{hi} (t,y,2R)$ is the local energy,
$$ e_{hi} (t,y,2R) := \int_{B(y,2R)} |\nabla u_{hi}(t,x)|^2 +
|u_{hi}(t,x)|^6\ dx.$$ We will denote the same quantity with
$u_{hi}$ replaced with $u$ by $e(t,y,2R)$.  We split this term
into the regions $e_{hi} (t,y,2R) \lesssim \eta_1$ and $e_{hi}
(t,y,2R) \gg \eta_1$.
\newline {\bf{Large Energy Regions:}}
\newline Consider first the large energy regions $e_{hi} (t,y,2R)
\gg \eta_1$.  By \eqref{low-small}, the same lower bound holds for
$e(t,y,2R)$.   For $(t,y)$ in such regions, we apply Proposition
\ref{reverse} to conclude that
$$ \int_{B(y, 2R)} |\nabla u (t,x)|^2 dx \leq \half e(t, y , 2R) + C(
\eta_1, \eta_2)  \int_{B(y, C(\eta_1 , \eta_2) R)} |u (t, x )|^6 dx $$
which implies
$$ \int_{B(y, 2R)} |\nabla u (t,x)|^2 dx \lesssim C(\eta_1 ,  \eta_2 ) \int_{B(y, C(\eta_1 ,
  \eta_2 ) R)} |u(t,x)|^6 dx $$  
and the same estimate is valid for $u_{hi}$ in light of
  \eqref{low-small}.
Thus we can bound the contribution to \eqref{smorawetz-5} of the large energy regions by
$$
 \lesssim C(\eta_1, \eta_2)
\frac{1}{J} \sum_{1 \leq R \leq 2^J} \frac{1}{R} \int_{I_0}
\int\int_{|x-y| \lesssim C(\eta_1, \eta_2) R} |u_{hi}(t,y)|^2 |u_{hi}(t,x)|^6\ dx dy dt.
$$
Shifting $R$ by $C(\eta_1, \eta_2)$, we can bound this by
\begin{equation}\label{high-energy-5}
 \lesssim (C(\eta_1,\eta_2))^2 \frac{1}{J} \sum_{1 \leq R \leq C(\eta_1,\eta_2) 2^J} A_R
\end{equation}
where
$$ A_R := \frac{1}{R} \int_{I_0} \int\int_{|x-y| \leq R}
|u_{hi}(t,y)|^2 |u_{hi}(t,x)|^6\ dx dy dt.$$
To bound this by $\eta_1^{1/100} (Y+W)$ (and not just by $O(Y)$) we exploit the
averaging\footnote{The key point here is that while the $A_R$ quantities have a factor of $1/R$,
the quantities in $Y$ have a larger factor of $1/|x-y|$.  After averaging, the latter factor begins
to dominate the former.} in $R$.
First observe that
$$ \sum_{1 \leq R \leq R'} A_R \lesssim
\int_{I_0} \int \int_{|x-y| \leq R'} \frac{|u_{hi}
  (t, y)|^2 |u_{hi} (t,x)|^6}{|x-y|}  dx dy dt$$
for all $R'$.  Averaging this over all $1 \leq R' \leq 2^J$, we see that
$$ \frac{1}{J} \sum_{1 \leq R' \leq 2^J} \sum_{1 \leq R \leq R'} A_R
\lesssim \frac{1}{J} \sum_{1 \leq R' \leq 2^J}
\int_{I_0} \int \int_{|x-y| \leq R'} \frac{|u_{hi}
  (t, y)|^2 |u_{hi} (t,x)|^6}{|x-y|}  dx dy dt \lesssim Y.$$
Now let $1 < J_0 < J$ be a parameter depending on $\eta_1, \eta_2$
to be chosen shortly. Observe that for each $1 \leq R \leq
2^{J-J_0}$ there are at least $J_0$ values of $R'$ which involve
that value of $R$ in the above sum.  Thus we have
$$ \frac{J_0}{J} \sum_{1 \leq R \leq 2^{J-J_0}} A_R \lesssim Y,$$
and thus the contribution to \eqref{high-energy-5} of the terms where $R \leq 2^{J-J_0}$ is
bounded by $O(\frac{(C(\eta_1,\eta_2))^2}{J_0} Y)$, which is acceptable if $J_0$ is chosen
sufficiently large depending on $\eta_1$, $\eta_2$.

It remains to control the expression
$$ \frac{(C(\eta_1 , \eta_2 ))^2}{J} \sum_{ 2^{J-J_0}
\leq R \leq C(\eta_1 ,
  \eta_2) 2^J} \frac{1}{R} \int_{I_0} \int \int_{|x-y| \leq R} |u_{hi}
  (t, y)|^2 |u_{hi} (t,x)|^6 dx dy dt .$$
This is bounded by
$$  \frac{C(\eta_1 , \eta_2, J_0 )}{J} \sup_{1 \leq R \leq C(\eta_1, \eta_2) 2^J}
  \frac{1}{R} \int_{I_0} \int \sup_{x \in \R^3} \left( \int_{B(x, R)} |u_{hi}
(t,y)|^2 dy \right) |u_{hi} (t,x)|^6 dx dt.$$
Using the energy bound \eqref{l6} this is in turn bounded by
\begin{align*}
& \lesssim \frac{C(\eta_1 , \eta_2, J_0 )}{J} \sup_{1 \leq R \leq C(\eta_1 , \eta_2) 2^J
}\frac{1}{R} \left( \int_{I_0} \sup_{x \in \R^3} \left( \int_{B(x, R)} |u_{hi} (t,y )|^2 dy
\right) dt \right)  \\
& \lesssim \frac{1}{J} C(\eta_1 , \eta_2, J_0 ) W\\
& \lesssim \eta_1^{1/100} W
\end{align*}
by choosing $J= J(\eta_1 , \eta_2 )$ sufficiently large.
\newline
{\bf{Small Energy Regions:}} \newline Now we deal with the contribution of the low energy regions:
$$ \frac{1}{J} \sum_{1 \leq R \leq 2^J} \frac{1}{R} \int_{I_0}
\int_{e_{hi} (t,y,2R) \lesssim \eta_1}
|u_{hi}(t,y)|^2 e_{hi} (t,y,2R)\ dy dt.$$
Observe that if $y$ is such that $e_{hi} (t,y,2R) \lesssim \eta_1$, then there exists $y' \in \frac{R}{100} \Z^3$
such that $y \in B(y',R) \subseteq B(y,2R) \subset B(y',3R)$, and thus
$e_{hi} (t,y,R) \lesssim \min(\eta_1, e_{hi} (t,y',3R))$.  Thus we can dominate the
above expression by
$$ \lesssim \frac{1}{J} \sum_{1 \leq R \leq 2^J} \frac{1}{R} \int_{I_0} \sum_{y' \in \frac{R}{100}\Z^3}
\min(\eta_1, e_{hi} (t,y',3R)) \int_{B(y',3R)} |u_{hi}(t,y)|^2 \ dy dt.$$
Now from \eqref{h1} we have, using $\min(\eta_1 ,
e_{hi}(\cdot))^{100/99} \leq \eta_1^{1/99} e_{hi} (\cdot)$, that
$$ \left(\sum_{y' \in \frac{R}{100}\Z^3} \min(\eta_1, e_{hi} (t,y',3R))^{100/99}\right)^{99/100}
\lesssim \eta_1^{1/100} (\sum_{y' \in \frac{R}{100}\Z^3} e_{hi} (t,y',3R))^{99/100} \lesssim \eta_1^{1/100},$$
so by H\"older's inequality we can bound the previous expression by $O(\eta_1^{1/100} W)$
as desired.
\end{proof}

In light of all the above estimates, we have thus shown that
\begin{equation}
\label{numberme}
 \int_{I_0} \int_{\R^3} |u_{hi}|^4\ dx dt \lesssim \eta_3^3 +
\eqref{smorawetz-1} + \eqref{smorawetz-2}
+ \eqref{smorawetz-3} + \eta_1^{1/100} W,
\end{equation}
since the $Y$ term on the left can be used to absorb the $\eta_1^{1/100} Y$ term which would otherwise
appear on the right.
It thus suffices to show that all of the quantities \eqref{smorawetz-1}-\eqref{smorawetz-3}
are $O(\eta_1)$, while the factor of $\eta_1^{1/100}$ allows the quantity $W$ to be estimated using
the weaker bound of $O(C_0 \eta_1)$.  (This is the main purpose of the reverse Sobolev inequality,
Proposition \ref{reverse}, in our argument.  The constant $C_0$ was defined in our bootstrap assumption \eqref{Q-bootstrap}.). 

As mentioned earlier, the terms \eqref{smorawetz-1}-\eqref{smorawetz-3} are roughly arranged in increasing order of difficulty, and $W$ is more difficult still.  To estimate any of these expressions, we of course need good spacetime
estimates on $u_{hi}$ and $u_{lo}$.  We do already have some estimates on these quantities \eqref{Q-bootstrap}, \eqref{h1}, \eqref{l6}, \eqref{low-small}, \eqref{hi-small-l2}, but it turns out that these are not directly sufficient to estimate
\eqref{smorawetz-1}-\eqref{smorawetz-3} and $W$.  Thus we shall first use the equation \eqref{nls} and Strichartz
estimates to bootstrap \eqref{Q-bootstrap} to yield further spacetime integrability; this will be the purpose
of the next section.

\section{Interaction Morawetz: Strichartz control}

In this section we establish the spacetime estimates we need in
order to bound \eqref{smorawetz-1}-\eqref{smorawetz-3} and $W$.
Ideally one would like to use \eqref{Q-bootstrap} to show that $u$
obeys the same estimates as a solution to the free Schr\"odinger
equation, however the quantity \eqref{Q-bootstrap} is
supercritical (it has roughly the scaling of $\dot H^{1/4}$,
instead of $\dot H^1$), and we must therefore accept some loss of
derivatives in the high frequencies\footnote{Recall that $u_{hi}$
and $u_{lo}$ were defined at the beginning of Section
\ref{interactionsetup}.}; a model example of a function $u$
obeying \eqref{Q-bootstrap} to keep in mind is a
pseudosoliton
solution where $u$
has magnitude $|u(t,x)| \sim N^{1/2}$ on the spacetime region $x =
O(N^{-1}), t = O(N)$ and has Fourier transform supported near the
frequency $N$ for some large $N \gg 1$.  We will however be able
to show that $u$ does behave like a solution to the free
Schr\"odinger equation modulo a high frequency forcing term which
is controlled in $L^2_t L^1_x$ but not in dual Strichartz spaces.

Recall that the constant $C_0$ which we use throughout these next three sections of the paper was defined in \eqref{Q-bootstrap}, our bootstrap hypothesis
on the $L^4$ norm of $u_{hi}$.  We begin by estimating the low frequency portion $u_{lo}$ of the
solution, which does behave like the free equation from the point
of view of Strichartz estimates:


\begin{proposition}[Low frequency estimate]\label{lowfreq-estimate}
For the functions $u_{hi}, u_{lo}$ defined at the start of Section \ref{interactionsetup}, we have
\begin{equation}\label{low-26}
\| \nabla u_{lo} \|_{L^2_t L^6_x(I_0 \times \R^3)} \lesssim C_0^{1/2} \eta_1^{1/2}
\end{equation}
and
\begin{equation}\label{low-4infty}
\| u_{lo} \|_{L^4_t L^\infty_x(I_0 \times \R^3)} \lesssim \eta_3^{1/2}.
\end{equation}
In fact, we have the slightly more general statement that \eqref{low-26}, \eqref{low-4infty} hold
for all $u_{\leq N}$ when $N \sim 1$.
\end{proposition}

Observe that the $L^4_t L^\infty_x$ estimate gains a power of $\eta_3$, which will be very useful for us in
overcoming certain losses of $R$ in the sequel.

\begin{proof}  We may assume that $N \geq 1$, since the case when $N < 1$
of course then follows.
Let $Z$ denote the quantity
\begin{equation}\label{X-def}
 Z := \| \nabla u_{\leq N} \|_{L^2_t L^6_x(I_0 \times \R^3)} + \| u_{\leq N} \|_{L^4_t L^\infty_x(I_0 \times \R^3)}.
\end{equation}
By Strichartz \eqref{strichartz}, \eqref{strichartz-components} we have
$$ Z \lesssim \| u_{\leq N} \|_{\dot S^1(I_0 \times \R^3)}
\lesssim \| u_{\leq N}(t_0) \|_{\dot H^1_x} + \| \nabla P_{\leq N}(|u|^4 u) \|_{L^2_t L^{6/5}_x(I_0 \times \R^3)}.$$
By \eqref{lowmedium-small} we have
$$ \| u_{\leq N}(t_0) \|_{\dot H^1_x} \lesssim \eta_3.$$
Now we consider the nonlinear term.  We split $u = u_{\leq N} + u_{>N}$ and use \eqref{5-schematic} to write
$$ \nabla P_{\leq N}(|u|^4 u) = \sum_{j=0}^5 \nabla P_{\leq N} \O( u_{> N}^j u_{\leq N}^{5-j} ).$$
Consider the $j=0$ term first.  We discard $P_{\leq N}$ and use the Leibnitz rule and
H\"older to estimate
\begin{align*} \| \nabla P_{\leq N} \O( u_{\leq N}^5 ) \|_{L^2_t L^{6/5}_x(I_0 \times \R^3)}
&\lesssim \| \O( u_{\leq N}^4 \nabla u_{\leq N} ) \|_{L^2_t L^{6/5}_x(I_0 \times \R^3)}\\
&\lesssim \| u_{\leq N} \|_{L^\infty_t L^6_x(I_0 \times \R^3)}^4 \| \nabla u_{\leq N} \|_{L^2_t L^6_x(I_0 \times \R^3)};
\end{align*}
by \eqref{lowmedium-small} and \eqref{X-def}
this term is thus bounded by $O(\eta_3^4 Z)$.

Now consider the $j=1$ term.  We argue similarly to estimate this term as
\begin{align*}
 \| \nabla P_{\leq N} \O( u_{\leq N}^4 u_{>N} ) \|_{L^2_t L^{6/5}_x(I_0 \times \R^3)}
&\lesssim \| \O( u_{\leq N}^4 \nabla u_{> N} ) \|_{L^2_t L^{6/5}_x(I_0 \times \R^3)}\\
&+ \| \O( u_{\leq N}^3 u_{> N} \nabla u_{\leq N} ) \|_{L^2_t L^{6/5}_x(I_0 \times \R^3)}\\
&\lesssim
\| u_{\leq N} \|_{L^4_t L^\infty_x(I_0 \times \R^3)}^2
\| u_{\leq N} \|_{L^\infty_t L^6_x(I_0 \times \R^3)}^2
\| \nabla u_{>N} \|_{L^\infty_t L^2_x(I_0 \times \R^3)}\\
&\quad +
\| u_{\leq N} \|_{L^\infty_t L^6_x(I_0 \times \R^3)}^3
\| u_{> N} \|_{L^\infty_t L^6_x(I_0 \times \R^3)}
 \| \nabla u_{\leq N} \|_{L^2_t L^6_x(I_0 \times \R^3)};
\end{align*}
using \eqref{l6}, \eqref{lowmedium-small}, \eqref{X-def} these terms
are thus bounded by $O(\eta_3^2 Z^2 + \eta_3^3 Z)$.

Now consider the $j=2,3,4,5$ terms.  This time we use Bernstein's inequality \eqref{bernstein}
(recalling that $N \sim 1$) and H\"older to estimate
\begin{align*}
 \| \nabla P_{\leq N} \O( u_{\leq N}^{5-j} u_{>N}^j ) \|_{L^2_t L^{6/5}_x(I_0 \times \R^3)}
&\lesssim \| \O( u_{\leq N}^{5-j} u_{> N}^j ) \|_{L^2_t L^1_x(I_0 \times \R^3)}\\
&\lesssim \| u_{> N} \|_{L^4_t L^4_x(I_0 \times \R^3)}^2 \| u_{>N} \|_{L^\infty_t L^6_x(I_0 \times \R^3)}^{j-2}
\| u_{\leq N} \|_{L^\infty_t L^6_x(I_0 \times \R^3)}^{5-j}.
\end{align*}
Applying \eqref{Q-bootstrap}, \eqref{l6} we can bound these terms by
$O(C_0^{1/2} \eta_1^{1/2})$ (we can do significantly better on the $j=2,3,4$ terms but
we will not exploit this).  Combining all these estimates we see that
$$ Z \lesssim \eta_3 + \eta_3^4 Z + \eta_3^2 Z^2 + \eta_3^3 Z + C_0^{1/2} \eta_1^{1/2};$$
by standard continuity arguments this then implies that $Z \lesssim C_0^{1/2} \eta_1^{1/2}$.
This proves \eqref{low-26}, but we did not achieve the $\eta_3$ gain in \eqref{low-4infty}.

To obtain \eqref{low-4infty}  we must refine the above analysis.  Let $u^0_{\leq N}$ be the
solution to the free Schr\"odinger equation with initial data $u^0_{\leq N}(t_0) = u_{\leq N}(t_0)$
for some $t_0 \in I_0$.  Then by \eqref{lowmedium-small} and Lemma \ref{4-lemma}
$$ \| u^0_{\leq N} \|_{L^4_t L^\infty_x(I_0 \times \R^3)} \lesssim \eta_3.$$
Thus it suffices to prove that
$$ \| u_{\leq N} - u^0_{\leq N} \|_{L^4_t L^\infty_x(I_0 \times \R^3)} \lesssim \eta_3^{1/2}.$$
We estimate the left-hand side as
$$ \sum_{N' \lesssim N} \| u_{N'} - u^0_{N'} \|_{L^4_t L^\infty_x(I_0 \times \R^3)}.$$
By Bernstein \eqref{bernstein} we may bound this by
$$ \sum_{N' \lesssim N} N' \| u_{N'} - u^0_{N'} \|_{L^4_t L^3_x(I_0 \times \R^3)}$$
which by interpolation can be bounded by
$$ \sum_{N' \lesssim N} N' \| u_{N'} - u^0_{N'} \|_{L^2_t L^6_x(I_0 \times \R^3)}^{1/2}
\| u_{N'} - u^0_{N'} \|_{L^\infty_t L^2_x(I_0 \times \R^3)}^{1/2}.$$
But from \eqref{lowmedium-small} we have
$$ \| u_{N'} - u^0_{N'} \|_{L^\infty_t L^2_x(I_0 \times \R^3)} \lesssim
(N')^{-1} \| u_{\leq N} - u^0_{\leq N} \|_{L^\infty_t \dot H^1(I_0 \times \R^3)} \lesssim
(N')^{-1} \eta_3$$
so we have
$$ \| u_{\leq N} - u^0_{\leq N} \|_{L^4_t L^\infty_x(I_0 \times \R^3)} \lesssim \eta_3^{1/2}
\sum_{N' \leq N} (N')^{1/2} \| u_{N'} - u^0_{N'} \|_{L^2_t L^6_x(I_0 \times \R^3)}^{1/2}.$$
Applying Strichartz \eqref{strichartz} we thus have
$$ \| u_{\leq N} - u^0_{\leq N} \|_{L^4_t L^\infty_x(I_0 \times \R^3)} \lesssim \eta_3^{1/2}
\sum_{N' \leq N} (N')^{1/2} \| P_{N'} (|u|^4 u) \|_{L^2_t L^{6/5}_x(I_0 \times \R^3)}^{1/2}.$$
But the preceding analysis already showed that
$$  \| P_{N'} (|u|^4 u) \|_{L^2_t L^{6/5}_x(I_0 \times \R^3)}
\lesssim \eta_3^4 Z + \eta_3^2 Z^2 + \eta_3^3 Z + C_0^{1/2} \eta_1^{1/2} \lesssim 1$$
(for instance), and the claim follows.
\end{proof}


Now we estimate the high frequencies.  We first need an estimate
on the high-frequency portion of the nonlinearity $|u|^4 u$.  It
turns out that this quantity cannot be easily estimated in a
single Strichartz norm, but must instead be decomposed into two
pieces estimated using separate space-time Lebesgue norms (cf. the
appearance of $M$ in Lemma \ref{disjointed-strichartz}).

\begin{proposition}\label{high-decomp}  We can decompose
$$ P_{hi}(|u|^4 u) = F + G$$
where $F$, $G$ are Schwartz functions with Fourier support in the region $|\xi| \gtrsim 1$ and
$$ \| \nabla F \|_{L^2_t L^{6/5}_x(I_0 \times \R^3)} \lesssim \eta_3^{1/2}$$
and
$$ \| G \|_{L^2_t L^1_x(I_0 \times \R^3)} \lesssim C_0^{1/2} \eta_1^{1/2}$$
\end{proposition}


Of the two pieces, $F$ is by far the better term; indeed, if $G$
were not present, then the Strichartz estimate \eqref{strichartz}
would be able to obtain $L^{10}_{t,x}$ bounds on $u_{hi}$. The
reader may in fact assume as a first approximation that $F$ is
negligible, and that the nonlinearity $P_{hi}(|u|^4 u)$ is
primarily in $L^2_t L^1_x$, which is not a dual  $\dot S^1$
Lebesgue norm. Note also that the $\eta_1^\half$ bound ultimately determines the $\eta_1$ on the
right side of \eqref{pm} at the end of Section \ref{flim-conclude}.

\begin{proof}
We split $u = u_{lo} + u_{hi}$, and then use \eqref{5-schematic}
to split
$$ P_{hi}( |u|^4 u ) = \sum_{j=0}^5 P_{hi} \O( u_{hi}^j u_{lo}^{5-j} ).$$
Consider the $j=0$ term first.  We have
\begin{align*}
 \| \nabla P_{hi} \O( u_{lo}^5 ) \|_{L^2_t L^{6/5}_x(I_0 \times \R^3)}
&\lesssim \| \O( u_{lo}^4 \nabla u_{lo} ) \|_{L^2_t L^{6/5}_x(I_0 \times \R^3)}\\
&\lesssim \| u_{lo} \|_{L^4_t L^\infty_x(I_0 \times \R^3)}^2 \| u_{lo} \|_{L^\infty_t L^6_x(I_0 \times \R^3)}^2
\|\nabla u_{lo} \|_{L^\infty_t L^2_x(I_0 \times \R^3)}
\end{align*}
which is $O(\eta_3^{1/2})$ (for instance) by Proposition \ref{lowfreq-estimate} and
\eqref{l6}, \eqref{h1}.  Thus this term may be placed as part of $F$.

Now consider the $j=1$ term.  We have
\begin{align*}
\| \nabla P_{hi} \O( u_{lo}^4 u_{hi} ) \|_{L^2_t L^{6/5}_x(I_0 \times \R^3)}
&\lesssim
\| \O( u_{lo}^4 \nabla u_{hi} ) \|_{L^2_t L^{6/5}_x(I_0 \times \R^3)}\\
&+ \| \O( u_{lo}^3 u_{hi} \nabla u_{lo}) \|_{L^2_t L^{6/5}_x(I_0 \times \R^3)} \\
&\lesssim
\| u_{lo} \|_{L^4_t L^\infty_x(I_0 \times \R^3)}^2 \| u_{lo} \|_{L^\infty_t L^6_x(I_0 \times \R^3)}^2
\|\nabla u_{hi} \|_{L^\infty_t L^2_x(I_0 \times \R^3)}\\
&+ \| u_{lo} \|_{L^\infty_t L^6_x(I_0 \times \R^3)}^3
\| u_{hi} \|_{L^\infty_t L^6_x(I_0 \times \R^3)} \| \nabla u_{lo} \|_{L^2_t L^6_x(I_0 \times \R^3)}.
\end{align*}
Applying Proposition \ref{lowfreq-estimate}, \eqref{l6},
\eqref{h1}, and \eqref{lowmedium-small} this expression is
$O(\eta_3^{3})$
and so this term may also be placed as
part of $F$.

Now consider the $j=2,3,4,5$ terms. We estimate
$$ \| P_{hi} \O( u_{hi}^j u_{lo}^{5-j} ) \|_{L^2_t L^1_x(I_0 \times \R^3)}
\lesssim \| u_{hi} \|_{L^4_t L^4_x(I_0 \times \R^3)}^2
\| u_{hi} \|_{L^\infty_t L^6_x(I_0 \times \R^3)}^{j-2}
\| u_{lo} \|_{L^\infty_t L^6_x(I_0 \times \R^3)}^{5-j}$$
which is $O(C_0^{1/2} \eta_1^{1/2})$ by \eqref{Q-bootstrap} and
\eqref{l6}.  Thus we may place this term as part of $G$.
\end{proof}

\begin{corollary}\label{high-26-cor} For every $N \geq 1$, we have
\begin{equation}\label{high-26}
\| u_N \|_{L^2_t L^6_x(I_0 \times \R^3)} \lesssim C_0^{1/2} N^{1/2} \eta_1^{1/2}.
\end{equation}
\end{corollary}

\begin{proof}
From Strichartz \eqref{strichartz} and Proposition \ref{high-decomp}
$$ \| u_N \|_{L^2_t L^6_x(I_0 \times \R^3)} \lesssim \| u_N(t_0) \|_{L^2_x(\R^3)} + \| P_N F \|_{L^2_t L^{6/5}_x(I_0 \times \R^3)}
+ \| P_N G \|_{L^2_t L^{6/5}_x(I_0 \times \R^3)}$$
The first term is certainly acceptable by \eqref{hi-small-l2}.  The second term is $O(\eta_3^{1/2} N^{-1})$ by
Proposition \ref{high-decomp}, and the third term is $O(C_0^{1/2} N^{1/2} \eta_1^{1/2})$ by Bernstein
\eqref{bernstein} and Proposition \ref{high-decomp}.  The claim follows.
\end{proof}

\section{Interaction Morawetz: Error estimates}


We now show that the comparatively easy terms \eqref{smorawetz-1},
\eqref{smorawetz-2}, \eqref{smorawetz-3} are indeed controlled by
$O(\eta_1)$, which is in turn bounded by the right hand side of
\eqref{Q-target}. The term $W$ is however significantly harder and
will be deferred to the next section.

\divider{The estimation of \eqref{smorawetz-1}.}

We have to show that
$$ \int_{I_0} \int\int_{|x-y| \leq 2R} \frac{|u_{hi}(t,y)|^2 |u_{lo}(t,x)|^5 |u_{hi}(t,x)|}{|x-y|}\ dx dy dt \lesssim \eta_1$$
for all $1 \leq R \leq 2^J$. This term will be fairly easy because
of the localization $|x-y| \leq 2R$. Observe that the kernel $
\frac{1}{|x|}$ has an $L^1_x$ norm of $O(R^2) \leq O(2^{2J})$ on
the ball $B(0,2R)$.  Thus by Young's inequality and Cauchy-Schwarz
we see that
$$ \int\int_{|x-y| \leq 2R} \frac{F(x) G(y)}{|x-y|}\ dx dy \lesssim R^2\| F\|_{L^2_x} \| G \|_{L^2_x}$$
 for any functions $F$,
$G$. In particular the expression to be estimated is bounded by
$$ \lesssim
2^{2J} \int_{I_0} \| u_{hi}(t) \|_{L^4_x}^2
\| |u_{lo}|^5 |u_{hi}|(t) \|_{L^2_x}\ dt.$$
We use H\"older and \eqref{hi-small-l2} to estimate
$$ \| |u_{lo}|^5 |u_{hi}|(t) \|_{L^2_x} \lesssim \| u_{lo}(t) \|_{L^\infty_x}^{5}
\| u_{hi}(t) \|_{L^2_x} \lesssim \eta_3 \| u_{lo}(t)
\|_{L^\infty_x}^{5}.$$ We dispose of three of the five factors of
$\|u_{lo}(t) \|_{L^\infty_x}^5$ by observing from \eqref{l6} and
Bernstein's inequality \eqref{bernstein} that $ \| u_{lo}(t)
\|_{L^\infty_x}^3 \lesssim \eta_3^3 \lesssim 1$. Combining all these estimates and
then using H\"older in time, we thus can bound the expression to
be estimated by
$$ 2^{2J} \eta_3
\| u_{hi} \|_{L^4_t L^4_x(I_0 \times \R^3)}^2 \| u_{lo} \|_{L^4_t
L^\infty_x(I_0 \times \R^3)}^2,$$ which is bounded by the right
hand side of \eqref{Q-target}  using  \eqref{Q-bootstrap} and
Proposition \ref{lowfreq-estimate} (note that $\eta_3$ will absorb
$2^{2J}$ since $J$ depends only on $\eta_1$ and $\eta_2$). This concludes the
treatment of \eqref{smorawetz-1}.

\divider{The estimation of \eqref{smorawetz-2}.}

We now handle \eqref{smorawetz-2}.  We have to show that
$$
\int_{I_0} \int\int_{|x-y| \leq 2R} |u_{hi}(t,y)|
|P_{hi}\O(u_{hi}^j u_{lo}^{5-j})(t,y)| |u_{hi}(t,x)| |\nabla
u_{hi}(t,x)| \ dx dy dt \lesssim \eta_1$$ for $j=0,1,2,3,4$ and $1
\leq R \leq 2^J$.

We begin by considering the cases $j=1,2,3,4$. We observe from
H\"older and \eqref{h1} that
$$
 \int_{B(y,2R)} |u_{hi}(t,x)| |\nabla u_{hi}(t,x)|\ dx \lesssim
R^{3/4} \| u_{hi}(t) \|_{L^4_x} \| \nabla u_{hi}(t) \|_{L^2_x}
\lesssim R^{3/4} \| u_{hi}(t) \|_{L^4_x},
$$
and hence it suffices to show that
$$
 R^{\frac{3}{4}} \int_{I_0} \| u_{hi}(t) \|_{L^4_x} \int_{\R^3}
|u_{hi}(t,y)| |P_{hi}\O(u_{hi}^j u_{lo}^{5-j})(t,y)|\ dy dt
\lesssim  \eta_1.
$$
 We use
H\"older (and the hypothesis $j=1,2,3,4$) to estimate this as
$$ \lesssim R^{\frac{3}{4}} \| u_{hi} \|_{L^4_t L^4_x(I_0 \times \R^3)}^3
\| u_{lo} \|_{L^4_t L^\infty_x(I_0 \times \R^3)} \| u_{hi} \|_{L^\infty_t L^6_x(I_0 \times \R^3)}^{j-1}
\| u_{lo} \|_{L^\infty_t L^6_x(I_0 \times \R^3)}^{4-j}.$$
The $L^\infty_t L^6_x(I_0 \times \R^3)$ factors are bounded by \eqref{l6}, leaving us to prove
$$  R^{\frac{3}{4}} \| u_{hi} \|_{L^4_t L^4_x(I_0 \times \R^3)}^3
\| u_{lo} \|_{L^4_t L^\infty_x(I_0 \times \R^3)} \lesssim
\eta_1.$$ 
 But this follows from \eqref{Q-bootstrap} and Proposition
\ref{lowfreq-estimate}, again using $\eta_3$ to wallop a positive
power of $R$.  

Finally, we consider the $j=0$ case of \eqref{smorawetz-2}, where we have to prove
$$
\int_{I_0} \int\int_{|x-y| \leq 2R} |u_{hi}(t,y)|
|P_{hi}\O(u_{lo}^{5})(t,y)| |u_{hi}(t,x)| |\nabla u_{hi}(t,x)| \
dx dy dt \lesssim \eta_1.$$ Here we use Cauchy-Schwarz and
\eqref{hi-small-l2} to crudely bound
$$
\int_{\R^3} |u_{hi}(t,x)| |\nabla u_{hi}(t,x)| \ dx
\lesssim \eta_3$$
and then use H\"older to reduce to showing that
$$ \eta_3 \| u_{hi} \|_{L^\infty_t L^2_x(I_0 \times \R^3)}
\| P_{hi} \O(u_{lo}^{5}) \|_{L^1_t L^2_x(I_0 \times \R^3)}
\lesssim \eta_1.$$
The factor $\| u_{hi} \|_{L^\infty_t L^2_x(I_0 \times \R^3)}$ is $O(\eta_3)$
by \eqref{hi-small-l2}.  For the second factor, we take advantage of
the high frequency localization (using \eqref{nabla-highfreq}) to write
\begin{align*}
 \| P_{hi} \O(u_{lo}^{5}) \|_{L^1_t L^2_x(I_0 \times \R^3)}
&\lesssim \| \nabla \O(u_{lo}^{5}) \|_{L^1_t L^2_x(I_0 \times \R^3)}\\
&\lesssim \| \O(u_{lo}^4 \nabla u_{lo})\|_{L^1_t L^2_x(I_0 \times \R^3)}\\
&\lesssim \| \nabla u_{lo} \|_{L^\infty_t L^2_x(I_0 \times \R^3)} \| u_{lo} \|_{L^4_t L^\infty_x(I_0 \times \R^3)}^4.
\end{align*}
The claim then follows from \eqref{h1} and Proposition \ref{lowfreq-estimate}.

\divider{The estimation of \eqref{smorawetz-3}.}

We now handle \eqref{smorawetz-3}.  We have to show that
\begin{equation}\label{s3-targ}
\int_{I_0} \int\int_{|x-y| \leq 2R}  |u_{hi}(t,y)| |P_{lo}\O(u_{hi}^5)(t,y)|
|u_{hi}(t,x)| |\nabla u_{hi}(t,x)|
\ dx dy dt \lesssim \eta_1.
\end{equation}
We begin by using H\"older to write
$$ \int_{B(y,2R)} |u_{hi}(t,x)| |\nabla u_{hi}(t,x)|\ dx
\lesssim R^{1/2} \| u_{hi}(t) \|_{L^3_x} \| \nabla u_{hi}(t) \|_{L^2_x}$$
and apply \eqref{h1} to estimate the left-hand side of \eqref{s3-targ} by
$$ \lesssim
R^{\frac{1}{2}} \int_{I_0} \| u_{hi}(t) \|_{L^3_x} \int
|u_{hi}(t,y)| |P_{lo}\O(u_{hi}^5)(t,y)|\ dy.$$ 
By H\"older and Bernstein \eqref{bernstein} we
have
\begin{align*}
\int |u_{hi}(t,y)| |P_{lo}\O(u_{hi}^5)(t,y)|\ dy
&\lesssim \| u_{hi}(t) \|_{L^3_x} \| P_{lo}\O(u_{hi}(t)^5) \|_{L^{3/2}_x}\\
&\lesssim  \| u_{hi}(t) \|_{L^{3}_x}  \| \O(u_{hi}(t)^5) \|_{L^1_x}\\
&\lesssim   \| u_{hi}(t) \|_{L^{3}_x} \| u_{hi}(t) \|_{L^4_x}^2 \| u_{hi}(t) \|_{L^6_x}^3.
\end{align*}
The terms $\| u_{hi}(t) \|_{L^6_x}^3$ are bounded by \eqref{l6}, so by a H\"older in time and
\eqref{Q-bootstrap} we have
\begin{equation}\label{m2-pause}
  \eqref{smorawetz-3} \lesssim
R^{\frac{1}{2}} \| u_{hi} \|_{L^4_t L^3_x(I_0 \times \R^3)}^2 \|
u_{hi} \|_{L^4_t L^4_x(I_0 \times \R^3)}^2 \lesssim R^{1/2}
C_0^{1/2} \eta_1^{1/2} \| u_{hi} \|_{L^4_t L^3_x(I_0 \times
\R^3)}^2.
\end{equation}
From the triangle inequality and H\"older we have
$$ \| u_{hi} \|_{L^4_t L^3_x(I_0 \times \R^3)} \lesssim \sum_{N \geq 1} \| u_N \|_{L^4_t L^3_x(I_0 \times \R^3)}
\lesssim \sum_{N \geq 1} \| u_N \|_{L^\infty_t L^2_x(I_0 \times \R^3)}^{1/2} \| u_N \|_{L^2_t L^6_x(I_0 \times \R^3)}^{1/2}.$$
From \eqref{hi-small-l2}, \eqref{nothigh} we have $\| u_N \|_{L^\infty_t L^2_x(I_0 \times \R^3)} \lesssim \min(\eta_3,N^{-1})$.
Applying \eqref{high-26} we thus have
$$ \| u_{hi} \|_{L^4_t L^3_x(I_0 \times \R^3)} \lesssim \sum_{N \geq 1} \min(\eta_3,N^{-1})^{1/2}
C_0^{1/4} N^{1/4} \eta_1^{1/4}  \lesssim \eta_3^{1/2} C_0^{1/4} \eta_1^{1/4}.$$
Inserting this estimate into \eqref{m2-pause} we see that \eqref{smorawetz-3}
is acceptable (again, the power of $\eta_3$ counteracts the loss in
$C_0$ and the presence of the power $R^{1/2}$).

Note that the first four factors on the right side of \eqref{numberme} have all in fact been shown to
be controlled with a positive power of $\eta_3$.

\section{Interaction Morawetz: A double Duhamel trick}\label{flim-conclude}

To conclude the proof of Proposition \ref{Morawetz}, we have to show that
$$ W \lesssim C_0 \eta_1,$$
or in other words that
$$
\frac{1}{R} \int_{I_0} (\sum_{x \in \frac{R}{100} \Z^3} (\int_{B(x,3R)} |u_{hi}(t,y)|^2\ dy)^{100})^{1/100} dt
\lesssim C_0 \eta_1
$$
for all $1 \leq R \leq C(\eta_1, \eta_2) 2^J$.  By duality we have
$$(\sum_{x \in \frac{R}{100} \Z^3} (\int_{B(x,3R)} |u_{hi}(t,y)|^2\ dy)^{100})^{1/100}
= \sum_{x \in \frac{R}{100} \Z^3} c(t,x) \int_{B(x,3R)} |u_{hi}(t,y)|^2\ dy$$
where $c(t,x) > 0$ are a collection of numbers which are almost summable in the sense that
\begin{equation}\label{c-bound}
 \sum_{x \in \frac{R}{100} \Z^3} c(t,x)^{100/99} = 1
\end{equation}
for all $t$.  Thus it suffices to show that
$$
\frac{1}{R} \int_{I_0} \sum_{x \in \frac{R}{100} \Z^3} c(t,x) \int_{B(x,3R)} |u_{hi}(t,y)|^2\ dy dt
\lesssim C_0 \eta_1.
$$
Let $\psi$ be a bump function adapted to $B(0,5)$ which equals 1 on $B(0,3)$.  Since
$$ \int_{B(x,3R)} |u_{hi}(t,y)|^2\ dy \lesssim \int |u_{hi}(t,y)|^2 \psi(\frac{y-x}{R})\ dy,$$
it suffices to show that
\begin{equation}\label{royal-pain}
\frac{1}{R} \int_{I_0} \sum_{x \in \frac{R}{100} \Z^3} c(t,x) \int_{\R^3} |u_{hi}(t,y)|^2
\psi(\frac{y-x}{R})\ dy dt
\lesssim C_0 \eta_1.
\end{equation}

In the proof of Lemma \ref{11.4}, we obtained spacetime control on
$u_{hi}$ by using the (forward-in-time) Duhamel formula \eqref{duhamel} followed
by Strichartz estimates.  This seems to be insufficient to prove
\eqref{royal-pain} (the best argument available seems to lose a
logarithm of the derivative); instead, we rely on both the
forward-in-time and backward-in-time Duhamel formulae \eqref{duhamel}, and argue
using the fundamental solution \eqref{fundamental-soln}. 
This will only lose a factor of $C_0$ (see \eqref{Q-bootstrap}), which is acceptable because of the gain of
$\eta_1^{1/100}$ which was 
obtained earlier in Lemma \ref{11.4}.

Let us write $I_0 = [t_-, t_+]$ for some times $-\infty < t_- < t_+ < \infty$.  We use Proposition
\ref{high-decomp} to decompose $P_{hi}(|u|^4 u) = F + G$.  We define the functions $u^\pm_{hi}$ to solve
the Cauchy problem
$$ (i\partial_t + \Delta) u^\pm_{hi} = F; \quad u^\pm_{hi}(t_\pm) = u_{hi}(t_\pm);$$
thus this is the same equation as $u_{hi}$ satisfies but without the term $G$.  From \eqref{duhamel}
observe the forward-in-time
Duhamel formula
$$ u_{hi}(t) = u^-_{hi}(t) - i \int_{t_- < s_- < t} e^{i(t-s_-)\Delta} G(s_-)\ ds_-$$
and the backward-in-time Duhamel formula
$$ u_{hi}(t) = u^+_{hi}(t) + i \int_{t < s_+ < t_+} e^{i(t-s_+)\Delta} G(s_+)\ ds_+.$$

Let us see how we would prove \eqref{royal-pain} if $u_{hi}$ were replaced by $u^\pm_{hi}$.
From Strichartz \eqref{strichartz}, \eqref{hi-small-l2}, and Proposition \ref{high-decomp} (discarding a derivative) we see that
$$ \| u^\pm_{hi} \|_{L^2_t L^6_x(I_0 \times \R^3)} \lesssim \eta_3^{1/2}.$$
But from H\"older we have
\begin{align*}
\frac{1}{R} \sum_{x \in \frac{R}{100} \Z^3} c(t,x)
&\int_{\R^3} |u^\pm_{hi}(t,y)|^2 \psi(\frac{y-x}{R}) \ dy\\
&\lesssim
\frac{1}{R} \sum_{x \in \frac{R}{100} \Z^3} c(t,x) R^2
(\int_{\R^3} |u^\pm_{hi}(t,y)|^6 \psi(\frac{y-x}{R}) \ dy)^{1/3} \\
&\lesssim  R
(\sum_{x \in \frac{R}{100} \Z^3} c(t,x)^{3/2})^{2/3}
(\sum_{x \in \frac{R}{100} \Z^3} \int_{\R^3} |u^\pm_{hi}(t,y)|^6 \psi(\frac{y-x}{R}) \ dy)^{1/3} \\
&\lesssim R (\sum_{x \in \frac{R}{100} \Z^3} c(t,x)^{100/99})^{99/100} \| u^\pm_{hi} \|_{L^6_x}^2\\
&= R \| u^\pm_{hi} \|_{L^6_x}^2
\end{align*}
and hence
$$
\frac{1}{R} \int_{I_0} \sum_{x \in \frac{R}{100} \Z^3} c(t,x) \int_{\R^3} |u_{hi}(t,y)|^2
\psi(\frac{y-x}{R})\ dy dt
\lesssim R \eta_3$$
which is acceptable if $\eta_3$ is sufficiently small (recall that $R \leq C(\eta_1, \eta_2) 2^J$ and $J$ depends only on
$\eta_1, \eta_2$).  Thus we see that \eqref{royal-pain} would be easy to prove if $u_{hi}$ were replaced by
$u^\pm_{hi}$.

It is now natural to use one of the two Duhamel formulae listed
above, and attempt to prove \eqref{royal-pain} for the integral
term.  This however turns out to be rather difficult. It will be
significantly easier if we use \emph{both} formulae
simultaneously.  More precisely, we re-arrange the above Duhamel
formulae as
$$ -i\int_{t_- < s_- < t} e^{i(t-s_-)\Delta}G(s_-)\ ds_- = u_{hi}(t) - u^-_{hi}(t)$$
and
$$ i\int_{t < s_+ < t_+} e^{i(t-s_+)\Delta}G(s_+)\ ds_+ = u_{hi}(t) - u^+_{hi}(t).$$
Then we multiply the first identity by the conjugate of the second to obtain
\begin{align*}
 - \int\int_{t_- < s_- < t < s_+ < t_+}&
(e^{i(t-s_-)\Delta}G(s_-))
(\overline{e^{i(t-s_+)\Delta}G(s_+))}\ ds_+
\ ds_- \\
&= |u_{hi}(t)|^2 - u^-_{hi}(t) \overline{u_{hi}(t)} - u_{hi}(t) \overline{u^+_{hi}(t)} + u^-_{hi}(t) \overline{u^+_{hi}(t)}.
\end{align*}
From the elementary pointwise estimates
\begin{align*}
|u^-_{hi}(t) \overline{u_{hi}(t)}| &\leq \frac{1}{4} |u_{hi}(t)|^2 + O(|u^-_{hi}(t)|^2) \\
|u_{hi}(t) \overline{u^+_{hi}(t)}| &\leq \frac{1}{4} |u_{hi}(t)|^2 + O(|u^+_{hi}(t)|^2) \\
|u^-_{hi}(t) \overline{u^+_{hi}(t)}| &\leq  O(|u^-_{hi}(t)|^2) + O(|u^+_{hi}(t)|^2)
\end{align*}
we thus have the pointwise inequality
\begin{equation}\label{dodge-bullet}
 |u_{hi}(t)|^2 \lesssim |\int\int_{t_- < s_- < t < s_+ < t_+}
e^{i(t-s_-)\Delta} G(s_-)
\overline{e^{i(t-s_+)\Delta} G(s_+)}\ ds_- ds_+| + |u^-_{hi}(t)|^2 + |u^+_{hi}(t)|^2.
\end{equation}
This should be compared with what one would obtain with a single Duhamel formula \eqref{duhamel}, namely something like
$$ |u_{hi}(t)|^2 \lesssim |\int\int_{t_- < s, s' < t}
e^{i(t-s)\Delta} G(s)
\overline{e^{i(t-s')\Delta} G(s')}\ ds ds'| + |u^-_{hi}(t)|^2.$$
This turns out to be an inferior formulation; the basic problem will be that the integral $\int_{t_- < s,s' < t}
\frac{ds ds'}{|s-s'|}$ is logarithmically divergent, whereas the integral
$\int_{t_- < s_- < t < s_+ < t_+} \frac{ds_- ds_+}{|s_- - s_+|}$ is not.

We now return to \eqref{royal-pain}, and insert \eqref{dodge-bullet}.  The latter two terms were already shown to
be acceptable.  So we are left with proving that
\begin{equation}\label{bullet-time}
\begin{split}
\frac{1}{R} |\int\int\int_{t_- < s_- < t < s_+ < t_+}
\sum_{x \in \frac{R}{100} \Z^3} c(t,x)
\int_{\R^3} &
e^{i(t-s_-)\Delta}G(s_-)(y)
\overline{e^{i(t-s_+)\Delta}G(s_+)}(y)\\
& \psi(\frac{y-x}{R})\ dy ds_- ds_+ dt|
\lesssim C_0\eta_1.
\end{split}
\end{equation}
To compute the $y$ integral, we use the following stationary phase estimate.

\begin{lemma}\label{thanks-eli}  For any $t_- < s_- < t < s_+ < t_+$, and any Schwarz functions $f_-(x), f_+(x)$, we have
\begin{align*}
 |\sum_{x \in \frac{R}{100} \Z^3} c(t,x) \int_{\R^3} e^{i(t-s_-)\Delta} &f_-(y)
\overline{e^{i(t-s_+)\Delta} f_+(y)}
 \psi(\frac{y-x}{R})\ dy| \\
&\lesssim |s_+ - s_-|^{-3/2} \min(\theta^{-3/2+3/100}, 1)
 \| f_- \|_{L^1_x} \| f_+\|_{L^1_x}
\end{align*}
where $\theta := \frac{|t-s_+| |t-s_-|}{R^2 |s_+ - s_-|}$.
\end{lemma}

\begin{proof}  Fix $t$.
We use the explicit formula for the fundamental solution \eqref{fundamental-soln}
to estimate the left-hand side as
\begin{align*} \lesssim &\frac{1}{|t-s_+|^{3/2} |t-s_-|^{3/2}}
\sum_{x \in \frac{R}{100} \Z^3} c(t,x)
\int_{\R^3}\int_{\R^3} \\
&(\int_{\R^3} e^{i(|y-x_-|^2/(t-s_-) - |y-x_+|^2/(t-s_+))}
\psi(\frac{y-x}{R}) dy) f_+(x_+) \overline{f_-(x_-)}\ dx_+ dx_-,
\end{align*}
so it suffices to show that
\begin{align*}
|\sum_{x \in \frac{R}{100} \Z^3} c(t,x) \int_{\R^3} &e^{i(|y-x_-|^2/(t-s_-) - |y-x_+|^2/(t-s_+))}
\psi(\frac{y-x}{R}) dy|\\
&\lesssim \frac{|t-s_+|^{3/2} |t-s_-|^{3/2}}{|s_+-s_-|^{3/2}} \min(\theta^{-3/2+3/100},1) \\
&=R^3 \min( \theta^{3/100}, \theta^{3/2} )
\end{align*}
for all $x_-, x_+$.  Making the change of variables $y = x + Rz$, this becomes
$$
|\sum_{x \in \frac{R}{100} \Z^3} c(t,x) I(x)|
\lesssim \min( \theta^{3/100}, \theta^{3/2} ).$$
where the integrals $I(x)$ are defined by
$$
I(x) := \int_{\R^3} e^{i \Phi_x(z)}
\psi(z) dz
$$
and $\Phi_x = \Phi_{x,R,x_-,x_+,s_-,s_+,t}$ is the phase
$$ \Phi_x(z) := |x-x_- + Rz|^2/(t-s_-) - |x-x_+ + Rz|^2/(t-s_+).$$
By the normalization of $c(t,x)$, it thus suffices to prove the bounds
$$ (\sum_{x \in \frac{R}{100} \Z^3} I(x)^{100})^{1/100}
\lesssim \min( \theta^{3/100}, \theta^{3/2} ).$$
We now divide into two cases, depending on the size of $\theta$.  First
suppose that $\theta \gg 1$.
Observe that the gradient of the phase $\Phi_x$ in $z$
is
$$ \nabla_z \Phi_x(z) = 2R (x - x_- + Rz)/(t-s_-) - 2R (x-x_+ + Rz)/(t-s_+)
= 2R (x+Rz - x_*) \frac{s_- - s_+}{(t-s_-)(t-s_+)}$$
where $x_* = x_*(s_-,s_+,t,x_-,x_+,R)$ is a quantity not depending
on $x$ or $z$.  In particular, in the region where
$$ |x-x_*| \gg \frac{|t-s_-| |t-s_+|}{R |s_+-s_-|} = R\theta \gg R$$
we can obtain an extremely good bound from the principle of nonstationary phase\footnote{
In other words, one can use repeated integration by parts; see \cite{stein:large}.  An alternate
approach in this lemma is to use a Gaussian cutoff $e^{-\pi |x|^2}$ instead of $\psi$, and
then compute all the integrals explicitly by contour integration (or equivalently by using
the ``Gaussian beam'' solutions of the Schr\"odinger equation).}, namely that
\begin{equation}\label{I-nonstationary}
 |I(x)| \lesssim (\frac{|x-x_*|}{R\theta})^{-100}.
\end{equation}
In the remaining cases where
$$ |x-x_*| \lesssim R\theta$$
(note that there are $O(\theta^3)$ such cases), we use the crude bound
$$ |I(x)| \lesssim 1,$$
and obtain the final estimate
$$ (\sum_{x \in \frac{R}{100} \Z^3} I(x)^{100})^{1/100}
\lesssim \theta^{3/100}$$
as desired.  Now consider the case $\theta \lesssim 1$.  In the region
$$ |x-x_*| \gg R$$
one can get very good bounds from nonstationary phase again, namely \eqref{I-nonstationary}.  There are only
$O(1)$ remaining values of $x$.  For each of these values we observe that the double derivative $\nabla^2_z \Phi_x$
is nondegenerate, indeed it is equal to $2R^2 \frac{s_- - s_+}{(t-s_-)(t-s_+)} = \frac{2}{\theta}$
times the identity matrix.  Thus
by the principle of stationary phase (see \cite{stein:large}) we have
$$ |I(x)| \lesssim \theta^{3/2}$$
in these cases.  Upon summing we obtain
$$ (\sum_{x \in \frac{R}{100} \Z^3} I(x)^{100})^{1/100} \lesssim \theta^{3/2}
$$
as claimed.
\end{proof}

Using Lemma \ref{thanks-eli}, we can estimate the left-hand side of \eqref{bullet-time} as
\begin{equation}\label{final-stretch}
\begin{split}
\lesssim
&\int\int\int_{t_- < s_- < t < s_+ < t_+}
\frac{1}{R} |s_+ - s_-|^{-3/2} \\
&\min((\frac{|t-s_+| |t-s_-|}{R^2 |s_+ - s_-|})^{-3/2+3/100}, 1)
\| G(s_-)\|_{L^1_x} \| G(s_+)\|_{L^1_x}\ ds_- ds_+ dt.
\end{split}
\end{equation}
Now we use the crucial time ordering $s_- < t < s_+$. An elementary computation (treating the cases
$s_+ - s_- < R^2$ and $s_+ - s_- \geq R^2$ separately) shows that
$$ \frac{1}{R} \int_{s_- < t < s_+}
|s_+ - s_-|^{-3/2} \min((\frac{|t-s_+| |t-s_-|}{R^2 |s_+ - s_-|})^{-3/2+3/100}, 1)\ dt
\lesssim \min( R|s_+ - s_-|^{-3/2}, R^{-1} |s_+ - s_-|^{-1/2} ).$$
The kernel $\min(R |s|^{-3/2}, R^{-1} |s|^{-1/2} )$ has an $L^1_s$ norm of\footnote{It is crucial
to note here that the powers of $R$ have cancelled out.  This seems to be a consequence of
dimensional analysis, although the presence of the frequency localization $|\xi| \gtrsim 1$
makes this analysis heuristic rather than rigorous.} $O(1)$.  Thus by
Young's inequality in time (and Cauchy-Schwarz in time), we can bound \eqref{final-stretch} by
$$ \| G\|_{L^2_t L^1_x(I_0 \times \R^3)}^2$$
which is $O(C_0 \eta_1)$ by Proposition \ref{high-decomp}, as desired.  This (finally!) concludes
the proof of Proposition \ref{Morawetz}.

\endprf



\section{Preventing energy evacuation}\label{evacuate-sec}



We now prove Proposition \ref{energy-travel}.  By the scaling
\eqref{scaling} we may take $N_{min} = 1$.

\subsection{The setup and contradiction argument}

Since the $N(t)$ can only take values
in a discrete set (the integer powers of 2), there thus exists a time
$t_{min} \in I_0$ such that
$$ N(t_{min}) = N_{min} = 1.$$
At this time $t = t_{min}$,
we see from \eqref{medfreq-largenergy} and \eqref{nabla-block} that
we have a substantial amount of mass\footnote{Note that we are
\emph{not} using the assumption that $u$ is Schwartz (and thus has finite $L^2$ norm) to obtain these estimates;
the bounds here are independent of the global $L^2$ norm of $u$, which may be very large even for fixed energy $E$
because the very low frequencies can simultaneously have small energy and large mass, and is also not preserved by the scale-invariance which we have exploited to normalize $N_{min} = 1$.} (and energy) at medium
frequencies:
\begin{equation}\label{0-medium-large}
\| P_{c(\eta_0) \leq \cdot \leq C(\eta_0)} u(t_{min}) \|_{L^2(\R^3)} \geq c(\eta_0).
\end{equation}
This should be contrasted with \eqref{nothigh}, which shows that there is not much mass at the frequencies much
higher than $C(\eta_0)$.

Our task is to prove \eqref{delta}.  Suppose for contradiction that this estimate failed; then
there exists a time $t_{evac} \in I_0$ for which $N(t_{evac}) \gg C(\eta_5)$.  If $C(\eta_5)$ is sufficiently
large, we then see from Corollary \ref{no-freq-dispersion} that energy has been almost
entirely evacuated from low and medium frequencies at time $t_{evac}$:
\begin{equation}\label{T-medium-verysmall}
\| P_{< 1/\eta_5} u(t_{evac}) \|_{\dot H^1(\R^3)} \leq \eta_5.
\end{equation}
Under the intuition that the $L^2$ mass density does not
rapidly adjust during the NLS evolution, \eqref{T-medium-verysmall} will be
contradicted if the low frequency $L^2$ mass \eqref{0-medium-large} sticks around $N_{min}$ until $t_{evac}$. We will validate this slow $L^2$ mass motion intuition by proving a frequency localized $L^2$
mass almost conservation law which leads to the contradiction, provided
that the $\eta_j$ are chosen small enough.


Having surveyed the argument, we carry out the details.  Fix $t_{evac}$; by time reversal symmetry
we may assume that $t_{evac} > t_{min}$.  (From \eqref{0-medium-large} it is clear that $t_{evac}$
cannot equal $t_{min}$).


From \eqref{0-medium-large} we have
$$
\| u_{> c(\eta_0)}(t_{min}) \|_{L^2(\R^3)} \geq \eta_1
$$
(for instance).  In particular, if we set
$$ P_{hi} := P_{\geq \eta_4^{100}}, P_{lo} := P_{< \eta_4^{100}}, u_{hi} :=
P_{hi} u, u_{lo} := P_{lo} u$$
then we have
\begin{equation}\label{0-hi}
\| u_{hi}(t_{min}) \|_{L^2(\R^3)} \geq \eta_1.
\end{equation}
Suppose we could show that
\begin{equation}\label{T-hi}
\| u_{hi}(t_{evac}) \|_{L^2(\R^3)} \geq \frac{1}{2} \eta_1.
\end{equation}
From \eqref{h1} and \eqref{nothigh} we thus have
$$
\| P_{\leq C(\eta_1)} u_{hi}(t_{evac}) \|_{L^2(\R^3)} \geq \frac{1}{4} \eta_1.
$$
By \eqref{nabla-highfreq} this implies that
$$
\| P_{\leq C(\eta_1)} u(t_{evac}) \|_{\dot H^1(\R^3)} \gtrsim c(\eta_1, \eta_4).
$$
But then Corollary \ref{no-freq-dispersion} implies that $N(t_{evac}) \lesssim C(\eta_1, \eta_4)$, which
contradicts \eqref{T-medium-verysmall} if $\eta_5$ is chosen sufficiently small.


It remains to prove \eqref{T-hi}.  We use the continuity method.
Suppose we have a time $t_{min} \leq t_* \leq t_{evac}$ for which
\begin{equation}\label{l2-prebootstrap}
\inf_{t_{min} \leq t \leq t_*} \| u_{hi}(t) \|_{L^2(\R^3)} \geq \frac{1}{2}
\eta_1.
\end{equation}
We will show that \eqref{l2-prebootstrap} can be bootstrapped to
\begin{equation}\label{l2-postbootstrap}
\inf_{t_{min} \leq t \leq t_*} \| u_{hi}(t) \|_{L^2(\R^3)} \geq \frac{3}{4}
\eta_1.
\end{equation}
This implies that the set of times $t_*$ for which
\eqref{l2-prebootstrap} holds is both open and closed in $[t_{min},
t_{evac}]$, which will imply \eqref{T-hi} as desired.  Note that the introduction of $\eta_4^{100} \ll N_{min} $ allows for the $L^2$ mass near $N_{min}$ at time $t_{min}$ to move toward low frequencies but the estimate \eqref{l2-postbootstrap} shows a portion of the mass stays above $\eta_4^{100}$ for $t \in [ t_{min}, t_{evac}].$ Note also that upon proving \eqref{l2-postbootstrap}, we will have established Lemma \ref{extralemma}.

It remains to derive \eqref{l2-postbootstrap} from
\eqref{l2-prebootstrap}.
The idea is to treat the $L^2$ norm of $u_{hi}(t)$, i.e.
$$ L(t) := \int_{\R^3} |u_{hi}(t,x)|^2\ dx,$$
as an almost conserved quantity.  From \eqref{0-hi} we have
$$ L(t_{min}) \geq \eta_1^2,$$
so by the Fundamental Theorem of Calculus it will suffice to show that
$$ \int_{t_{min}}^{t_*} |\partial_t L(t)|\ dt \leq \frac{1}{100} \eta_1^2.$$
From \eqref{local-mass-conserv} and \eqref{mass-cancel} we have
\begin{align*}
\partial_t L(t) &= 2 \int \{ P_{hi}( |u|^4 u ), u_{hi} \}_\mass\ dx \\
&= 2
 \int \{ P_{hi}( |u|^4 u ) - |u_{hi}|^4 u_{hi}, u_{hi} \}_\mass\ dx.
\end{align*}
Thus it will suffice to show that
\begin{equation}\label{l2-increment}
\int_{t_{min}}^{t_*} 
 |\int \{ P_{hi}( |u|^4 u ) - |u_{hi}|^4 u_{hi}, u_{hi} \}_\mass\ dx|
 dt\leq \frac{1}{100} \eta_1^2.
\end{equation}
The proof of \eqref{l2-increment} is accomplished using a quintilinear
analysis of various interactions using three inputs: $\dot S^1$ Strichartz
estimates on low frequencies, $L^4_{xt}$ estimates via frequency localized
interaction Morawetz and Bernstein estimates on medium and higher frequencies,
and $\dot{S}^1$ Strichartz estimates on very high frequencies.



\subsection{Spacetime estimates for high, medium, and low frequencies}


To prove \eqref{l2-increment} we need a number of spacetime bounds on $u$, which we now discuss
in this subsection.
Observe by the discussion following \eqref{T-hi}
that the hypothesis \eqref{l2-prebootstrap} implies in particular that
$$
N(t) \leq C(\eta_1,\eta_4) \hbox{ for all } t_{min} \leq t \leq t_*.
$$
This in turn can be combined with Corollary \ref{bounded-snake-ok} to obtain the useful Strichartz bounds
\begin{equation}\label{vhf-bounds}
\| u \|_{\dot S^1([t_{min}, t_*] \times \R^3)} \leq C(\eta_1, \eta_3, \eta_4).
\end{equation}
Because of the dependence of the right-hand side on $\eta_4$, this bound is only useful for us when
there is a power of $\eta_5$ present.  In all other circumstances,
we resort instead to Proposition \ref{Morawetz}, which gives the $L^4_{t,x}$ bounds\footnote{It is important
here to note that while the $L^2$ control on $u_{hi}$ only extends from $t_{min}$ up to $t_*$, the $L^4_{t,x}$
control on $u_{\geq N}$ extends all the way up to $t_{evac}$.  This allows us to access the evacuation hypothesis
\eqref{T-medium-verysmall} to provide useful new control, especially at low frequencies, in the time interval
$[t_{min},t_*]$.  This additional control will be crucial for us to obtain the desired almost conservation law
on the mass of $u_{hi}$, thus closing the bootstrap and allowing $t_*$ to extend all the way up to $t_{evac}$,
at which point we can declare a contradiction.}
\begin{equation}\label{boundedl4}
\int_{t_{min}}^{t_{evac}} \int |P_{\geq N} u(t,x)|^4\ dx dt \lesssim
\eta_1 N^{-3}
\end{equation}
whenever $N < c(\eta_3)$.  
Roughly speaking,
the bound \eqref{boundedl4} is better than \eqref{vhf-bounds} for medium frequencies, but
\eqref{vhf-bounds} is superior for very high frequencies, and Lemma \ref{l2l1-bound} below will be superior for
very low frequencies.

It turns out that we also need some bounds on the low frequencies such
as $u_{\leq\eta_4}$; the estimate \eqref{boundedl4} is inappropriate for this purpose
because $N^{-3}$ diverges as $N \to 0$.  One can modify Proposition \ref{lowfreq-estimate} to obtain
some reasonable control, but it turns out that the constants obtained by that estimate
are not strong enough to counteract the losses arising from \eqref{vhf-bounds}, and we need
a stronger version of Proposition \ref{lowfreq-estimate} which takes advantage of
the evacuation hypothesis \eqref{T-medium-verysmall}, which asserts among other things
that $u_{\leq \eta_4}$ has extremely small
energy at time $t_{evac}$.  Because of this hypothesis we expect $u_{\leq \eta_4}$ to
have extremely small energy at all other times in
$[t_{min},t_{evac}]$ (i.e. we expect bounds gaining an $\eta_5$
instead of just an $\eta_3$).
Of course there is a little bit of energy leaking from the high
frequencies to the low, but fortunately the $L^4_{t,x}$ bound on the
high frequencies will limit\footnote{It may seem surprising that the
$L^4_{t,x}$ bound, which is supercritical, can lead to control on
critical quantities such as the energy.  The point is that once one
localizes in frequency, the distinctions between subcritical,
critical, and supercritical quantities become less relevant, as one
can already see from Bernstein's inequality \eqref{bernstein}.  In
this section the entire analysis is localized around the frequency
$N_{min} = 1$, and so that supercritical norms (such as $L^4_{t,x}$
or $L^\infty_t L^2_x$) can begin to play a useful role.} how far the
high frequency energy can penetrate to the very low modes.
This intuition is made rigorous as follows:

\begin{lemma}\label{l2l1-bound}  With the above notation and assumptions (in particular, assuming
\eqref{T-medium-verysmall} and \eqref{boundedl4})
we have
\begin{equation}\label{lo-piece-eq}
\| P_{\leq N} u \|_{\dot S^1([t_{min},t_{evac}] \times \R^3)} \lesssim \eta_5 + \eta_4^{-3/2} N^{3/2}
\end{equation}
for all $N \leq \eta_4$.  
\end{lemma}

One should think of the $C\eta_5$ term on the right-hand side of
\eqref{lo-piece-eq} as the energy coming from the low modes of
$u(t_{evac})$, while the $\eta_4^{-3/2} N^{3/2}$ term comes from the
nonlinear corrections generated by the high modes of $u(t)$ for
$t_{min} \leq t \leq t_{evac}$. Note the very strong decay of
$N^{3/2}$ as $N \to 0$; this means that the high modes cannot project
their energy very far into the low modes.  This estimate should be compared with
Proposition \ref{lowfreq-estimate}.  This estimate begins to deteriorate if $N$ gets
too close to $\eta_4$; we avoid this problem by making the high-low frequency decomposition
$u = u_{hi} + u_{lo}$ at $\eta_4^{100}$ instead of $\eta_4$.  Note this bound is
superior to either \eqref{vhf-bounds} or \eqref{boundedl4} at low frequencies.

\begin{proof}
As usual, we rely on the
continuity method, although now we will evolve\footnote{The arguments in this section
seem to rely in an essential way on evolving both forwards and backwards in time
simultaneously; compare with the ``double Duhamel trick'' in Section \ref{flim-conclude}.}
\emph{backwards} in time from
$t_{evac}$, rather than forwards from $t_{min}$.
Let $C_0$ be a large absolute constant (not depending on any of the
$\eta_j$) to be chosen later.  Let $\Omega \subseteq [t_{min},t_{evac})$
denote the set of all times $t_{min} \leq t < t_{evac}$ such that we
have the bounds 
\begin{equation}\label{t2}
\| P_{\leq N} u \|_{\dot S^1([t,t_{evac}] \times \R^3)} \leq C_0
\eta_5 + \eta_0 \eta_4^{-3/2} N^{3/2}
\end{equation}
for all $N \leq \eta_4$.  To show \eqref{lo-piece-eq}, it will
clearly suffice to show that $t_{min} \in \Omega$ (the additional factor
of $\eta_0$ is useful for the continuity method but will be discarded
for the final estimate \eqref{lo-piece-eq}).  In particular, we observe
from \eqref{t2} that we have
\begin{equation}\label{t2-crude}
\| P_{\leq N} u \|_{\dot S^1([t,t_{evac}] \times \R^3)} \lesssim \eta_0
\end{equation}
for all $N \leq \eta_4$.


First we show that $t \in \Omega$ for $t$ sufficiently close to
$t_{evac}$.  We use \eqref{strichartz-def}, H\"older and Sobolev to
estimate
\bas
\| P_{\leq N} u \|_{\dot S^1([t,t_{evac}] \times \R^3)} &\lesssim \|
\nabla P_{\leq N} u \|_{L^\infty_t L^2_x([t,t_{evac}] \times \R^3)} +
 \| \nabla u \|_{L^2_t L^6_x([t,t_{evac}] \times \R^3)} \\
&\lesssim \| \nabla P_{\leq N} u (t_{evac}) \|_{L^2} + C |t_{evac}-t| \|
\nabla \partial_t u \|_{L^\infty_t L^2_x(I_0 \times \R^3)} \\
&\quad +  |t_{evac}-t|^{1/2} \| \nabla u \|_{L^\infty_t L^6_x(I_0
\times \R^3)}.
\end{align*}
Since $u$ is Schwartz, the latter two norms are finite (though
perhaps very large).  By \eqref{T-medium-verysmall} we thus have the
estimate
$$
\| P_{\leq N} u \|_{\dot S^1([t,t_{evac}] \times \R^3)} \lesssim \eta_5 +
C(I_0,u) |t_{evac}-t| + C(I_0,u) |t_{evac}-t|^{1/2}.$$
We thus see that the bound \eqref{t2} holds for $t$ sufficiently
close to $t_{evac}$, if $C_0$ is chosen large enough (but not depending
on $I_0$, $u$ or any of the $\eta_j$.)

Now suppose that $t \in \Omega$, so that \eqref{t2} holds for all $N \leq \eta_4$.  We shall
bootstrap \eqref{t2} to
\begin{equation}\label{t2-boot}
\| P_{\leq N} u \|_{\dot S^1([t,t_{evac}] \times \R^3)} \leq
\frac{1}{2} C_0 \eta_5 + \frac{1}{2} \eta_0 \eta_4^{-3/2} N^{3/2}
\end{equation}
for all $N \leq \eta_4$.
If this claim is true, this would imply (since $u$ is Schwartz) that
$\Omega$ is both open and closed, and so we have $t_{min} \in \Omega$ as desired.

It remains to deduce \eqref{t2-boot} from \eqref{t2}.  For the rest
of the proof, all spacetime norms will be restricted to $[t,
t_{evac}] \times \R^3$.

Fix $N \leq \eta_4$.  By \eqref{nls} and Lemma \ref{disjointed-strichartz}
we have
$$
\| P_{\leq N} u \|_{\dot S^1([t,t_{evac}] \times \R^3)} \lesssim
\|P_{\leq N} u(t_{evac})\|_{\dot H^1(\R^3)}
+ C \sum_{m=1}^M \| \nabla F_m \|_{L^{q'_m}_t L^{r'_m}_x([t,t_{evac}] \times \R^3)}$$
for some dual $L^2$-admissible exponents $(q'_m,r'_m)$, and some decomposition
$P_{\leq N}(|u|^4 u) = \sum_{m=1}^M F_m$ which we will give shortly.

From \eqref{T-medium-verysmall} we have
$$  \|P_{\leq N} u(t_{evac})\|_{\dot H^1(\R^3)} \lesssim \eta_5,$$
which is acceptable for \eqref{t2-boot} if $C_0$ is large enough.

Now consider the nonlinear term  $P_{\leq N} (|u|^4 u)$.  We split
$u$ into high and low frequencies $u = u_{\leq \eta_4} + u_{>\eta_4}$,
where of course $u_{\leq \eta_4} := P_{\leq \eta_4} u$ and $u_{>\eta_4} := P_{>\eta_4} u$,
and use \eqref{5-schematic} to decompose
$$ P_{\leq N} (|u|^4 u) = \sum_{j=0}^5 F_j,$$
where $F_j := P_{\leq N} \O(u_{>\eta_4}^j u_{\leq \eta_4}^{5-j})$.
We now treat the various terms separately.

\divider{Case 1. Estimation of $F_2, F_3, F_4, F_5$.}

We estimate these terms in
$L^2_t L^{6/5}_x$ norm.  We use Bernstein's inequality \eqref{sobolev} to bound these terms by
$$ C N^{3/2} \| \O(u_{>\eta_4}^j u_{\leq \eta_4}^{5-j}) \|_{L^2_t L^1_x([t,t_{evac}] \times \R^3)},$$
which by H\"older can be estimated by
$$ C N^{3/2} \| u_{>\eta_4} \|_{L^\infty_t L^6_x([t,t_{evac}] \times \R^3)}^{j-2}
\| u_{\leq \eta_4} \|_{L^\infty_t L^6_x([t,t_{evac}] \times \R^3)}^{5-j}
 \| u_{> \eta_4} \|_{L^4_t
L^4_x([t,t_{evac}] \times \R^3)}^2.$$
Applying \eqref{l6} and Sobolev, as well as \eqref{boundedl4}, we
obtain a bound of
$$ C \eta_1^\half N^{3/2} \eta_4^{-3/2},$$
which is acceptable for \eqref{t2-boot} if $\eta_1$ is sufficiently small. It is at this step that the small
constant $\eta_1$ appearing in \ref{pm} is used to close the bootstrap.

\divider{Case 2a. Estimation of $F_1$ when $N \ll \eta_4$.}

Now consider $F_1$ term
\begin{equation}\label{scheme}
\| \nabla P_{\leq N} \O( u_{>\eta_4} u_{\leq \eta_4}^4 )
\|_{L^{q'_1}_t L^{r'_1}_x}.
\end{equation}
Suppose first that $N < c \eta_4$.  Then this expression vanishes
unless at least one of the four $u_{\leq \eta_4}$ factors has frequency $> c \eta_4$.  Thus
we can essentially write \eqref{scheme} as\footnote{Strictly speaking, one can only write
\eqref{scheme} as a sum of such terms, where some of the $u_{\leq \eta_4}$ factors in $u_{\leq \eta_4}^3$ must be replaced
by either $P_{> c\eta_4} u_{\leq \eta_4}$ or $P_{\leq c\eta_4} u_{\leq \eta_4}$.  As these projections are bounded on every space
under consideration we ignore this technicality.  Similarly for other such decompositions in this Lemma.}
$$
\| \nabla P_{\leq N} \O( u_{\leq \eta_4}^3 (P_{> c\eta_4} u_{\leq \eta_4})
u_{> \eta_4} ) \|_{L^2_t L^{6/5}_x([t,t_{evac}] \times \R^3)},
$$
where we have chosen $(q'_1,r'_1) = (2, 6/5)$.
By \eqref{boundedl4}, the function $P_{\geq c \eta_4} u_{\leq \eta_4}$ obeys essentially the same
$L^4_{t,x}$ estimates as $u_{> \eta_4}$, and so this term is acceptable by
repeating the arguments used to deal with $F_2,F_3,F_4,F_5$.

\divider{Case 2b. Estimation of $F_1$ when $N \sim \eta_4$.}

Now consider the case when $N \geq c \eta_4$.  In this case we choose
$(q'_1,r'_1) = (1,2)$ and use \eqref{nabla-block}, to bound
\eqref{scheme} by
$$ C \eta_4 \| \O( u_{> \eta_4} u_{\leq \eta_4}^4 ) \|_{L^1_t L^2_x([t,t_{evac}] \times \R^3)}$$
which we estimate using H\"older by
$$ C \eta_4 \| u_{\leq \eta_4} \|_{L^4_t L^\infty_x([t,t_{evac}] \times \R^3)}^4 \| u_{> \eta_4}
\|_{L^\infty_t L^2_x([t,t_{evac}] \times \R^3)}. $$
From \eqref{h1} we have
$$ \| u_{> \eta_4} \|_{L^\infty_t L^2_x([t,t_{evac}] \times \R^3)} \leq C \eta_4^{-1},$$
while from \eqref{strichartz-components} and \eqref{t2-crude} we have
$$ \| u_{\leq \eta_4} \|_{L^4_t L^\infty_x([t,t_{evac}] \times \R^3)}
\lesssim \| u_{\leq \eta_4} \|_{\dot S^1([t,t_{evac}] \times \R^3)}
\lesssim \eta_0.$$
Thus we can bound \eqref{scheme} by $O(\eta_0^4)$, which is
acceptable for \eqref{t2-boot} since $N \sim \eta_4$.  This concludes the estimation of $F_1$.

\divider{Case 3. Estimation of $F_0$.}

It remains to consider the $F_0$ term; we set $(q'_0,r'_0) = (1,2)$ and estimate
$$ \| \nabla P_{\leq N}\O( u_{\leq \eta_4}^5 )
\|_{L^1_t L^2_x([t,t_{evac}] \times \R^3)}.$$
We split $u_{\leq \eta_4} = u_{< \eta_5} + u_{\eta_5 \leq \cdot \leq \eta_4}$, and consider
any term which involves the very low frequencies $u_{<\eta_5}$; schematically, this is
$$ \| \nabla P_{\leq N}\O( u_{\leq \eta_4}^4 u_{<\eta_5} )
\|_{L^{1}_t L^{2}_x([t,t_{evac}] \times \R^3)}.$$
For this case we discard the $P_{\leq N}$, and apply Lemma \ref{leibnitz-holder}
to estimate this term by
$$
\lesssim \| u_{\leq \eta_4} \|_{\dot S^1([t,t_{evac}] \times \R^3)}^4
\| u_{< \eta_5} \|_{\dot S^1([t,t_{evac}] \times \R^3)}.$$
Applying \eqref{t2} we can bound this by
$$ \lesssim (C_0 \eta_5 + \eta_0)^4 (C_0 \eta_5 + \eta_0 \eta_4^{-3/2} \eta_5^{3/2}) \lesssim C_0 \eta_0^4 \eta_5$$
which is acceptable.

We can thus discard all the terms involving $u_{< \eta_5}$, and reduce to estimating
\begin{equation}\label{qr}
 \| \nabla P_{\leq N}\O( u_{\eta_5 \leq \cdot \leq \eta_4}^5 )
\|_{L^1_t L^2_x([t,t_{evac}] \times \R^3)}.
\end{equation}
By Bernstein \eqref{sobolev}, \eqref{nabla-lowfreq} we may
estimate this by
\begin{align*}
\lesssim N^{3/2} \| \O( u_{\eta_5 \leq \cdot \leq \eta_4}^5 ) \|_{L^1_t L^{3/2}_x([t,t_{evac}] \times \R^3)}
&= N^{3/2} \| u_{\eta_5 \leq \cdot \leq \eta_4} \|_{L^5_t L^{15/2}_x([t,t_{evac}] \times \R^3)}^5.
\end{align*}
But from \eqref{strichartz-components}, Bernstein \eqref{bernstein}, \eqref{nabla-block}, and \eqref{t2} we have
\bas
\| u_{\eta_5 \leq \cdot \leq  \eta_4} \|_{L^5_t L^{15/2}_x} &
\leq \sum_{\eta_5 \leq N' \leq \eta_4} \|P_{N'} u\|_{L^5_t L^{15/2}_x([t,t_{evac}] \times \R^3)} \\
&\lesssim \sum_{\eta_5 \leq N' \leq \eta_4}
(N')^{-3/10} \|\nabla P_{N'} u\|_{L^5_t L^{30/11}_x([t,t_{evac}] \times \R^3)} \\
&\lesssim \sum_{\eta_5 \leq N' \leq \eta_4}
(N')^{-3/10} \| P_{N'} u\|_{\dot S^1([t,t_{evac}] \times \R^3)} \\
&\lesssim \sum_{\eta_5 \leq N' \leq \eta_4} (N')^{-3/10}
(C_0 \eta_5 + \eta_0 \eta_4^{-3/2} (N')^{3/2}) \\
&\lesssim \eta_0 \eta_4^{-3/10},
\end{align*}
so that \eqref{qr} is estimated by $O(\eta_0^5 \eta_4^{-3/2} N^{3/2})$ which is
acceptable if $\eta_0$ is small
enough.  This proves \eqref{t2-boot}, which closes the bootstrap.
\end{proof}

\subsection{Controlling the localized $L^2$ mass increment}

Now we have enough estimates to prove \eqref{l2-increment}.  We first
rewrite
$$ P_{hi}( |u|^4 u ) - |u_{hi}|^4 u_{hi} = P_{hi}( |u|^4 u -
|u_{hi}|^4 u_{hi} - |u_{lo}|^4 u_{lo}) {\bf{+}} P_{hi}(|u_{lo}|^4
u_{lo}) - P_{lo} (|u_{hi}|^4 u_{hi}),$$
and so it will suffice to consider the three quantities
\begin{align}
\int_{t_{min}}^{t_*} &|\int \overline{u_{hi}} P_{hi}(|u|^4 u -
|u_{hi}|^4 u_{hi} - |u_{lo}|^4 u_{lo})\ dx| dt \label{q1}\\
\int_{t_{min}}^{t_*} &|\int \overline{u_{hi}} P_{hi}(|u_{lo}|^4
u_{lo})\ dx| dt \label{q2}\\
\int_{t_{min}}^{t_*} &|\int \overline{u_{hi}} P_{lo}(|u_{hi}|^4
u_{hi})\ dx| dt \label{q3}
\end{align}

We now estimate the three quantities \eqref{q1}, \eqref{q2},
\eqref{q3}.


\divider{Case 1. Estimation of \eqref{q1}.}

We move the self-adjoint operator $P_{hi}$ onto $\overline{u_{hi}}$,
and then apply the pointwise estimate (cf. \eqref{5-schematic})
$$ ||u|^4 u - |u_{hi}|^4 u_{hi} - |u_{lo}|^4 u_{lo}|
\lesssim |u_{hi}|^4 |u_{lo}| + |u_{hi}| |u_{lo}|^4$$
to bound \eqref{q1} by
$$
\eqref{q1} \lesssim \int_{t_{min}}^{t_*} \int |P_{hi} u_{hi}| (|u_{hi}|^4 |u_{lo}|
+ |u_{hi}| |u_{lo}|^4)\ dx dt.
$$
For notational convenience, we will ignore the $P_{hi}$ projection and write $P_{hi} u_{hi}$ as $u_{hi}$;
strictly speaking this is not quite accurate but as $P_{hi} u_{hi}$ obeys all the same estimates
as $u_{hi}$, and we have already placed absolute values everywhere, this is a harmless modification.
We can now write our bound for \eqref{q1} as
\begin{equation}\label{q1-bound}
\eqref{q1} \lesssim  \int_{t_{min}}^{t_*} \int |u_{hi}|^5 |u_{lo}|
+ |u_{hi}|^2 |u_{lo}|^4\ dx dt.
\end{equation}

\divider{Case 1a.  Contribution of $|u_{hi}|^2 |u_{lo}|^4$.}

Consider first the contribution of $|u_{hi}|^2 |u_{lo}|^4$; we have to show that
\begin{equation}
  \label{twofour}
  \int_{t_{min}}^{t_*} \int |u_{hi}|^2 |u_{lo}|^4 dx dt \ll \eta_1^2.
\end{equation}
 By
H\"older we can bound this contribution by
$$ \eqref{twofour} \lesssim \| u_{hi} \|_{L^\infty_t L^2_x([t_{min},t_*] \times \R^3)}^2 \| u_{lo} \|_{L^4_t
L^\infty_x(  [t_{min},t_*] \times \R^3)}^4.$$
From \eqref{lo-piece-eq} we have
$$ \| u_{lo} \|_{L^4_t L^\infty_x([t_{min},t_*] \times \R^3)}\lesssim \eta_4^{3/2}$$
while from \eqref{h1}, \eqref{lo-piece-eq} we have
\bas
\| u_{hi} \|_{L^\infty_t L^2_x([t_{min},t_*] \times \R^3)} &\leq \| P_{\geq \eta_4}u_{hi} \|_{L^\infty_t
L^2_x([t_{min},t_*] \times \R^3)} + \sum_{\eta_4^{100} \leq N \leq\eta_4} \| P_N u_{hi}
\|_{L^\infty_t L^2_x([t_{min},t_*] \times \R^3)} \\
&\lesssim \eta_4^{-1} \| u \|_{L^\infty_t \dot H^1_x([t_{min},t_*] \times \R^3)} +
\sum_{\eta_4^{100} \leq N \leq C \eta_4} N^{-1} \| P_N \nabla u_{hi}
\|_{L^\infty_t L^2_x([t_{min},t_*] \times \R^3)} \\
&\lesssim \eta_4^{-1} + \sum_{\eta_4^{100} \leq N \leq \eta_4} N^{-1} \|
P_N u_{hi} \|_{\dot S^1([t_{min},t_*] \times \R^3)} \\
&\lesssim \eta_4^{-1} + \sum_{\eta_4^{100} \leq N \leq \eta_4} N^{-1}
(\eta_5 + \eta_4^{-3/2} N^{3/2})\\
&\lesssim \eta_4^{-1} +  \eta_5 \eta_4^{-100} +  \eta_4^{-1} \\
&\lesssim \eta_4^{-1}.
\end{align*}
We thus obtain a bound of $O(\eta_4^4)$, which is acceptable.

\divider{Case 1b.  Contribution of $|u_{hi}|^5 |u_{lo}|$.}

It remains to control the contribution of $|u_{hi}|^5 |u_{lo}|$; in other words, we need to show
\begin{equation}
  \label{fiveone}
  \int_{t_{min}}^{t_*} \int |u_{hi}|^5 |u_{lo}| dx dt \ll \eta_1^2.
\end{equation}
This estimate will also be useful in controlling \eqref{q3}.

We will split $u_{lo}$ further into somewhat low, and very low,
frequency pieces:
$$ u_{lo} = P_{> \eta_5^{1/2}} u_{lo} + P_{\leq \eta_5^{1/2}} u.$$
The contribution of the very low frequency piece to \eqref{fiveone} can be bounded by Sobolev
embedding \eqref{sobolev} and \eqref{l6} as
\begin{align*}
  \int_{t_{min}}^{t_*} \int |u_{hi}|^5 |u_{\leq \eta_5^{1/2}}| dx dt
&\lesssim \| u_{hi} \|_{L^5_t L^5_x([t_{min},t_*] \times \R^3)}^5 \| u_{\leq \eta_5^{1/2}} \|_{L^\infty_t L^\infty_x([t_{min},t_*] \times \R^3)} \\
&\lesssim C(\eta_4) \| \nabla u \|_{L^5_t L^{30/11}_x([t_{min},t_*] \times \R^3)}^5
\eta_5^{1/4} \| \nabla u \|_{L^\infty_t L^2_x([t_{min},t_*] \times \R^3)} \\
&\lesssim C(\eta_4) \eta_5^{1/4} \| u \|_{\dot S^1([t_{min},t_*] \times \R^3)}^5,
\end{align*}
which is acceptable by \eqref{vhf-bounds}.  Hence we only need to consider the somewhat
low frequencies $P_{> \eta_5^{1/2}} u_{lo}$.
By H\"older, we obtain the bound
\begin{equation}\label{low-vlow}
\||u_{hi}|^{5}|P_{> \eta_5^{1/2}} u_{lo}|\|_{L^{1}_{t,x}([t_{min},t_*] \times \R^3)}\leq C \| u_{hi}
\|_{L^{10}_t L^5_x([t_{min},t_*] \times \R^3)}^5 \| P_{> \eta_5^{1/2}} u_{lo}
\|_{L^2_t L^\infty_x([t_{min},t_*] \times \R^3)}.
\end{equation}
From Bernstein \eqref{bernstein} and \eqref{lo-piece-eq} we have
\bas
\| P_{> \eta_5^{1/2}} u_{lo} \|_{L^2_t L^\infty_x([t_{min},t_*] \times \R^3)}
&\leq \sum_{\eta_5^{1/2} < N \leq \eta_4^{100}} \| P_N u \|_{L^2_t
L^\infty_x([t_{min},t_*] \times \R^3)} \\
&\leq \sum_{\eta_5^{1/2} < N \leq \eta_4^{100}} N^{-1/2} \| \nabla P_N u
\|_{L^2_t L^6_x([t_{min},t_*] \times \R^3)} \\
&\leq \sum_{\eta_5^{1/2} < N \leq \eta_4^{100}} N^{-1/2} \| P_N u \|_{\dot
S^1([t_{min},t_*] \times \R^3)} \\
&\leq \sum_{\eta_5^{1/2} < N \leq \eta_4^{100}} N^{-1/2} ( \eta_5 +
\eta_4^{-3/2} N^{3/2} ) \\
&\leq C \eta_4^{-3/2} \eta_4^{100}.
\end{align*}
To estimate $\|u_{hi} \|_{L^{10}_t L^5_x([t_{min},t_*] \times \R^3)}$, we split $u_{hi}$ into
the higher frequencies $u_{>\eta_4}$ and the medium frequencies $u_{\eta_4^{100} \leq \cdot \leq \eta_4}$.
For the higher frequencies we use\footnote{Note that this application of \eqref{boundedl4} does not require the small constant $\eta_1$. 
} \eqref{boundedl4}, \eqref{l6}, and H\"older to obtain
 \begin{align*}
 \| u_{> \eta_4} \|_{L^{10}_t L^5_x([t_{min},t_*] \times \R^3)} &\lesssim
 \| u_{> \eta_4} \|_{L^4_t L^4_x([t_{min},t_*] \times \R^3)}^{2/5} \| u_{> \eta_4} \|_{L^\infty_t
 L^6_x([t_{min},t_*] \times \R^3)}^{3/5} \\
 &\lesssim \eta_4^{-3/10},  \\
 \end{align*}
while for the medium frequencies we instead use Bernstein \eqref{bernstein}, \eqref{nabla-block},
\eqref{strichartz-components} and Lemma \ref{l2l1-bound} to estimate
\begin{align*}
\| u_{\eta_4^{100} \leq \cdot \leq \eta_4} \|_{L^{10}_t L^5_x([t_{min},t_*] \times \R^3)}
&\lesssim \sum_{\eta_4^{100} \leq N \leq \eta_4} \| u_N \|_{L^{10}_t L^5_x([t_{min},t_*] \times \R^3)} \\
&\lesssim \sum_{\eta_4^{100} \leq N \leq \eta_4} N^{-3/10} \| \nabla u_N \|_{L^{10}_t L^{30/13}_x([t_{min},t_*] \times \R^3)}\\
&\lesssim \sum_{\eta_4^{100} \leq N \leq \eta_4} N^{-3/10} \| u_N \|_{\dot S^1([t_{min},t_*] \times \R^3)} \\
&\lesssim \sum_{\eta_4^{100} \leq N \leq \eta_4} N^{-3/10} ( \eta_5 + \eta_4^{-3/2} N^{3/2} ) \\
&\lesssim \eta_4^{-3/10}.
\end{align*}
Inserting these bounds into \eqref{low-vlow} we obtain a bound of
$\eta_4^{-3} \eta_4^{100}$ for \eqref{fiveone}, which is acceptable.

\divider{Case 2. Estimation of \eqref{q2}.}

Because of the presence of the $P_{hi}$ projection, one of
the $u_{lo}$ terms must have frequency $\geq c \eta_4^{100}$.
We then move $P_{hi}$ over to the $u_{hi}$, bounding
\eqref{q2} as a sum of terms which are essentially of the form\footnote{Actually,
some of the $u_{lo}$ factors in $|u_{lo}|^4$ may have to be replaced by either
$P_{\geq c \eta_4^{100}} u_{lo}$ or $P_{< c \eta_4^{100}} u_{lo}$, but this will make
no difference to the estimates.}
$$
\int_{t_{min}}^{t_*} \int |P_{hi} u_{hi}| |P_{\geq c\eta_4^{100}} u_{lo}|
|u_{lo}|^4\ dx dt.$$
Now observe that $P_{hi} u_{hi}$ and  $P_{\geq c\eta_4^{100}} u_{lo} = P_{lo} u_{\geq c\eta_4^{100}}$
satisfy essentially the same estimates as $u_{hi}$, so this expression can be shown
to be acceptable by a minor modification of \eqref{twofour}.

\divider{Case 3. Estimation of \eqref{q3}.}

We move projections around, using the identity $P_{lo} u_{hi} = P_{hi}
u_{lo}$, to write \eqref{q3} as
\begin{equation}
\label{recastq3}
\int_{t_{min}}^{t_*} |\int \overline{P_{hi} u_{lo}} |u_{hi}|^4
u_{hi} dx| dt.
\end{equation}
Thus, we are concerned here with a term involving five $u_{hi}$ factors and one $P_{hi} u_{lo}$ factor.  But this
is basically \eqref{fiveone}, which has already been shown to be acceptable.  (We have $P_{hi} u_{lo}$
instead of $u_{lo}$ but the reader may verify that the $P_{hi}$ is harmless since it does not
destroy any of the estimates of $u_{lo}$).

This proves \eqref{l2-increment}, and the proof of Proposition
\ref{energy-travel} is complete.  This (finally!) concludes the proof of Theorem \ref{main}.

\endprf

\section{Remarks}\label{remarks-sec}

We make here some miscellaneous remarks concerning certain variants of Theorem \ref{main}.

The global well-posedness result in Theorem \ref{main} was asserted with regard to finite energy solutions $u$
in the class $C^0_t \dot H^1_x \cap L^{10}_{t,x}$, in that the solution existed and was unique in this class,
and depended continuously on the initial data (cf. Lemma \ref{perturb}).  However, the uniqueness result can
be strengthened, in the sense that the solution constructed by Theorem \ref{main} is in fact the only such
solution in the class $C^0_t \dot H^1_x$ (without the assumption of finite $L^{10}_{t,x}$ norm).  This type of
``unconditional well-posedness'' result was first obtained in \cite{kato}, \cite{katounique} (see also \cite{twounique},
\cite{FPT_NLSunique}); the result in \cite{kato}, \cite{katounique} was phrased for the sub-critical Schr\"odinger equation but can
be extended to the critical setting thanks to the endpoint Strichartz estimates in \cite{tao:keel} (or
Lemma \ref{disjointed-strichartz}).  For the convenience of the reader we sketch the ideas of this argument, which
are essentially in \cite{katounique}, \cite{twounique}, \cite{FPT_NLSunique}) here; we are indebted to
Thierry Cazenave for pointing out the relevance of the endpoint Strichartz estimate to the $\dot H^1$-critical uniqueness problem.

Let $u_0$ be finite energy initial data, and let $u \in C^0_t \dot H^1_x \cap L^{10}_{t,x}$ be the (global)
solution to \eqref{nls} constructed in Theorem \ref{main} with this initial data, thus $u(0) = u_0$.   Suppose for
contradiction that we have another (local or global) solution $v \in C^0_t \dot H^1_x$ to \eqref{nls} with initial $u_0$,
in the sense that $v$ verifies the (Duhamel) integral formulation of \eqref{nls},
$$ v(t) = e^{it\Delta} u_0 - i \int_0^t e^{i(t-s)\Delta}(|v|^4 v(s))\ ds.$$
Note that $v \in C^0_t \dot H^1_x \subseteq C^0_t L^6_x$ by Sobolev embedding, so in particular the non-linearity $|v|^4 v$
is locally integrable, and the right-hand side of the above formula makes sense distributionally at least.  We now
claim that $u \equiv v$ on the entire time interval for which $v$ is defined.  Actually, we shall just show that
$u \equiv v$ for all times $t$ in a sufficiently small neighborhood $I$ of 0; one can then extend this to the whole time
interval by a continuity argument and time translation invariance.

To prove the claim, we write $v = u+w$ and observe that $w$ obeys a difference equation, which we write in
integral form as
$$ w(t) = -i \int_0^t e^{i(t-s)\Delta}(|u+w|^4(u+w)(s) - |u|^4 u(s))\ ds.$$
Let $\eps > 0$ be a small number to be chosen shortly.
Note that $w \in C^0_t \dot H^1_x \subseteq C^0_t L^6_x$ and $w(0) = 0$, so in particular we can ensure that
$\| w \|_{L^\infty_t L^6_x(I \times \R^3)} \leq \eps$ by choosing $I$ sufficiently small.  Also, from the Strichartz analysis
$u$ has finite $\dot S^1$ norm, and in particular it has finite $L^8_t L^{12}_x$ norm. Thus 
we can also ensure that $\| u \|_{L^8_t L^{12}_x(I \times \R^3)} \leq \eps$ by choosing $I$ sufficiently small.
Now we use \eqref{5-schematic} to write the equation for $w$ as
$$ w(t) = \int_0^t e^{i(t-s)\Delta}(O(|w(s)|^5) + O(|u(s)|^4 |w(s)|))\ ds.$$
We apply Lemma \ref{disjointed-strichartz} with $k=0$ to conclude in particular that
$$ \| w \|_{L^2_t L^6_x(I \times \R^3)}
\leq C \| |w|^5 \|_{L^2_t L^{6/5}_x(I \times \R^3)} + C \| |u|^4 |w| \|_{L^1_t L^2_x(I \times \R^3)}.$$
From our choice of $I$ and H\"older's inequality we see in particular that
$$  \| w \|_{L^2_t L^6_x(I \times \R^3)} \leq C \eps^4 \| w \|_{L^2_t L^6_x(I \times \R^3)}.$$
Note that the $L^2_t L^6_x(I \times \R^3)$ norm of $w$ is finite since $w \in C^0_t L^6_x$.  If we choose $\eps$
sufficiently small, we then conclude that $w$ vanishes identically on $I \times \R^3$.  One can then extend this vanishing to
the entire time interval for which $v$ is defined by a standard continuity argument which we omit.

We now briefly discuss possible extensions to Theorem \ref{main}.  One obvious extension to study would be
the natural analogue of Theorem \ref{main} in higher dimensions $n > 3$, with the equation \eqref{nls} replaced by
its higher-dimensional energy-critical counterpart
$$ iu_t + \Delta u = |u|^{\frac{4}{n-2}} u.$$
The four-dimensional case $n=4$ seems particularly tractable since the non-linearity is cubic\footnote{Note added in proof: The four-dimensional case has been handled by a very recent preprint of Ryckman and Visan \cite{rv}, using a modification of the methods here.  The case of dimensions five and higher has also been very recently settled (Visan, personal communication).}.
In higher dimensions $n \geq 5$ one no longer expects a regularity result since the non-linearity is not smooth
when $u$ vanishes, however one might still hope for a global well-posedness result in the energy space (especially since
this is already known to be true for small energies, see \cite{cwI}).
In the radial case, such a result was obtained in four dimensions in \cite{borg:scatter}, \cite{borg:book} and
more recently in general dimension in \cite{tao}, and so it is reasonable to conjecture that one in fact has
global well-posedness in the energy space for all dimensions $n \geq 3$ and all finite energy data,
in analogy with Theorem \ref{main}.  However, extending our arguments here to
the higher dimensional setting is far from automatic, even in the four-dimensional case; all the Strichartz
numerology changes, of course, but
also the interaction Morawetz inequality behaves in a somewhat different manner in higher dimensions (since the
quantity $\Delta \frac{1}{|x|}$ is no longer a Dirac mass, but instead a fractional integral potential).
However, it seems that other parts of the argument, such as the induction on energy machinery, the localization of
minimal-energy blowup solutions, and the energy evacuation arguments based on frequency-localized
approximate mass conservation laws, do have a good chance of extending to this setting.  We will not pursue these matters
in detail here.

Another natural extension would be to add a lower order non-linearity to \eqref{nls}, for instance combining the
pure power quintic non-linearity $|u|^4 u$ with a pure power cubic non-linearity $|u|^2 u$.  Heuristically, we do
not expect such lower order terms to affect the global well-posedness and regularity of the equation (especially if
those terms have the same defocusing sign as the top order term), although they may cause some difficulty in
obtaining a scattering result (especially if one adds a non-linearity of the form $|u|^{p-1} u$
for very low $p$, such as $p \leq 1 + \frac{4}{n} = \frac{7}{3}$ or $p \leq 1 + \frac{2}{n} = \frac{5}{3}$).
However, these lower order terms do create some non-trivial difficulties in our argument, which relies heavily on
scale-invariance.  One may need to add some lower order terms (such as the $L^2$ mass) to the energy $E$, or to
the definition of the quantity $M(E)$, in order to salvage the induction on energy argument in this setting.
Again, we will not pursue these matters here\footnote{Note added in proof: the lower order terms have been successfully treated by Xiaoyi Zhang (personal communication), relying on this result and perturbation theory.}.

As remarked in Remark \ref{huge-size}, our final bound $M(E)$ for the global $L^{10}_{t,x}$ norm of $u$
in terms of the energy $E$ is extremely bad; this is due to our extremely heavy reliance on the induction
on energy hypothesis (Lemma \ref{induction}) in our argument\footnote{For a different approach to the problem which yields an unspecified bound uniform in the energy, under the assumption that one 
is given a global $L^{10}_{x,t}$ bound for the solution which is {\em not} uniform in the energy, 
see \cite{keraani}.}.  We do not expect our bounds to be anywhere close
to best possible.  Indeed, any simplification of this argument would almost surely lead to less use of the
induction hypothesis, and consequently to a better bound on $M(E)$.  For recent progress in this direction in
the radial case (in which no induction hypothesis is used at all, leading to a bound on $M(E)$ which is
merely exponential in $E$), see \cite{tao}.

The global existence and
scattering result obtained here has analogs for the
critical nonlinear Klein-Gordon equation
$- \frac{1}{2c^2}u_{tt} +  \Delta u = -|u|^4 u + \frac{m^2c^2}{2} u$ (see
introduction for references.). As we remarked earlier, there are
some important differences between the methods
employed for the Klein-Gordon equation and those we use here.   In particular,
it  is not at all clear how our arguments might help show that the
space-time bounds
for the nonlinear Klein-Gordon equation are uniform in the non-relativistic limit $c \to \infty$, even though one heuristically expects the non-linear Klein-Gordon equation to converge in some sense to the non-linear
Schr\"odinger equation in this regime with suitable normalizations and assumptions on the data.  One major difficulty in extending
our arguments to the relativistic case is that we have no analogue of the interaction Morawetz inequality
\eqref{morawetz-interaction} (or any localized variants) for the Klein-Gordon equation.  For small energy data, uniform bounds on the solution
are available in the nonrelativistic limit (see remarks in \cite{nakrelII}), but
for general solutions such bounds do not
seem available. (See also \cite{nakrelI} and references therein for further results
on the subcritical problem.)

\end{document}